\documentclass[a4paper,12pt]{article}
\usepackage{amsfonts}
\usepackage{amstext}
\usepackage{amsthm}
\usepackage{amsmath}
\usepackage{amssymb}
\usepackage{latexsym}
\usepackage[latin1]{inputenc}
\newcommand{\bl}{\hfill\rule{2mm}{2mm}}
\newcommand{\R}{\Bbb{R}}

\newcommand{\s}{\Bbb{S}}
\newtheorem{teor}{Theorem}[section]
\newtheorem{propo}{Proposition}[section]

\newtheorem{lema}{Lemma}[section]
\newtheorem{cor}{Corollary}[section]
\newtheorem{rem}{Remark}[section]

\newcommand{\n}{\noindent}
\newcommand{\vs}{\vspace}

\begin{document}

\title{Extremal maps in best constants vector theory\\ {\Large Part I: Duality and Compactness}
}
\author{\textbf{Ezequiel R. Barbosa, Marcos Montenegro \footnote{\textit{E-mail addresses}:
ezequiel@mat.ufmg.br (E. R. Barbosa), montene@mat.ufmg.br (M.
Montenegro)}}\\ {\small\it Departamento de Matem\'{a}tica,
Universidade Federal de Minas Gerais,}\\ {\small\it Caixa Postal
702, 30123-970, Belo Horizonte, MG, Brazil}}


\maketitle

\markboth{abstract}{abstract}
\addcontentsline{toc}{chapter}{abstract}

\hrule \vspace{0,2cm}

\n {\bf Abstract}

We develop a comprehensive study on sharp potential type Riemannian Sobolev inequalities of order $2$ by means of a local geometric Sobolev inequality of same kind and suitable De Giorgi-Nash-Moser estimates. In particular we discuss questions like continuous dependence of optimal constants and existence and compactness of extremal maps. The main obstacle arising in the present setting lies at fairly weak conditions of regularity assumed on potential functions.

\vspace{0.5cm}
\hrule\vspace{0.2cm}

\n {\small {\bf Content}}

\n {\small 1. Introduction, overview and main results..................................................................................................\ 2\\
2. Miscellaneous and central tools...............................................................................................................\ 14\\
3. Continuous dependence of optimal constants..........................................................................................\ 37\\
4. The vector duality...................................................................................................................................\ 47\\
5. Compactness of extremal maps...............................................................................................................\ 48\\
6. Examples and counter-examples.............................................................................................................\ 51\\
Appendix: Concentration, existence and regularity lemmas\\
References}

\vspace{0.5cm}
\hrule\vspace{0.2cm}

\pagebreak

\section{Introduction, overview and main results}

Over the past forty years, a lot of attention has been paid to so called sharp Riemannian Sobolev inequalities. A vast literature has led to the building of a rich theory known as best constants theory which has connected areas as analysis, geometry and topology. Such inequalities have for instance played an important role in geometric analysis, specially in the study of non-collapse of the Ricci flow along time (cf. \cite{Hsu}, \cite{Ye}, \cite{Zhang}), existence and multiplicity of solutions for the Yamabe problem (cf. \cite{Au3}, \cite{HV11}, \cite{S}), and Riemannian isoperimetric inequalities (cf. \cite{D4}). Particularly, a number of results have shown the strong influence exerted by the geometry and/or topology on various questions in this field. We refer to Aubin \cite{Au4}, Druet and Hebey \cite{DH}, Hebey \cite{H3} for surveys in book form and other references therein.

The goal of the present paper is to extend part of the well-known best constants scalar $L^2$-theory to a fairly general vector perspective. The meaning of the word ``general" will soon clarify. Precisely, we wish to

\begin{itemize}

\item[$(\alpha.1)$] establish necessary and sufficient topological conditions to the continuity of optimal constants on all involved parameters: potential functions and Riemannian metrics;

\item[$(\alpha.2)$] exhibit a dichotomy adapted to the present setting: existence of extremal maps versus explicit values of optimal constants;

\item[$(\alpha.3)$] discuss the compactness problem of extremal maps related to families of potential functions and Riemannian metrics.

\end{itemize}

\n A major challenge here is to address $(\alpha.1)$-$(\alpha.3)$ under rather weak assumptions of regularity. As we shall see, such assumptions prevent us for instance to argue with Euler-Lagrange equations which are essential in the development of significant part of the scalar theory.

Before we go further and precise these points, a little bit of notation and background should be introduced.

\n {\bf 1.1 Sharp Euclidean Sobolev inequalities}. Sobolev inequalities are among the
most famous and useful functional inequalities in analysis. They express a strong
integrability or regularity property for a function in terms of some integrability
property for some derivatives of this function. The typical {\bf sharp Euclidean $L^2$-Sobolev inequality} states that, for any $u \in C_0^\infty(\R^n)$,

\begin{equation} \label{E.1}
\left (\int _{\mathbb{R}^n}|u|^{2^*}\; dx\right
)^{2/2^*}\leq A_0(n) \int_{\mathbb{R}^n}|\nabla u|^2\; dx\,,
\end{equation}
where $n \geq 3$, $2^* = \frac{2n}{n-2}$ and, by definition, $A_0(n)$ is the best
possible constant in this inequality. For the great majority of applications,
it is not necessary to know more about the Sobolev embedding, apart maybe from
explicit bounds on $A_0(n)$. However, in some circumstances one is interesting to know the exact
value of $A_0(n)$ as in the computation of the ground state energy in some physical models. Most often, the
determination of $A_0(n)$ relies on the identification of extremal functions to (\ref{E.1}), i.e. non-zero functions satisfying equality in (\ref{E.1}). The natural space to look for extremal functions is the usual Euclidean Sobolev space

\[
{\cal D}^{1,2}(\R^n) := \{u \in L^{2^*}(\R^n):\; |\nabla u| \in L^2(\R^n)\}
\]
endowed with the norm

\[
||u||_{{\cal D}^{1,2}(\R^n)} := \left(\int_{\R^n} |\nabla u|^2\; dx\right)^{1/2}\; .
\]
The corresponding extremal functions to (\ref{E.1}) are given precisely by

\[
u_0(x) = a w(b(x - x_0))
\]
with $a, b \in \mathbb{R}$, $b \not = 0$, $x_0 \in \mathbb{R}^n$ and

\begin{equation}\label{fe1a}
w(x) = \left(1 + |x|^{2}\right)^{- n/2^*} \; .
\end{equation}
This fact is due to Aubin \cite{Au2}, Rodemich (unpublished) and Talenti \cite{Ta}.

A natural extension which is directly linked to estimates of ground state energy related to some reaction-diffusion systems is as follows.

Let $k \geq 1$ be an integer number and $n \geq 3$. Consider a positive continuous function $F: \mathbb{R}^k \rightarrow \mathbb{R}$ homogeneous of degree $2^*$, i.e. $F(\lambda t) = \lambda^{2^*} F(t)$ for $\lambda > 0$ and $t \in \R^k$. Part of the vector notations below are followed from \cite{HH3}. Denote by $C_{0,k}^\infty(\R^n)$ the space $C_0^\infty(\R^n) \times \ldots \times C_0^\infty(\R^n)$ of smooth $k$-maps with compact support in $\R^n$.

The {\bf sharp potential type Euclidean $L^2$-Sobolev inequality} states that, for any $U \in C_{0,k}^\infty(\R^n)$,

\begin{equation} \label{Euc}
\left(\int_{\mathbb{R}^n} F(U) \; dx\right)^{2/2^*} \leq \mathcal{A}_0(n,F) \int
_{\mathbb{R}^n}|\nabla U|^2 \; dx\; ,
\end{equation}
where

\[
\int_{\mathbb{R}^n} |\nabla U|^2\; dx := \sum_{i=1}^k \int_{\mathbb{R}^n} |\nabla u_i|^2\; dx\; ,\ \ U = (u_1, \ldots ,u_k)\; ,
\]
and $\mathcal{A}_0(n,F)$ is the best possible constant in this inequality. Obviously, $\mathcal{A}_0(n,F)$ is well defined due to the conditions assumed on $F$. The name given to (\ref{Euc}) comes from physical motivations where $F$ represents a potential function in a specific model. In this case, we have the notion of extremal map as being a non-zero $k$-map satisfying equality in (\ref{Euc}). The adequate space to seek extremal maps is the $k$-vector Euclidean Sobolev space $\mathcal{D}^{1,2}_k(\mathbb{R}^n) := \mathcal{D}^{1,2}(\mathbb{R}^n) \times \ldots \times \mathcal{D}^{1,2}(\mathbb{R}^n)$ equipped with the norm

\[
||\nabla U||_{\mathcal{D}_k^{1,2}(\mathbb{R}^n)} := \left( \int_{\mathbb{R}^n} |\nabla U|^2\; dx
\right)^{1/2}\; .
\]
The set of extremal maps to (\ref{Euc}) and the value of $\mathcal{A}_0(n,F)$ can fortunately be found. This set consists precisely of $k$-maps of the form $U_0 = t_0 u_0$, where $t_0 \in \mathbb{S}^{k-1} := \{t \in \R^k:\; |t| = 1\}$ is a maximum point of $F$ and $u_0 \in \mathcal{D}^{1,2}(\mathbb{R}^n)$ is an extremal function to (\ref{E.1}), and so $\mathcal{A}_0(n,F) = F(t_0)^{2/2^*} A_0(n)$. For completeness, we include in this paper a proof of these claims. It is important to emphasize that this characterization of extremal maps to (\ref{Euc}) is essential to achieve our goals.\\

\n {\bf 1.2 Sharp Riemannian Sobolev inequalities}. Considerable efforts have been spent on the investigation of sharp Riemannian Sobolev inequalities during a few decades. Part of the developments obtained is currently known in the literature as the $AB$ program. Below, we briefly comment this program for $L^2$-Sobolev inequalities.

Let $(M,g)$ be a smooth compact Riemannian manifold of dimension $n \geq 3$ and $\beta : M \rightarrow \R$ be a positive continuous function. By using a standard unity partition argument and (\ref{E.1}), one gets constants $A, B \in \R$ such that, for any $u \in C^\infty(M)$,

\begin{equation} \label{AB}
\left(\int_{M} |u|^{2^*}\; dv_g\right)^{2/2^*} \leq A \int_{M} |\nabla_g u|^2 \; dv_g + B \int_{M} \beta(x) u^2\; dv_g\; ,
\end{equation}
where $dv_g$ and $\nabla_g$ denote, respectively, the Riemannian volume element and the gradient operator of $g$. Some basic notations and definitions related to (\ref{AB}) are now introduced.

The {\bf first Riemannian $L^2$-Sobolev best constant} is defined by

\[
A_0(n,\beta,g) := \inf \{ A \in \R:\; \mbox{ there exists} \hspace{0,18cm} B \in \R \hspace{0,18cm} \mbox{such that (\ref{AB})} \hspace{0,18cm} \mbox{is valid}\}\; .
\]

The {\bf first sharp Riemannain $L^2$-Sobolev inequality} states that there exists a constant $B \in \R$ such that, for any $u \in C^\infty(M)$,

\begin{equation} \label{A-opt}
\left(\int_{M} |u|^{2^*}\; dv_g\right)^{2/2^*} \leq A_0(n,\beta,g) \int_{M} |\nabla_g u|^2 \; dv_g + B \int_{M} \beta(x) u^2\; dv_g\; .
\end{equation}
In this case, one defines the {\bf second Riemannian $L^2$-Sobolev best constant} by

\[
B_0(n,\beta,g) := \inf \{ B \in \R:\; (\ref{A-opt}) \hspace{0,18cm} \mbox{is valid} \}\; ,
\]
and the {\bf second sharp Riemannian $L^2$-Sobolev inequality} as the saturated version of (\ref{A-opt}), i.e.

\begin{equation} \label{B-opt}
\left(\int_{M} |u|^{2^*}\; dv_g\right)^{2/2^*} \leq A_0(n,\beta,g) \int_{M} |\nabla_g u|^2 \; dv_g + B_0(n,\beta,g) \int_{M} \beta(x) u^2\; dv_g\; .
\end{equation}
Note that (\ref{B-opt}) is sharp with respect to both the first and second best constants in the sense that none of them can be lowered. In a natural way, it arises then the notion of extremal functions as being non-zero functions satisfying (\ref{B-opt}) with equality. The appropriate space to look for extremal functions is the Riemannian Sobolev space $H^{1,2}(M)$ defined as the completion of $C^{\infty}(M)$ under the norm

\[
||u||_{H^{1,2}(M)} := \left( \int_{M} |\nabla_g u|^2\; dv_g  + \int_{M} u^2\; dv_g \right)^{1/2}\; .
\]

The scalar $AB$ program in the $L^2$ environment consists of various questions involving the optimal constants $A_0(n,\beta,g)$ and $B_0(n,\beta,g)$ and the sharp Sobolev inequalities (\ref{A-opt}) and (\ref{B-opt}).

Some of them are:

\begin{itemize}

\item[(a)] Is it possible to find the exact values and/or bounds for $A_0(n,\beta,g)$ and $B_0(n,\beta,g)$?

\item[(b)] Is the sharp $L^2$-Sobolev inequality (\ref{A-opt}) valid?

\item[(c)] Do $A_0(n,\beta,g)$ and $B_0(n,\beta,g)$ depend continuously on $\beta$ and $g$ in some topology?

\item[(d)] Does the sharp $L^2$-Sobolev inequality (\ref{A-opt}) admit any extremal function?

\item[(e)] Is the set of extremal functions to (\ref{B-opt}) with unit $L^{2^*}$-norm compact in some topology?

\end{itemize}

\n These questions were initially addressed in the classical case $\beta \equiv 1$ and important complete or partial answers were given during the penultimate decade. In 2001, targeting the study of questions surrounding (d), the scalar $AB$ program began to be focused in the general setting by Hebey and Vaugon \cite{HV2} in connection with the notion of critical function. We highlight below some partial answers known so far.

\n {\it On the question (a)}: It is well known since Aubin \cite{Au2} that $A_0(n,\beta,g) = A_0(n)$. In particular, $A_0(n,\beta,g)$ depends only on the dimension $n$. Unlike, $B_0(n,\beta,g)$ depends also on the variables $\beta$ and $g$ as can be seen from the scaling relation $B_0(n,\lambda \beta, \mu g) = \lambda^{-1} \mu^{-1} B_0(n,\beta,g)$ valid for any constants $\lambda, \mu > 0$. Explicit values of $B_0(n,1,g)$ have been computed in only specific situations. For instance, Aubin \cite{Au2} showed that on the round unit $n$-sphere $\mathbb{S}^n$,

\[
B_0(n,1,g)=\omega_n^{-\frac{2}{n}}\; ,
\]
where $\omega_n$ stands for the volume of $\mathbb{S}^n$. On upper
bounds, Hebey and Vaugon \cite{H3} (see also \cite{HEBEY} and \cite{HV11}) discovered that on the product manifold $\mathbb{S}^1
\times \mathbb{S}^{n-1}$ induced with the usual Euclidean product
metric,

\[
B_0(n,1,g) \leq \frac{1+(n-2)^2}{n(n-2)\omega _n^{2/n}}\;,
\]
and on the canonical Euclidean projective space $\mathbb{P}^n$,
\[
B_0(n,1,g) \leq \frac{n+2}{(n-2)\omega _n^{2/n}}\;.
\]
Unfortunately, the explicit value of $B_0(n, \beta, g)$ is no known in general. However, one knows that, by taking $u \equiv 1$ in (\ref{B-opt}),

\[
B_0(n,\beta,g) \geq \frac{v_g(M)^{2/2^*}}{\int_M \beta(x)\; dv_g}\; .
\]

\n Besides, Hebey and Vaugon \cite{HV2} proved the following geometric lower bound for $n\geq 4$:

\begin{equation} \label{G-E-s}
B_0(n,\beta,g) \beta(x) \geq \frac{n-2}{4(n-1)} A_0(n) S_g(x)
\end{equation}

\n for all $x \in M$, where $S_g$ denotes the scalar curvature of the metric $g$.

\n {\it On the question (b)}: This question was raised by Aubin \cite{Au2} and answered positively by Hebey and Vaugon \cite{HV1} in a pioneer work in this field which introduced important ideas based on concentration analysis of solutions of elliptic PDEs.

\n {\it On the question (c)}: The geometric continuity of $B_0(n,1,g)$ has recently been discussed in \cite{EM6}. Let $M$ be a compact differentiable manifold of dimension $n$. Denote by ${\cal M}^n$ the space of smooth Riemannian metrics on $M$ equipped with the $C^2$-topology. When $n \geq 4$, the authors showed that the map $g \in {\cal M}^n \mapsto B_0(n,1,g)$ is continuous and, moreover, the $C^2$-topology is sharp among all $C^k$-topologies.

\n {\it On the question (d)}: When $(M,g)$ is the round unit $n$-sphere $(\s^n,h)$, the extremal functions to (\ref{B-opt}) are all classified. Precisely, modulo a constant factor and/or an isometry on $\s^n$, they are given by $u \equiv 1$ or $u = (\beta - \cos r)^{-n/2^*}$, where $\beta > 1$ is a constant and $r$ is the distance to some point on $\s^n$. By Hebey \cite{H2}, if $n \geq 4$ and $g$ is a conformal metric to $h$, then

\[
B_0(n,1,g) = \frac{n-2}{4(n-1)} A_0(n) \max \limits_{M} S_g
\]

\n and there exist extremal functions to (\ref{B-opt}) if, and only if, modulo a constant scale factor, $g$ and $h$ are isometric, in which case all the extremal functions are known. For an arbitrary compact Riemannian manifold $(M,g)$ of dimension $n$, Djadli and Druet \cite{DjDr} proved that, at least, one of the assertions holds for $n \geq 4$ and $\beta \equiv 1$:

\begin{itemize}

\item[(i)] (\ref{B-opt}) admits an extremal function;

\item[(ii)] $B_0(n,1,g) = \frac{n-2}{4(n-1)} A_0(n) \max \limits_{M} S_g$.

\end{itemize}

\n Applying this result, one easily checks that an extremal function always exists when either $S_g \leq 0$ or $S_g$ is constant on $M$ (cf. \cite{DjDr}). Examples showing that all possible situations in (i)-(ii) can occur are also known. For instance, both (i) and (ii) hold for the round unit $n$-sphere $\mathbb{S}^n$, (i) fails and (ii) holds for $\mathbb{S}^n$ endowed with a conformal metric non-isometric to the canonical metric, whereas (i) holds and (ii) fails for certain quotients of the round unit $n$-sphere $\mathbb{S}^n$. Closely related questions have also been discussed in the works \cite{Barbosa}, \cite{EM} and \cite{HV2}. In this last one, the duality (i)-(ii) has been extended to arbitrary positive functions $\beta$. Precisely, in their Theorem 3, it is shown that, at least, one of the claims below holds when $n \geq 4$:

\begin{itemize}

\item[(i)*] (\ref{B-opt}) admits an extremal function;

\item[(ii)*] $B_0(n,\beta,g) \beta(x_0) = \frac{n-2}{4(n-1)} A_0(n) S_g(x_0)$ for some $x_0 \in M$.

\end{itemize}

\n {\it On the question (e)}: The compactness problem of extremal functions has been addressed in \cite{EM6}, \cite{DjDr} and \cite{HV2}. Denote by $E(\beta, g)$ the set of extremal functions to (\ref{B-opt}) with unit $L^{2^*}$-norm. In \cite{DjDr}, Djadli and Druet proved that $E(1,g)$ is compact in the $C^0$-topology (or is $C^0$-compact) whenever $n \geq 4$ and

\begin{equation} \label{1}
B_0(n,1,g) > \frac{n-2}{4(n-1)} A_0(n) \max \limits_{M} S_g\; .
\end{equation}
Inspired in their proof, we established in \cite{EM6} the compactness of $E(1,g)$ when the geometry varies. Precisely, let $G \subset {\cal M}^n$ be a subset such that (\ref{1}) holds for all $g \in G$. Then, $\bigcup_{g \in G} E(1,g)$ is $C^0$-compact whenever $n \geq 4$ and $G$ is compact in ${\cal M}^n$. Following similar ideas to those ones of \cite{DjDr}, Hebey and Vaugon \cite{HV2} extended the compactness result to positive functions $\beta$ and Riemannian metrics $g$ satisfying

\[
B_0(n,\beta,g) \beta(x) > \frac{n-2}{4(n-1)} A_0(n) S_g(x)
\]
for all $x \in M$. Note, however, that when equality occurs in (\ref{1}), the compactness can fail in all dimensions as can easily be seen from the round unit $n$-sphere $\mathbb{S}^n$ due to the noncompactness of its conformal group.

The purpose of this work is to investigate sharp Riemannian Sobolev inequalities modeled on vector $L^2$-Sobolev spaces. Despite being linked to some questions of analytical interest, as for example the computation of the ground state energy, Sobolev inequalities on vector spaces have been little explored in the literature. A specially important class of such inequalities is so called potential type Riemannian Sobolev inequalities which naturally extend the above-mentioned Riemannian $L^2$-Sobolev inequalities.

We begin with some basic notations and definitions geared to the vector context. Let $(M,g)$ be a smooth compact Riemannian manifold of dimension $n \geq 3$ and $k \geq 1$ be an integer number. Consider positive continuous functions $F: \mathbb{R}^k \rightarrow \mathbb{R}$ and $G: M \times \mathbb{R}^k \rightarrow \mathbb{R}$ with $F$ homogeneous of degree $2^*$ and $G$ homogeneous of degree $2$ on the second variable. Again for physical reasons, $F$ and $G$ are called potential functions. Denote by $C^\infty_k(M)$ the space $C^\infty(M) \times \ldots \times C^\infty(M)$ of smooth $k$-maps on $M$. Using (\ref{AB}), one easily discovers constants $\mathcal{A}, \mathcal{B} \in \mathbb{R}$ such that, for any $U \in C^\infty_k(M)$,

\begin{equation} \label{AB-v}
\left( \int_M F(U) \; dv_g \right)^{2/2^*} \leq \mathcal{A} \int_M |\nabla_g U|^2 \; dv_g + \mathcal{B} \int_M G(x,U)\; dv_g\;,
\end{equation}
where

\[
\int_M |\nabla_g U|^2\; dv_g := \sum_{i=1}^k \int_M |\nabla_g u_i|^2\; dv_g\; , \ \ U = (u_1, \ldots ,u_k)\; .
\]

In parallel to the notations and definitions corresponding to the scalar $AB$ program, we have following definitions related to (\ref{AB-v}).

The {\bf first Riemannian $L^2$-Sobolev best constant} is defined by

\[
\mathcal{A}_0(n,F,G,g) := \inf \{ \mathcal{A} \in \R:\; \mbox{ there exists} \hspace{0,18cm} \mathcal{B} \in \R \hspace{0,18cm} \mbox{such that (\ref{AB-v})} \hspace{0,18cm} \mbox{is valid}\}\; .
\]

The {\bf first sharp potential type Riemannain $L^2$-Sobolev inequality} states that there exists a constant $\mathcal{B} \in \R$ such that, for any $U \in C^\infty_k(M)$,

\begin{equation} \label{A-opt-v}
\left( \int_M F(U) \; dv_g\right)^{2/2^*} \leq \mathcal{A}_0(n,F,G,g) \int_M |\nabla_g U|^2 \; dv_g + \mathcal{B}
\int_M G(x,U) \; dv_g\; ,
\end{equation}
in which case one defines the {\bf second Riemannian $L^2$-Sobolev best constant} by

\[
\mathcal{B}_0(n,F,G,g) := \inf \{ \mathcal{B} \in \R:\; (\ref{A-opt-v}) \hspace{0,18cm} \mbox{is valid} \}\; ,
\]
and the {\bf second sharp potential type Riemannian $L^2$-Sobolev inequality} as the saturated version of (\ref{A-opt-v}), i.e.

\begin{equation} \label{B-opt-v}
\left( \int_{M} F(U)\; dv_g \right)^{2/2^*} \leq \mathcal{A}_0(n,F,G,g) \int_{M} |\nabla_g U|^2 \; dv_g +
\mathcal{B}_0(n,F,G,g) \int_{M} G(x,U)\; dv_g\; .
\end{equation}

\n In a natural way, an extremal map to (\ref{B-opt-v}) is defined as a non-zero $k$-map that realizes equality in (\ref{B-opt-v}). The adequate space to search extremal maps is the Riemannian $k$-vector Sobolev space $H_k^{1,2}(M) := H^{1,2}(M) \times \ldots \times H^{1,2}(M)$ induced with the norm

\[
||U||_{H_k^{1,2}(M)} := \left( \int_{M} |\nabla_g U|^2\; dv_g +
\int_{M} |U|^2\; dv_g \right)^{1/2},
\]
where

\[
\int_{M} |U|^2\; dv_g := \sum_{i=1}^k \int_{M} u_i^2\; dv_g\; ,\ \ U = (u_1, \ldots ,u_k) \; .
\]
Of course, the above definitions generalize the corresponding scalar ones. Indeed, when $k = 1$ we have, modulo constant factors, $F(t) = |t|^{2^*}$ and $G(x,t) = \beta(x) t^2$. However, for $k \geq 2$, there exist many examples of positive homogeneous functions only continuous as can be seen from

\[
F(t) =
\left\{
\begin{array}{ll}
|t|^{2^*} f(t/|t|), & {\rm if}\ t \neq 0\\
0, & {\rm if}\ t = 0
\end{array}
\right. ,
\]

\n where $f : \mathbb{S}^{k-1} \rightarrow \R$ is a positive $C^0$ function. Some canonical examples of potential functions are:

\[
F(t) = |t|_q^{2^*}\; ,\ \ F(t) = ||t||^{2^*}\; ,
\]

\[
G(x,t) = \beta(x) |t|_q^2\; ,\ \ G(x,t) = \langle A(x) t, t \rangle\; ,\ \ G(x,t) = \beta(x) ||t||^2\;,
\]

\n where $|t|_q := (\sum_{i=1}^k |t_i|^q)^{1/q}$ with $q > 0$, $\beta \in C^0(M)$ is a positive function, $\langle \cdot, \cdot \rangle$ stands for the usual Euclidean inner product, $A(x)$ are positive symmetric $k \times k$ matrices with continuous entries on $M$ and $|| \cdot ||$ denotes an arbitrary norm on $\R^k$.

A successful best constants theory on sharp potential type Riemannian Sobolev inequalities relies on answers to a number of questions related to the optimal constants $\mathcal{A}_0(n,F,G,g)$ and $\mathcal{B}_0(n,F,G,g)$ and to the sharp Sobolev inequalities (\ref{A-opt-v}) and (\ref{B-opt-v}). In the same spirit of the scalar $AB$ program, we inquire:

\begin{itemize}

\item[(A)] Is it possible to find the explicit values and/or bounds for $\mathcal{A}_0(n,F,G,g)$ and $\mathcal{B}_0(n,F,G,g)$?

\item[(B)] Is the sharp potential type $L^2$-Sobolev inequality (\ref{A-opt-v}) valid?

\item[(C)] Do $\mathcal{A}_0(n,F,G,g)$ and $\mathcal{B}_0(n,F,G,g)$ depend continuously on $F$, $G$ and $g$ in some topology?

\item[(D)] Does the sharp potential type $L^2$-Sobolev inequality (\ref{B-opt-v}) admit any extremal map?

\item[(E)] Is the set of $F$-normalized extremal maps to (\ref{B-opt-v}) compact in some topology?  A map $U$ is said to be $F$-normalized if $||F(U)||_{L^1(M)} = 1$.

\end{itemize}

In \cite{HH3}, Hebey addressed (A), (B), (D) and (E) to a vector extension closely related to the scalar situation. Precisely, he considered the case when $F(t) = |t|_{2^*}^{2^*}$ and $G(x,t) = \langle A(x) t, t \rangle$, where $A(x) = (a_{ij}(x))$ represent positive symmetric $k \times k$ matrices with smooth entries on $M$, and proved that $\mathcal{A}_0(n,F,G,g) = A_0(n)$ and (\ref{A-opt-v}) is always valid. When $n \geq 4$, one has

\begin{equation} \label{G-I}
\mathcal{B}_0(n,F,G,g) a_{ii}(x) \geq \frac{n-2}{4(n-1)} A_0(n) S_g(x)
\end{equation}
for all $x \in M$ and $i = 1, \ldots, k$ and also arises the following dichotomy:

\begin{itemize}

\item[(I)] (\ref{B-opt-v}) admits an extremal map;

\item[(II)] $\mathcal{B}_0(2,F,G,g) a_{ii}(x_0) = \frac{n-2}{4(n-1)} A_0(n) S_g(x_0)$ for some $x_0 \in M$ and some $i = 1, \ldots, k$.

\end{itemize}

\n Moreover, when the inequality (\ref{G-I}) is strict for all $x \in M$ and $i = 1, \ldots, k$, Hebey proved the $C^0$ compactness of the set
of extremal maps normalized by the $L^{2^*}$-norm. Recently, Druet, Hebey and V\'{e}tois \cite{DHVe} have studied the bounded stability property, something stronger than compactness, of solutions for a class of potential type elliptic systems. Precisely, their main result ensures the bounded stability of the set of solutions of system

\[
- \Delta_g u_i + \sum_{j=1}^k a_{ij}(x) u_j = |U|^{2^* - 2} u_i \ \ {\rm on}\ \ M,\ \ i=1,\ldots,k\ ,
\]

\n where $\Delta_g$ stands for the Laplace-Beltrami operator associated to the metric $g$, under the following condition:

\[
A(x) < \frac{n-2}{4(n-1)} S_g(x) I_k
\]

\n for all $x \in M$ in the sense of bilinear forms, where $I_k$
denotes the $k \times k$ identity matrix. For the study of analytical stability of the system above, which also leads to compactness, under the opposite inequality

\[
A(x) > \frac{n-2}{4(n-1)} S_g(x) I_k\,,
\]

\n see the Druet and Hebey's paper \cite{DH2009}.

In the general case, through a simple argument for reducing the
scalar case, answers to (A) and (B) can be obtained directly from
the scalar $AB$ program. As easily will be checked,
$\mathcal{A}_0(n,F,G,g) = M_F^{2/2^*}A_0(n)$, (\ref{A-opt-v}) is
always valid and the lower estimate holds for $n \geq 4$:

\begin{equation} \label{G-E}
\mathcal{B}_0(n,F,G,g) G(x,t_0) \geq \frac{n-2}{4(n-1)} M_F^{2/2^*} A_0(n) S_g(x)
\end{equation}
for all $x \in M$ and $t_0 \in \mathbb{S}^{k-1}$ with $F(t_0) = M_F$, where $M_F = \max \limits_{|t|=1} F(t)$. These claims are precisely justified in Section 2. Note that the inequality (\ref{G-E}) simultaneously extends the corresponding scalar (\ref{G-E-s}) and vector (\ref{G-I}). Remark also that $\mathcal{A}_0(n,F,G,g)$ depends only on the dimension $n$ and the potential function $F$, while $\mathcal{B}_0(n,F,G,g)$ depends further on the potential function $G$ and the metric $g$ as shows the scaling relation $\mathcal{B}_0(n,\theta F, \lambda G,\mu g) = \theta^{2/2^*} \lambda^{-1} \mu^{-1} \mathcal{B}_0(n,F,G,g)$ valid for any constants $\theta, \lambda, \mu > 0$. It is important, however, to emphasize that the remaining questions are far more delicate than (A) and (B) for several reasons. The first one concerns the weak regularity condition satisfied by the involved potential functions. Indeed, as seen previously, while $F(t)$ and $G(x,t)$ are always of $C^1$ class on $t$ when $k = 1$, the most of these functions are only continuous when $k \geq 2$. This leads to an interesting contrast between extremal functions and extremal maps from the smoothness view point. In fact, an extremal function to (\ref{B-opt}) is always of $C^1$ class since, modulo a constant scale factor, it solves the elliptic equation

\begin{equation} \label{E.2}
- \Delta_g u + \lambda \beta(x) u = |u|^{2^*-2} u \ \ {\rm on}\ \ M
\end{equation}
for some constant $\lambda > 0$. Nevertheless, extremal maps to (\ref{B-opt-v}) do not need to be smooth and much less satisfying any equations, unless $F(t)$ and $G(x,t)$ are of $C^1$ class on $t$. In this specific situation, an extremal map $U=(u_1,\ldots, u_k)$ solves, up to a scaling, the potential type elliptic system

\begin{equation} \label{S.1}
-\Delta_g u_i + \frac{\lambda}{2} \frac{\partial G(x,U)}{\partial t_i} = \frac{1}{2^*} \frac{\partial F(U)}{\partial t_i} \ \ {\rm on}\ \ M,\ \ i=1,\ldots,k
\end{equation}
for some constant $\lambda > 0$, and so, by elliptic regularity results, each component $u_i$ is of $C^1$ class. This dichotomy leads to an immediate impact on existence and compactness problems of extremal maps, because the development of the scalar $AB$ program is completely based on concentration analysis of critical points of energy type functionals whose Euler-Lagrange equations are like (\ref{E.2}). In short, when $k \geq 2$ we are facing a new situation where are prevented of directly using variational techniques and arguments of any kind involving PDEs. Another important difference concerns the positivity of solutions in PDEs. While an extremal function has defined sign on $M$, generally components of extremal maps have undeterminate sign, even when $F(t)$ and $G(x,t)$ are smooth functions with symmetry properties. This difference occurs because, contrary to elliptic equations where the maximum principle can be applied, there is no maximum principle for potential type elliptic systems. Since positivity plays an essential role when $k=1$, the maps context becomes more involved than the functions one. All these comments suggest that the questions (A)-(E) are embedded in a broader structure whose understanding is a challenge.

We shall provide answers to (A)-(E) into the general perspective proposed here. Part of our main results are established for dimensions $n \geq 5$, as for instance:

\begin{itemize}

\item[$(\beta.1)$] the continuity of the map

\[
(F,G,g) \mapsto \mathcal{B}_0(n,F,G,g)
\]

\n on appropriate spaces;

\item[$(\beta.2)$] the duality between existence of extremal maps to (\ref{B-opt-v}) and the explicit value of $\mathcal{B}_0(n,F,G,g)$;

\item[$(\beta.3)$] a compactness result of extremal maps related to families of continuous potential functions and Riemannian metrics.

\end{itemize}

\n In particular, our results extend in a comprehensive way much of the scalar $AB$ program developed so far. We also expect some theoretical tools used in this work may be useful in related contexts within the nonlinear analysis on manifolds.

In a nutshell, our strategy to overcome the above-described obstacles is based on the following tools:

\begin{itemize}

\item[$\star$] De Giorgi-Nash-Moser type estimates applied to solutions of systems like (\ref{S.1}) (Proposition \ref{P.4}) and to minimizers of non-smooth functionals (Proposition \ref{P.5}) with constants depending only on the respective $C^0$-norms of $F$ and $G$ on $\mathbb{S}^{k-1}$ and $M \times \mathbb{S}^{k-1}$;

\item[$\star$] An approximation scheme of potential functions $F$ and $G$ in the $C^0_{loc}$-topology by $C^1$ functions of same kind  (Proposition \ref{P.6});

\item[$\star$] A local geometric potential type Sobolev inequality uniformly satisfied for suitable families of potential functions $(F_\alpha)$ and $(G_\alpha)$  (Proposition \ref{P.7}).

\end{itemize}

\n These ingredients are essential in the proof of $(\beta.1)$-$(\beta.3)$ and, in addition, the proof of $(\beta.2)$ and $(\beta.3)$ rely on the continuity result $(\beta.1)$.\\

\n {\bf 1.3 Main theorems}. Our most striking results are summarized in the next three theorems. But first we need a few notations.

Let $M$ be a compact differentiable manifold of dimension $n$. Denote by ${\cal M}^n$ the space of smooth Riemannian metrics on $M$ endowed with the $C^2$-topology. For $k \geq 1$, we represent by $C^0_k(M)$ the space $C^0(M) \times \ldots \times C^0(M)$ of continuous $k$-maps on $M$ equipped with the usual product topology, by ${\cal F}_k$ the set of positive continuous functions on $\R^k$ homogeneous of degree $2^*$ with the induced topology of $C_{loc}^0(\R^k)$ and by ${\cal G}_k$ the set of positive continuous functions on $M \times \R^k$ homogeneous of degree $2$ on the second variable with the induced topology of $C_{loc}^0(M \times \R^k)$. For $F \in {\cal F}_k$, $G \in {\cal G}_k$ and $g \in {\cal M}^n$, we denote by ${\cal E}_k(F,G,g)$ the set of $F$-normalized extremal maps to (\ref{B-opt-v}). Finally, given $g \in {\cal M}^n$, denote by $L^{2^*}_k(M)$ the Riemannian $k$-vector Lebesgue space $L^{2^*}(M) \times \ldots \times L^{2^*}(M)$ endowed with the usual product topology.

\begin{teor} {\rm(Continuity).} \label{Teo.1}
Let $M$ be a compact differentiable manifold of dimension $n \geq 5$. For any integer $k \geq 1$, the map $(F,G,g) \in \mathcal{F}_k  \times \mathcal{G}_k \times {\cal M}^n \mapsto \mathcal{B}_0(n,F,G,g)$ is continuous.
\end{teor}

\begin{teor} {\rm(Duality).} \label{Teo.2}
Let $M$ be a compact differentiable manifold of dimension $n \geq 5$. For any integer $k \geq 1$ and $(F,G,g) \in \mathcal{F}_k  \times \mathcal{G}_k \times {\cal M}^n$, at least, one of the following assertions holds:

\begin{itemize}

\item[(I)*] (\ref{B-opt-v}) admits an extremal map;

\item[(II)*] $\mathcal{B}_0(n,F,G,g) G(x_0,t_0) = \frac{n-2}{4(n-1)} \mathcal{A}_0(n,F) S_g(x_0)$ for some $x_0 \in M$ and some maximum point $t_0$ of $F$ on $\s^{k-1}$.

\end{itemize}

\end{teor}

It follows from Theorem \ref{Teo.1} that the set of triples $(F,G,g)$ satisfying the condition (II)* of Theorem \ref{Teo.2} is closed in $\mathcal{F}_k  \times \mathcal{G}_k \times {\cal M}^n$. An easy consequence to Theorem \ref{Teo.2} is that sharp and saturated inequalities like (\ref{B-opt-v}) always admits extremal maps when $n \geq 5$ and $(M,g)$ has nonpositive scalar curvature. Another possible application, which follows from Theorem \ref{Teo.2} and the resolution of the Yamabe problem is that if $n \geq 5$, $G(x,t)$ does not depend on $x$, and $(M,g)$ has constant scalar curvature, then (\ref{B-opt-v}) admits extremal maps (see Section 2).

\begin{teor} {\rm(Compactness).} \label{Teo.3}
Let $M$ be a compact differentiable manifold of dimension $n \geq 5$. Let $k \geq 1$ be an integer number, $((F_\alpha,G_\alpha,g_\alpha))$ be a sequence converging to $(F,G,g)$ in $\mathcal{F}_k  \times \mathcal{G}_k \times {\cal M}^n$ and $(U_\alpha)$ be a sequence of extremal maps with $U_\alpha \in {\cal E}_k(F_\alpha,G_\alpha,g_\alpha)$. If the triple $(F,G,g)$ satisfies

\[
\mathcal{B}_0(n,F,G,g) G(x,t_0) > \frac{n-2}{4(n-1)} \mathcal{A}_0(n,F) S_{g}(x)
\]

\n for all $x \in M$ and all maximum point $t_0$ of $F$ on $\s^{k-1}$, then $(U_\alpha)$ is $L^{2^*}_k$-compact. Moreover, if $((F_\alpha,G_\alpha,g_\alpha))$ converges in $C^1_{loc}(\R^k) \times C^0(M, C^1_{loc}(\R^k)) \times {\cal M}^n$, then $(U_\alpha)$ is $C^0_k$-compact.
\end{teor}

For dimensions $n \geq 5$, Theorems \ref{Teo.2} and \ref{Teo.3} extend in a unified framework Theorem 1 of \cite{DjDr} ($k = 1$ and $\beta \equiv 1$), Theorem 3 of \cite{HV2} ($k = 1$ and $\beta$ general) and Theorem 0.1 of \cite{HH3} ($k \geq 2$, $F(t) = |t|_{2^*}^{2^*}$ and $G(x,t) = \langle A(x) t, t \rangle$). In addition, Theorem \ref{Teo.1} and \ref{Teo.3} generalize, respectively, Theorem 1.1 and Corollary 1.2 of \cite{EM6} ($k = 1$ and $\beta \equiv 1$). We emphasize that Theorem \ref{Teo.1} is new even when $k = 1$. In order to highlight this specific case, denote by $C^0_+(M)$ the cone of positive continuous functions on $M$ endowed with the usual uniform topology. Noting that Proposition \ref{P.7} is true in dimension $n = 4$ when $k = 1$ (see \cite{CHV}), we can state:

\begin{cor} {\rm(Continuity for $k=1$).} \label{Cor.1}
Let $M$ be a compact differentiable manifold of dimension $n \geq 4$. Then the map $(\beta,g) \in C_+^0(M) \times {\cal M}^n \mapsto B_0(n,\beta,g)$ is continuous.
\end{cor}

\begin{cor} {\rm(Compactness for $k=1$).} \label{Cor.2}
Let $M$ be a compact differentiable manifold of dimension $n \geq 4$. Let $((\beta_\alpha,g_\alpha))$ be a sequence convergent to $(\beta, g)$ in $C_+^0(M) \times {\cal M}^n$ and $(u_\alpha)$ be a sequence with $u_\alpha \in E(\beta_\alpha,g_\alpha)$. If the couple $(\beta, g)$ satisfies

\[
B_0(n,\beta,g) \beta(x) > \frac{n-2}{4(n-1)} A_0(n) S_{g}(x)
\]

\n for all $x \in M$, then the sequence $(u_\alpha)$ is $C^0$-compact.
\end{cor}

\vs{0.3cm}

\n {\bf 1.4 Organization of the paper}. Section 2 is central from the view point of solution to various gaps arising in the vector setting. In it, we answer partially the questions (A) and (B) of Subsection 1.2, present an overview on sharp potential type Sobolev inequalities and develop all theoretical tools needed to the proof of the main theorems. For instance, Propositions \ref{P.2} and \ref{P.3} are directly linked to (A) and (B). The main tools are expressed in Propositions \ref{P.4}, \ref{P.5}, \ref{P.6} and \ref{P.7}. Proposition \ref{P.4} deals with De Giorgi-Nash-Moser type estimates to solutions of potential type systems like (\ref{S.1}) with constants depending only on quantities related to the continuity of $F$ and $G$. Proposition \ref{P.5} provides such estimates for minimizers of (only continuous) non-smooth functionals. Proposition \ref{P.6} states that positive homogeneous continuous functions can always be approximated by $C^1$ functions of same kind. Finally, Proposition \ref{P.7} deals with a local potential type Sobolev inequality involving scalar curvature which is closely related to a geometric Sobolev inequality obtained by Druet in \cite{D3}. The proof is inspired in this work and is based on a thorough concentration analysis on least energy critical points of certain smooth energy functional. Here, the complete knowledge of extremal maps to (\ref{Euc}), given in Proposition \ref{P.1}, plays an important role. Proposition \ref{P.7} is the only place where we use the assumption on the dimension $n \geq 5$. When $k \geq 2$, we believe in its validity up to $n = 4$ such as in the case $k=1$, see remark made on the page $2535$ of \cite{CHV}. Propositions \ref{P.4}, \ref{P.6} and \ref{P.7} are fundamental in the proofs of all theorems and Proposition \ref{P.5} is essential only in the proof of Theorem \ref{Teo.3}. Part of the ideas employed in the proof of Theorem \ref{Teo.1} is also exploited in the proofs of Theorems \ref{Teo.2} and \ref{Teo.3} and, therefore, these last ones will be made more concise. The proofs of Theorems \ref{Teo.1}, \ref{Teo.2} and \ref{Teo.3} will stretch from Sections 3 to 5. The general outlines of the proofs of Proposition \ref{P.7} and Theorems \ref{Teo.1} and \ref{Teo.2} are similar in essence and are inspired in key ideas developed in the scalar theory. We refer the reader to the works \cite{AuLi}, \cite{DjDr}, \cite{D2}, \cite{HV1} and \cite{HV2} where various approaches are presented. Although, as already mentioned, important difficulties arise in the vector setting mainly caused by our will to treat problems in a very general perspective, the common core of proofs can be summarized in four steps:

\begin{itemize}

\item[$(\gamma.1)$] to suppose, by contradiction, that each statement fails;

\item[$(\gamma.2)$] to go in search of least energy critical points $(U_\alpha)$ to suitable smooth energy functionals;

\item[$(\gamma.3)$] to perform blow-up and concentration analysis on $(U_\alpha)$;

\item[$(\gamma.4)$] to deduce a contradiction via integral estimates obtained around a concentration point;

\end{itemize}

\n most of the rest is just technique. Propositions \ref{P.4}, \ref{P.6} and \ref{P.7} are directly linked to the steps $(\gamma.3)$, $(\gamma.2)$ and $(\gamma.4)$, respectively. Section 6 is devoted to examples and counterexamples exploring the necessity of topology $C^2$ on the space of Riemannian metrics ${\cal M}^n$ for continuity of the map $(F,G,g) \in \mathcal{F}_k  \times \mathcal{G}_k \times {\cal M}^n \mapsto \mathcal{B}_0(n,F,G,g)$, the dichotomy (I)*-(II)* provided in Theorem \ref{Teo.2} and the compactness and non-compactness of extremal maps to (\ref{B-opt-v}). For convenience of the reader, we organize in an appendix some important basic results such as non-smooth vector versions of the Br\'{e}zis-Lieb lemma \cite{BrLi} and of a compactness-concentration principle \cite{Li}, among others needed in proofs.\\

\n {\bf 1.5 Further remarks and open problems}. Before starting the proofs, we make a few final comments about the results and sketch some open problems. First, Theorem \ref{Teo.1} is indeed a sharp result from the topological view point since the $C^2$-topology condition on ${\cal M}^n$ is the possible weaker one among all $C^k$-topologies. A simple counter-example showing such a claim is provided in Section 6. There is a number of open interesting problems concerning with the dichotomy (I)*-(II)* stated in Theorem \ref{Teo.2}. A first question is to know if all possible situations can hold for $k \geq 2$ as in the scalar case. If yes, to exhibit examples and to contrast them with those well-known ones when $k=1$. If (I)* fails (i.e. ${\cal E}_k(F,G,g) = \emptyset$), it is still possible to ask if ${\cal E}_k(F,G,\tilde{g}) \neq \emptyset$ for some conformal metric $\tilde{g}$. The vector dichotomy seems to be more intricate than the corresponding scalar one. This claim is reinforced by examples as the round unit $n$-sphere $\s^n$ and functions $G(x,t)$ depending only on $t$ where (I)* happens and by far (II)* is not clear to be true (see Section 2). According to some answers, other questions can possibly be explored toward a proper understanding of (I)*-(II)*. Theorem \ref{Teo.3} is also sharp in the same sense of Theorem \ref{Teo.1}.

A natural problem on a compact Riemannian manifold of dimension $n \geq 4$ not conformal to the round sphere is if the sequence $(U_\alpha)$ of Theorem \ref{Teo.3} remains compact when, for each $\alpha > 0$,

\[
\mathcal{B}_0(n,F_\alpha,G_\alpha,g_\alpha) G_\alpha(x_\alpha,t_\alpha) = \frac{n-2}{4(n-1)} \mathcal{A}_0(n,F_\alpha) S_{g_\alpha}(x_\alpha)
\]

\n for some $x_\alpha \in M$ and some maximum point $t_\alpha$ of $F_\alpha$ on $\s^{k-1}$. This one seems quit difficult and is open for $k \geq 1$.

Regarding the case of $3$-dimensional compact manifolds, all questions treated here are open. One may expect a sort of low-dimensional phenomenon. As first pointed out by Br\'{e}zis and Nirenberg in \cite{BrNi}, an elliptic type situation may change drastically when passing from low dimensions, in this case $n = 3$, to high dimensions, say $n \geq 4$. When dealing with sharp Sobolev inequalities on compact Riemannian manifolds, this phenomenon happened in the works of Hebey \cite{H2} and Druet, Hebey and Vaugon \cite{DHV}.

Another potentially interesting direction is the development of a theory on sharp potential type Riemannian $L^p$-Sobolev inequalities for $p \neq 2$ and for a class of potential functions $F$ and $G$ satisfying weak regularity conditions. The corresponding questions (A)-(E) in this new setting are widely addressed in \cite{EM8}. Unlike the case $p=2$ (see (II)* of Theorem \ref{Teo.2}), we establish general results of existence and compactness of extremal maps without involving any conditions on the maximum points of $F$ on $\mathbb{S}_p^{k-1} := \{t \in \R^k:\; |t|_p^p = \sum_{i=1}^k |t_i|^p = 1\}$. The absence of such an influence is not so surprising for $p \neq 2$ according to the scalar $AB$ program to $L^p$-Sobolev inequalities. Moreover, we investigate geometric and topological obstructions to validity in function of geometric and/or topological properties of $M$ and the range of values of $p$ in the same spirit of the scalar $L^p$-theory.\\

\section{Miscellaneous and central tools}

Let $k \geq 1$ be an integer number. Throughout this work, we will adopt the notations ${\cal F}_k$ and ${\cal G}_k$ for the space of positive continuous functions, respectively, from $\R^k$ onto $\R$ and from $M \times \R^k$ onto $\R$, and homogeneous of degrees $2^*$ and $2$. We also carry all notations of best constants and Sobolev spaces introduced in the previous section.\\

\n {\bf 2.1 Sharp potential type Sobolev inequalities}. In this subsection we collect some facts related to sharp potential type Sobolev inequalities which follow in a simple way. We begin with the characterization of extremal maps to sharp Euclidean potential type Sobolev inequalities:

\begin{propo}\label{P.1}
For $n \geq 3$ and $F \in {\cal F}_k$, we have:

\begin{itemize}

\item[(a)] $\mathcal{A}_0(n,F) = M_F^{2/2^*} A_0(n)$, where $M_F = \max \limits_{|t|=1} F(t)$,

\item[(b)] $\mathcal{A}_0(n,F)$ is achieved exactly by maps of the form $U_0 = t_0 u_0$, where $t_0 \in \s^{k-1}$ is a maximum point of $F$ and $u_0$ is an extremal function to $A_0(n)$.

\end{itemize}
\end{propo}

\n {\bf Proof of Proposition \ref{P.1}.} The $2^*$-homogeneity of $F$ yields

\[
F(t)\leq M_{F} \left( \sum _{i=1}^{k} t_i^2 \right)^{2^*/2}
\]
for all $t \in \mathbb{R}^n$. By Minkowski's inequality, for any $U \in \mathcal{D}^{1,2}_k(\mathbb{R}^n)$,

\begin{eqnarray}
\left( \int_{\mathbb{R}^n} F(U)\; dx \right)^{2/2^*} &\leq& M_F^{2/2^*} \left( \int_{\mathbb{R}^n} \left( \sum_{i=1}^{k} u_i^{2}\right)^{2^*/2}\; dx \right)^{2/2^*} \label{des-1}\\
&\leq&  M_F^{2/2^*} \sum_{i=1}^{k} \left(\int _{\mathbb{R}^n}|u_i|^{2^*}\; dx \right)^{2/2^*} \nonumber \\
&\leq& M_F^{2/2^*} A_0(n) \sum_{i=1}^k \int _{\mathbb{R}^n} |\nabla u_i|^{2}\; dx \nonumber\\
&=& M_F^{2/2^*} A_0(n) \int _{\mathbb{R}^n} |\nabla U|^{2}\; dx\; , \nonumber
\end{eqnarray}

\n and this inequality clearly implies

\[
\mathcal{A}_0(n,F) \leq M_F^{2/2^*} A_0(n)\; .
\]

\n Consider now the map $U_0 = t_0 u_0 \in \mathcal{D}^{1,2}_k(\mathbb{R}^n)$ where $t_0 \in \mathbb{S}^{k-1}$ is such that $M_F = F(t_0)$ and $u_0 \in \mathcal{D}^{1,2}(\mathbb{R}^n)$ is an extremal function to $A_0(n)$. Since

\begin{eqnarray*}
\left( \int_{\mathbb{R}^N} F(U_0)\; dx \right)^{2/2^*} &=& M_F^{2/2^*} \left( \int_{\mathbb{R}^n} |u_0|^{2^*} dx \right)^{2/2^*} \\
&=& M_F^{2/2^*} A_0(n) \int_{\mathbb{R}^n} |\nabla u_0|^{2}\; dx\\
&=& M_F^{2/2^*} A_0(n) \int _{\mathbb{R}^n} |\nabla t_0 u_0|^{2}\; dx \\
&=& M_F^{2/2^*} A_0(n) \int _{\mathbb{R}^n} |\nabla U_0|^{2}\; dx\; ,
\end{eqnarray*}

\n we have

\[
\mathcal{A}_0(n,F) = M_F^{2/2^*} A_0(n)\; .
\]

\n Note that $\mathcal{A}_0(n,F)$ is achieved by maps $U_0 = t_0 u_0$ as chosen above. It remain to show that an arbitrary extremal map $U$ to $\mathcal{A}_0(n,F)$ can always be placed in this form. In fact, $U$ satisfies (\ref{des-1}) with equality in all steps. Remark also that the second equality corresponds to Minkowski's inequality. So, there exist $t \in \mathbb{S}^{k-1}$ and a function $u \in \mathcal{D}^{1,2}(\mathbb{R}^n)$ such that $U = t u$. Finally, from the first equality, one gets $F(t) = M_F$, and from the third one, one easily checks that $u$ is an extremal function to $A_0(n)$.\; \bl \\

\begin{rem} \label{R.1} An alternative way to attack the problem of classifying all extremal maps for Euclidean potential type Sobolev inequalities in the smooth case is to reduce to the radial situation, via an usual symmetrization argument, and after identifying the solutions of the corresponding Euler-Lagrange system.
\end{rem}

\begin{rem} \label{R.2} It follows directly from the part (b) of Proposition \ref{P.1} that the extremal maps to (\ref{Euc}) are precisely given by

\[
U_0(x) = a t_0 w(b(x - x_0)),
\]

\n where $a, b \in \mathbb{R} \setminus \{0\}$, $x_0 \in \mathbb{R}^n$, $t_0$ is a maximum point of $F$ on $\s^{k-1}$ and $w$ is as in (\ref{fe1a}).
\end{rem}

In the next two results, we select some basic facts about sharp potential type Riemannian Sobolev inequalities.

\begin{propo}\label{P.2}
Let $(M,g)$ be a smooth compact Riemannian manifold of dimension $n \geq 3$, $F \in {\cal F}_k$ and $G \in {\cal G}_k$. The following assertions hold:

\begin{itemize}

\item[(a)] $\mathcal{A}_0(n,F,G,g) = \mathcal{A}_0(n,F)$,

\item[(b)] There exists a constant $\mathcal{B} \in \R$ such that the sharp inequality (\ref{A-opt-v}) is valid.

\end{itemize}

\end{propo}

\n {\bf Proof of Proposition \ref{P.2}.} Consider constants $\mathcal{A}, \mathcal{B}$ such that, for any $U \in H^{1,2}_k(M)$,

\begin{equation}\label{a}
\left( \int _M F(U)\; dv_g  \right)^{2/2^*} \leq \mathcal{A} \int_M|\nabla_g U|^2\; dv_g + \mathcal{B} \int_M G(x,U)\; dv_g\; .
\end{equation}

\n Choosing now maps of the form $U = t_0 u$ in (\ref{a}), with $t_0 \in \mathbb{S}^{k-1}$ such that $F(t_0) = M_F$ and $u \in H^{1,2}(M)$, and dividing both sides by $M_F^{2/2^*}$, one gets

\[
\left( \int _M |u|^{2^*}\; dv_g  \right)^{2/2^*} \leq M_F^{-2/2^*}\mathcal{A} \int_M|\nabla_g u|^2\; dv_g + M_F^{-2/2^*} M_G \mathcal{B} \int_M u^2\; dv_g\; ,
\]

\n where $M_G = \max \limits_{(x,t) \in M \times \s^{k-1}} G(x,t)$. Since $A_0(n,1,g) = A_0(n)$, one then has

\[
\mathcal{A}_0(n,F,G,g) \geq M_F^{2/2^*} A_0(n) = \mathcal{A}_0(n,F)\; .
\]

\n Conversely, since (\ref{A-opt}) is always valid, we have

\begin{eqnarray}
\left( \int_{M} F(U)\; dv_g \right)^{2/2^*} &\leq& M_F^{2/2^*} \left( \int_{M} \left( \sum_{i=1}^{k} u_i^{2}\right)^{2^*/2}\; dv_g \right)^{2/2^*} \label{des-2} \\
&\leq&  M_F^{2/2^*} \sum_{i=1}^{k} \left(\int_{M} |u_i|^{2^*}\; dv_g \right)^{2/2^*} \nonumber \\
&\leq& M_F^{2/2^*} A_0(n) \sum_{i=1}^k \int_{M} |\nabla_g u_i|^{2}\; dv_g + M_F^{2/2^*} B \sum_{i=1}^k \int_{M} u_i^{2}\; dv_g \nonumber \\
&=& M_F^{2/2^*} A_0(n) \int_{M} |\nabla_g U|^{2}\; dv_g + M_F^{2/2^*} B \int_{M} |U|^{2}\; dv_g \nonumber
\end{eqnarray}

\n for any $U \in H^{1,2}_k(M)$. Therefore, it follows simultaneously from (\ref{des-2}) that

\[
\mathcal{A}_0(n,F,G,g) \leq M_F^{2/2^*} A_0(n) = \mathcal{A}_0(n,F)
\]

\n and the conclusion of part (b).\bl \\

\begin{rem} \label{R.3} The part (a) of Proposition \ref{P.2} furnishes complete information about the best constant $\mathcal{A}_0(n,F,G,g)$. In particular, it neither depends on the potential function $G$ nor on the geometry of the manifold $M$. Furthermore, the map $(F,G,g) \in \mathcal{F}_k  \times \mathcal{G}_k \times {\cal M}^n \mapsto \mathcal{A}_0(n,F,G,g)$ is clearly continuous. Note also that, thanks to the part (b) of Proposition \ref{P.2}, the best constant $\mathcal{B}_0(n,F,G,g)$ is always well-defined.
\end{rem}

\begin{propo}\label{P.3}
Let $(M,g)$ be a smooth compact Riemannian manifold of dimension $n \geq 4$, $F \in {\cal F}_k$ and $G \in {\cal G}_k$. The following assertions hold:

\begin{itemize}

\item[(a)] For any $x \in M$ and $t_0 \in \mathbb{S}^{k-1}$ with $F(t_0) = M_F$,

\[
\mathcal{B}_0(n,F,G,g) G(x,t_0) \geq \frac{n-2}{4(n-1)} \mathcal{A}_0(n,F) S_g(x)\, ,
\]

\item[(b)] If $G(x,t)$ does not depend on the variable $x$ and $(M,g)$ has constant scalar curvature, then (\ref{B-opt-v}) possesses extremal maps.

\end{itemize}

\end{propo}

An immediate consequence is:

\begin{cor}\label{C.1}
Let $M = \mathbb{S}^{n}$ be the round unit $n$-sphere of dimension $n \geq 4$, $F \in {\cal F}_k$ and $G \in {\cal G}_k$. Then, (\ref{B-opt-v}) possesses extremal maps whenever $G(x,t)$ does not depend on the variable $x$.
\end{cor}

\n {\bf Proof of Proposition \ref{P.3}.} The proof of the part (a) is simple. It suffices to take maps of the form $U = t_0 u$ in the inequality (\ref{B-opt-v}) as in the previous proof. In fact, this choice yields

\[
\left( \int_{M} |u|^{2^*}\; dv_g \right)^{2/2^*} \leq A_0(n) \int_{M} |\nabla_g u|^2 \; dv_g + M_F^{-2/2^*}
\mathcal{B}_0(n,F,G,g) \int_{M} G(x,t_0) u^2\; dv_g\; .
\]

\n Set $\beta(x) = G(x,t_0)$. From the definition of $B_0(n,\beta,g)$, we have

\[
B_0(n,\beta,g) \leq M_F^{-2/2^*} \mathcal{B}_0(n,F,G,g)\; .
\]

\n So, the part (a) follows readily from the geometric estimate (\ref{G-E-s}) and Proposition \ref{P.1}. Already the part (b) can be obtained from a combination among the part (a), the solution of the Yamabe problem and Theorem \ref{Teo.2}. Indeed, thanks to the works of Aubin \cite{Au3} and Schoen \cite{S}, there exists a positive solution $u_0 \in C^{\infty}(M)$, normalized by unit $L^{2^*}$-norm, of the equation

\[
- \frac{4(n-1)}{n-2} \Delta_g u + S_g u = \mu_g(M) u^{2^*-1} \ \ {\rm on}\ \ M\; ,
\]

\n where $\mu_g(M)$ denotes the Yamabe invariant associated to $(M,g)$. For Riemannian manifolds $(M,g)$ non-conformal to the round unit $n$-sphere $\s^n$, one knows that

\[
\mu_g(M) < \frac{4(n-1)}{(n-2)A_0(n)}\, .
\]

\n Letting $U_0 = t_0 u_0$, where $t_0 \in \mathbb{S}^{k-1}$ is such that $F(t_0) = M_F$, we derive from the Yamabe equation,

\[
\left( \int_M F(U_0) \; dv_g\right)^{2/2^*} > \mathcal{A}_0(n,F) \int_M |\nabla_g U_0|^2 \; dv_g + \frac{n-2}{4(n-1) G(t_0)} \mathcal{A}_0(n,F) S_g \int_M G(U_0) \; dv_g\; .
\]

\n But this inequality implies that

\[
\mathcal{B}_0(n,F,G,g) G(t_0) > \frac{n-2}{4(n-1)} \mathcal{A}_0(n,F) S_{g}
\]

\n for all maximum point $t_0$ of $F$ on $\mathbb{S}^{k-1}$, so that (\ref{B-opt-v}) admits an extremal map by Theorem \ref{Teo.2}. Finally, if $(M,g)$ is conformal to the round unit $n$-sphere $\s^n$, we have

\[
\mu_g(M) = \frac{4(n-1)}{(n-2)A_0(n)}\; .
\]

\n If the inequality above holds for all $t_0$, evoking once again Theorem \ref{Teo.2}, it follows the existence of an extremal map to (\ref{B-opt-v}). Otherwise, if equality happens for some $t_0$, then the corresponding map $U_0$ is an extremal map since

\[
\left( \int_M u_0^{2^*} \; dv_g\right)^{2/2^*} = A_0(n) \int_M |\nabla_g u_0|^2 \; dv_g + \frac{n-2}{4(n-1)} A_0(n) S_g \int_M u_0^2 \; dv_g
\]

\n yields

\[
\left( \int_M F(U_0) \; dv_g\right)^{2/2^*} = \mathcal{A}_0(n,F) \int_M |\nabla_g U_0|^2 \; dv_g + \frac{n-2}{4(n-1) G(t_0)} \mathcal{A}_0(n,F) S_g \int_M G(U_0) \; dv_g\, .
\]

\bl\\

\n {\bf 2.2 De Giorgi-Nash-Moser estimates for systems and variants}. The main goal of this subsection is to establish De Giorgi-Nash-Moser estimates for weak solutions of potential type elliptic systems and for minimizers of non-differentiable energy functionals on $H^{1,2}_k(M)$. Such estimates play a strategic role in all work.

We begin by dealing with weak solutions $U = (u_1, \ldots , u_k)$ of the following potential type elliptic system:

\begin{equation} \label{syssa}
- \mathcal{A} \Delta_{g} u_i + \frac{1}{2} \mathcal{B} \frac{\partial G(x,U)}{\partial t_i} = \frac{1}{2^*} \frac{\partial
F(U)}{\partial t_i} \ \ {\rm in}\ \ \Omega \subset M \,,\ \ i = 1, \ldots ,k,
\end{equation}

\n where $\Omega$ is an open of the $n$-dimensional smooth compact Riemannian manifold $(M,g)$, $\mathcal{A}$ and $\mathcal{B}$ are positive constants and $F \in \mathcal{F}_k$ and $G \in \mathcal{G}_k$ are functions of $C^1$ class.

\begin{propo} \label{P.4}
Let $M$ be a smooth compact Riemannian manifold of dimension $n \geq 3$ and $\lambda$, $\mathcal{A}_1$ and $C_1$ be positive constants. Consider a metric $g \in {\cal M}^n$, constants $\mathcal{A}$ and $\mathcal{B}$ and $C^1$ potential functions $F \in \mathcal{F}_k$ and $G \in \mathcal{G}_k$ satisfying

\[
g_{ij} \xi_i \xi_j \geq \lambda |\xi|^2,\ \ ||g_{ij}||_{C^0} \leq \lambda^{-1}
\]

\n for all $\xi \in \R^n$ and $i,j = 1, \ldots, n$ in a fixed coordinates system on $M$,

\[
\mathcal{A} \geq \mathcal{A}_1,\ \ \mathcal{B} \geq 0\; ,
\]

\[
F(t) \leq C_1 |t|^{2^*}
\]

\n for all $t \in \R^k$. Then, given constants $q > 2^*$ and $p, K > 0$, there exists a positive constant $C_0$, depending only on $n$, $q$, $p$, $K$, $\lambda$, $\mathcal{A}_1$, and $C_1$ such that for any $\delta > 0$, any point $x_0 \in M$, any domain $\Omega \subset M$ and any solution $U \in H^{1,2}_k(\Omega)$ of the system (\ref{syssa}) satisfying $U = 0$ on $B_g(x_0, 2\delta) \cap \partial \Omega$ and

\[
||U||_{L_k^{q}(B_g(x_0, 2\delta))} \leq K\, ,
\]

\n one has

\begin{equation} \label{lim3}
\sup_{B_g(x_0, \delta) \cap \Omega} |U| \leq C_0 \delta^{-n/p}  \left( \int_{B_g(x_0, 2\delta) \cap \Omega} |U|^{p}\; dv_g \right)^{1/p}\, .
\end{equation}

\end{propo}

\begin{rem} \label{R.4} The strength of this proposition lies in the fact of the constant $C_0$ to rely continuously on $F$ and not on partial derivatives of this function.
\end{rem}

\n {\bf Proof of Proposition \ref{P.4}.} A idea of the proof is to seek a differential inequality satisfied by the scalar function $|U|$ in terms of $F$ and $G$. Surprisingly, we discover such an inequality thanks to homogeneity properties. First, for each $\varepsilon > 0$, we set
\[
v_\varepsilon := \left( \sum_{i=1}^k u_i^2 + \varepsilon \right)^{1/2}\, .
\]

\n For any nonnegative function $\varphi \in C^1(M)$, we have

\begin{eqnarray*}
\int_{M} \nabla_g v_\varepsilon \cdot \nabla_g \varphi\; dv_g &=& \int_{M} v_\varepsilon^{-1} (\sum_{i=1}^k u_i \nabla_g u_i) \cdot \nabla_g \varphi\; dv_g \\
&=& \int_{M} \sum_{i=1}^k \nabla_g u_i \cdot \left( \nabla_g (u_i v_\varepsilon^{-1} \varphi ) - \varphi \nabla_g (u_i v_\varepsilon^{-1} ) \right)\; dv_g \\
&=& \sum_{i=1}^k \int_{M} \nabla_g u_i \cdot \nabla_g (u_i v_\varepsilon^{-1} \varphi)\; dv_g  \\
&& - \int_{M} v_\varepsilon^{-3} \varphi \left( v_\varepsilon^{2} \sum_{i=1}^k |\nabla_g u_i|^2 - |\sum_{i=1}^k u_i \nabla_g u_i|^2 \right) \; dv_g\, . \end{eqnarray*}

\n Using now the system (\ref{syssa}) on the first right-hand side integral and the Cauchy-Schwarz inequality on the second right-hand side integral, we deduce that

\begin{eqnarray*}
\int_{M} \nabla_g v_\varepsilon \cdot \nabla_g \varphi\; dv_g &\leq& \sum_{i=1}^k \int_{M} \nabla_g u_i \cdot \nabla_g (u_i v_\varepsilon^{-1} \varphi)\; dv_g\\
&=& \int_{M} \sum_{i=1}^k \left( \frac{1}{2^*} \mathcal{A}^{-1} \frac{\partial F(U)}{\partial t_i} u_i - \frac{1}{2} \mathcal{A}^{-1}
\mathcal{B} \frac{\partial G(x,U)}{\partial t_i} u_i \right) (v_\varepsilon^{-1} \varphi)\; dv_g \\
&=& \mathcal{A}^{-1} \int_{M} v_\varepsilon^{-1} F(U) \varphi\; dv_g - \mathcal{A}^{-1}
\mathcal{B} \int_{M} v_\varepsilon^{-1} G(x,U) \varphi\; dv_g\, .
\end{eqnarray*}

\n Letting $\varepsilon \rightarrow 0$ in the preceding inequality, we derive

\[
\int_{M} \nabla_g |U| \cdot \nabla_g \varphi\; dv_g \leq \mathcal{A}^{-1} \int_{M} |U|^{-1} F(U) \varphi\; dv_g - \mathcal{A}^{-1}
\mathcal{B} \int_{M} |U|^{-1} G(x,U) \varphi\; dv_g\, .
\]

\n Summarizing, we have shown that

\begin{equation} \label{R}
- \mathcal{A} \Delta_g |U| + \mathcal{B} |U|^{-1} G(x,U) \leq |U|^{-1} F(U)\, .
\end{equation}

\n So, the additional assumptions on $\mathcal{A}$, $\mathcal{B}$ and $F$ applied to (\ref{R}) readily yield

\[
- \Delta_g |U| \leq K_1 |U|^{2^* - 1}\ \ {\rm in}\ \ \Omega\, ,
\]

\n where $K_1 = \mathcal{A}_1^{-1} C_1$. The conclusion then follows from classical De Giorgi-Nash-Moser estimates when studying inequations of kind

\[
- \Delta_{g} u + c u \leq f\ \ {\rm in}\ \ \Omega
\]

\n with $c = - K_1 |U|^{2^* - 2}$ and $f = 0$. (see Serrin \cite{Se} or Trudinger \cite{T}).\bl \\

We now focus De Giorgi-Nash-Moser estimates for minimizers of non-differentiable energy functionals. Given a $n$-dimensional smooth compact Riemannian manifold $(M,g)$ and potential functions $F \in \mathcal{F}_k$ and $G \in \mathcal{G}_k$, we consider the functional

\begin{equation} \label{E.7}
J(U) = \mathcal{A}_0(n,F) \int_M |\nabla_g U|^2\; dv_g + \mathcal{B} \int_M G(x,U)\; dv_g
\end{equation}

\n on $\Lambda = \{U \in H^{1,2}_k(M) :\, \int_M F(U) dv_g = 1\}$ and their possible minimizers $U$. Note that since $F$ and $G$ are only continuous, we cannot differentiate $J$ at $U$. In particular, their possible minimizers do not satisfy systems like (\ref{syssa}). Even so, it is still possible to establish the following result:

\begin{propo} \label{P.5}
Let $M$ be a smooth compact Riemannian manifold of dimension $n \geq 3$ and $\lambda$ and $C_1$ be positive constants. Consider a metric $g \in {\cal M}^n$, a constant $\mathcal{B} \geq 0$ and potential functions $F \in \mathcal{F}_k$ and $G \in \mathcal{G}_k$ satisfying

\[
g_{ij} \xi_i \xi_j \geq \lambda |\xi|^2,\ \ ||g_{ij}||_{C^0} \leq \lambda^{-1}
\]

\n for all $\xi \in \R^n$ and $i,j = 1, \ldots, n$ in a fixed coordinates system on $M$, and

\[
F(t) \leq C_1 |t|^{2^*}
\]

\n for all $t \in \R^k$. Then, given constants $q > 2^*$ and $p, K_1, K_2 > 0$, there exists a positive constant $C_0$, depending only on $n$, $q$, $p$, $K_1$, $K_2$, $\lambda$ and $C_1$ such that for any $\delta > 0$, any point $x_0 \in M$ and any minimizer $U \in \Lambda$ of the functional (\ref{E.7}) satisfying $J(U) \leq K_1$ and

\[
||U||_{L_k^{q}(B_g(x_0, 2\delta))} \leq K_2\, ,
\]

\n one has

\begin{equation} \label{lim4}
\sup_{B_g(x_0, \delta)} |U| \leq C_0 \delta^{-n/p}  \left( \int_{B_g(x_0, 2\delta)} |U|^{p}\; dv_g \right)^{1/p}\, .
\end{equation}

\end{propo}

\n {\bf Proof of Proposition \ref{P.5}.} Let $U \in \Lambda$ be a minimizer of $J$. Given $\varphi \in C^1(M)$, we introduce the function

\[
h_\varepsilon(t) = \frac{\mathcal{A}_0(n,F) \int_M |\nabla_g (U + t v_\varepsilon^{-1} \varphi U)|^2\; dv_g + \mathcal{B} \int_M G(x,U + t v_\varepsilon^{-1} \varphi U)\; dv_g}{\left( \int_M F(U + t v_\varepsilon^{-1} \varphi U)\; dv_g \right)^{2/2^*}}
\]

\n for $t \in (- \delta, \delta)$, where $v_\varepsilon$ is provided in the previous proof. By homogeneity, $h_\varepsilon(t)$ can be rewritten as

\[
h_\varepsilon(t) = \frac{\mathcal{A}_0(n,F) \int_M |\nabla_g (U + t v_\varepsilon^{-1} \varphi U)|^2\; dv_g + \mathcal{B} \int_M (1 + t v_\varepsilon^{-1} \varphi)^2 G(x,U)\; dv_g}{\left( \int_M (1 + t v_\varepsilon^{-1} \varphi)^{2^*} F(U)\; dv_g \right)^{2/2^*}}\, .
\]

\n Note that $h_\varepsilon$ is differentiable at $t=0$ and, moreover, $h'_\varepsilon(0) = 0$ since $t=0$ is a minimum point of $h_\varepsilon$. Then,

\[
h'_\varepsilon(0) = 2 \mathcal{A}_0(n,F) \int_M \nabla_g U \cdot \nabla_g (v_\varepsilon^{-1} \varphi U)\; dv_g + 2 \mathcal{B} \int_M v_\varepsilon^{-1} G(x,U) \varphi\; dv_g - 2 J(U) \int_M v_\varepsilon^{-1} F(U) \varphi\; dv_g = 0\, ,
\]

\n so that

\[
\mathcal{A}_0(n,F) \int_M \nabla_g U \cdot \nabla_g (v_\varepsilon^{-1} \varphi U)\; dv_g + \mathcal{B} \int_M v_\varepsilon^{-1} G(x,U) \varphi\; dv_g = J(U) \int_M v_\varepsilon^{-1} F(U) \varphi\; dv_g\, .
\]

\n On the other hand, mimicking the proof of Proposition \ref{P.4}, we get

\[
\mathcal{A}_0(n,F) \int_{M} \nabla_g v_\varepsilon \cdot \nabla_g \varphi\; dv_g \leq \sum_{i=1}^k \int_{M} \nabla_g u_i \cdot \nabla_g (u_i v_\varepsilon^{-1} \varphi)\; dv_g = \int_M \nabla_g U \cdot \nabla_g (v_\varepsilon^{-1} \varphi U)\; dv_g
\]

\n for all nonnegative function $\varphi \in C^1(M)$. Since $\mathcal{B} \geq 0$, we arrive at

\[
\mathcal{A}_0(n,F) \int_{M} \nabla_g v_\varepsilon \cdot \nabla_g \varphi\; dv_g \leq J(U) \int_M v_\varepsilon^{-1} F(U) \varphi\; dv_g\, .
\]

\n So, when $\varepsilon \rightarrow 0$, we derive

\[
\mathcal{A}_0(n,F) \int_{M} \nabla_g |U| \cdot \nabla_g \varphi\; dv_g \leq J(U) \int_M |U|^{-1} F(U) \varphi\; dv_g
\]

\n for all nonnegative function $\varphi \in C^1(M)$. The conclusion then follows as in the proof of Proposition \ref{P.4}.\bl \\

\n {\bf 2.3 An approximation scheme}. For the developing of a potential type best constants general theory as proposed here, we need to know if the space of positive homogeneous $C^1$ functions is dense in the space of positive homogeneous continuous functions. In this subsection we shall introduce a simple approximation scheme that allow us to answer affirmatively this question.

\begin{propo} \label{P.6}
Let $M$ be a compact differentiable manifold of dimension $n \geq 2$, $k \geq 1$ be an integer and $H_0 : M \times \R^k \rightarrow \R$ be a positive continuous function homogeneous of degree $p > 1$ on the second variable. Then, for any $\varepsilon > 0$, there exists a positive function $H : M \times \R^k \rightarrow \R$ smooth on the first variable and of class $C^1$ and homogeneous of degree $p$ on the second variable such that

\begin{equation} \label{E.28}
|H(x,t) - H_0(x,t)| \leq \varepsilon |t|^p
\end{equation}

\n for all $x \in M$ and $t \in \R^k$.

\end{propo}

\n {\bf Proof of Proposition \ref{P.6}.} The construction of $H$ is quite simple. First, note that the space $C^\infty_{loc}(M \times \R^k)$ is dense in the space $C^0_{loc}(M \times \R^k)$. So, let $(\tilde{H}_\alpha) \subset C^\infty(M \times \R^k)$ be a sequence of functions converging to $H_0$ in $C^0_{loc}(M \times \R^k)$. Since $H_0$ is positive on $M \times \s^{k-1}$, it follows that $(\tilde{H}_\alpha)$ is positive too on $M \times \s^{k-1}$ for $\alpha$ large. Define now the sequence $(H_\alpha)$ by

\[
H_\alpha(x,t) =
\left\{
\begin{array}{ll}
|t|^{p} \tilde{H}_\alpha(x,\frac{t}{|t|}), & {\rm for}\ x \in M \ {\rm and} \ t \neq 0\\
0, & {\rm for}\ x \in M \ {\rm and} \ t = 0
\end{array}
\right. \, .
\]

\n Clearly, the functions $H_\alpha$'s are positive and homogeneous of degree $p$ on the second variable. Moreover, they also converge to $H_0$ in $C^0(M \times \s^{k-1})$. In particular, for any $\varepsilon > 0$, there exists $\alpha_0 \geq 1$ such that $H = H_{\alpha_0}$ satisfies (\ref{E.28}). Since $p > 1$, it easily follows that $H$ is of class $C^\infty$ on the first variable and of class $C^1$ on the second one. This ends the proof. \bl\\

\begin{rem} \label{R.5} Although the proof above guarantees the existence of suitable approximations which are smooth on the first variable, only $C^0$ regularity suffices for our purposes.
\end{rem}

\begin{rem} \label{R.6} As a consequence of Proposition \ref{P.6}, for any $\varepsilon > 0$ there exists a positive continuous function $H : M \times \R^k \rightarrow \R$ of class $C^1$ and homogeneous of degree $p$ on the second variable such that

\begin{equation} \label{E.29}
(1 - \varepsilon) H_0(x,t) \leq H(x,t) \leq (1 + \varepsilon) H_0(x,t)
\end{equation}

\n on $M \times \R^k$. This inequality will be frequently used later.
\end{rem}

\n {\bf 2.4 A local geometric potential type Sobolev inequality}. In this subsection, we extend a local geometric Sobolev inequality established by Druet in \cite{D3} to the vector context. Such an inequality will play an essential role in the proofs of Theorems \ref{Teo.1}, \ref{Teo.2} and \ref{Teo.3}.

\begin{propo} \label{P.7}
Let $M$ be a compact differentiable manifold of dimension $n \geq 4$ and $x_0$ a point of $M$. Consider convergent sequences of metrics $(g_\alpha) \subset {\cal M}^n$ and of potential functions $(F_\alpha) \subset \mathcal{F}_k$ and $(G_\alpha) \subset \mathcal{G}_k$. In particular, there exist positive constants $\lambda$, $C_1$, $C_2$, $c_1$ and $c_2$, independent of $\alpha$, satisfying

\begin{equation} \label{H.1}
({g_\alpha})_{ij} \xi_i \xi_j \geq \lambda |\xi|^2,\ \ ||({g_\alpha})_{ij}||_{C^2} \leq \lambda^{-1}
\end{equation}

\n for all $\xi \in \R^n$ and $i,j = 1, \ldots, n$, in a fixed coordinates system on $M$, and

\begin{equation} \label{H.2}
c_1 |t|^{2^*} \leq F_\alpha(t) \leq C_1 |t|^{2^*},\ \ c_2 |t|^2 \leq G_\alpha(x,t) \leq C_2 |t|^2
\end{equation}

\n for all $x \in M$ and $t \in \R^k$. Then, for any $\varepsilon > 0$ there exists a radius $r_0 > 0$, depending only on $\varepsilon$, $\lambda$, $C_1$, $C_2$, $c_1$ and $c_2$, such that for any map $U \in C^1_{0,k}(B_{g_\alpha}(x_0,r_0))$,

\[
\left ( \int_M F_\alpha(U)\; dv_{g_\alpha}  \right )^{2/2^*} \leq  \mathcal{A}_0(n,F_\alpha) \int_M |\nabla_{g_\alpha} U|^2\; dv_{g_\alpha} + \mathcal{B}_\varepsilon(F_\alpha,G_\alpha,g_\alpha) \int_M G_\alpha(x,U)\; dv_{g_\alpha}\,,
\]

\n where $B_g(x_0,r_0)$ stands for the ball of radius $r_0$ centered at $x_0$ in relation to a Riemannian metric $g$, $C^1_{0,k}(B_g(x_0,r_0))$ denotes the space $C^1_0(B_g(x_0,r_0)) \times \ldots \times C^1_0(B_g(x_0,r_0))$ of compactly supported $k$-maps on $B_g(x_0,r_0)$ of class $C^1$, and

\[
\mathcal{B}_\varepsilon(F,G,g) := \frac{n-2}{4(n-1)} \frac{\mathcal{A}_0(n,F)}{m_{F,G}(x_0)} Scal_g(x_0) + \varepsilon\; ,
\]

\n is defined for potential functions $F$ and $G$ and Riemannian metrics $g$, where

\[
m_{F,G}(x) := \min \limits_{t \in X_F} G(x,t)
\]

\n with $X_F = \{ t \in \mathbb{S}^{k-1}:\; F(t) = M_F\}$.

\end{propo}

We must initially present the PDEs setting involved in the proof of this proposition and that also will be useful in the proof of other results in this work. Suppose by contradiction that the conclusion fails. In this case, there exist $\varepsilon_0 > 0$ and a sequence of positive numbers $(r_\alpha)$ with $r_\alpha \rightarrow 0$ as $\alpha \rightarrow + \infty$, such that, unless re-indexing,

\begin{equation} \label{I.1}
\lambda_\alpha : = \inf \limits_{U \in C^1_{0,k}(B_{{g_\alpha}}(x_0,r_\alpha)) \setminus \{0\}} \frac{\mathcal{A}^\alpha \int_M |\nabla_{g_\alpha} U|^2\; dv_{g_\alpha} + \mathcal{B}^\alpha \int_M G_\alpha(x,U)\; dv_{g_\alpha}}{\left( \int_M F_\alpha(U)\; dv_{g_\alpha} \right)^{2/2^*}} < 1\; ,
\end{equation}

\n where

\[
\mathcal{A}^\alpha = \mathcal{A}_0(n,F_\alpha)\; ,
\]

\[
\mathcal{B}^\alpha := \frac{n-2}{4(n-1)} \frac{\mathcal{A}_0(n,F_\alpha)}{m_{F_\alpha,G_\alpha}(x_0)} Scal_{g_\alpha}(x_0) + \varepsilon_0\; .
\]

\n Note that one of the new sequences $(g_\alpha)$, $(F_\alpha)$ or $(G_\alpha)$ can possibly be constant. In any case, let $g \in {\cal M}^n$, $F \in \mathcal{F}_k$ and $G \in \mathcal{G}_k$ be the respective limits of these sequences.

By Remark \ref{R.3} of Proposition \ref{P.2}, it follows that

\begin{equation} \label{lim1}
\mathcal{A}^\alpha \rightarrow \mathcal{A}_0(n,F)
\end{equation}

\n and

\begin{equation} \label{lim2}
\mathcal{B}^\alpha \rightarrow \frac{n-2}{4(n-1)} \frac{\mathcal{A}_0(n,F)}{m_{F,G}(x_0)} Scal_{g}(x_0) + \varepsilon_0
\end{equation}

\n as $\alpha \rightarrow +\infty$.

By Proposition \ref{P.6}, we can assume that the potential functions $F_\alpha$ and $G_\alpha$ in (\ref{I.1}) are of class $C^1$, converge to $F$ and $G$, respectively, and satisfy (\ref{H.2}) with possibly different constants.

By (\ref{I.1}) and Proposition \ref{exist}, it follows that $\lambda_\alpha$ is achieved by a map $U_\alpha = (u_{\alpha}^1, \ldots, u_{\alpha}^k) \in H^{1,2}_k(M)$ for each $\alpha > 0$. In particular, the maps $U_\alpha$ satisfy the potential systems

\begin{equation} \label{S.2}
-\mathcal{A}^\alpha \Delta_{g_\alpha} u_{\alpha}^i +\frac{1}{2}\mathcal{B}^\alpha \frac{\partial G_\alpha}{\partial t_i}(x,U_{\alpha}) = \frac{\lambda_{\alpha}}{2^*}\frac{\partial F_\alpha}{\partial t_i}(U_{\alpha}) \ \ {\rm in} \ B_{g_\alpha}(x_0,r_{\alpha})\, ,\ \ i=1, \ldots, k
\end{equation}

\n with

\[
U_{\alpha} = 0 \ \ {\rm on} \ \partial B_{g_\alpha}(x_0,r_{\alpha})
\]

\n and

\[
\int_M F_\alpha(U_{\alpha})\; dv_{g_\alpha} = 1\, .
\]

\n Moreover, the maps $U_\alpha$ are of class $C^1$, by Proposition \ref{Reg}. Although the proofs of existence and regularity of minimizers is rather classical, for completeness we include them in an appendix.

Our goal now is to study the behavior of $(U_\alpha)$ as $\alpha$ increases. By (\ref{H.2}) and H\"{o}lder's inequality, we get

\begin{eqnarray*}
\int_{B_{g_\alpha}(x_0,r_{\alpha})} |U_{\alpha}|^2\; dv_{g_\alpha} &\leq& vol_{g_\alpha}(B_{g_\alpha}(x_0,r_{\alpha}))^{2/n} \left( \int_{B_{g_\alpha}(x_0,r_{\alpha})} |U_{\alpha}|^{2^*}\; dv_{g_\alpha} \right)^{2/2^*} \\
&\leq& c_1^{-2/2^*} \left( \int_{B_{g_\alpha}(x_0,r_{\alpha})} F_\alpha(U_{\alpha})\; dv_{g_\alpha} \right)^{2/2^*} vol_{g_\alpha}(B_{g_\alpha}(x_0,r_{\alpha}))^{2/n} \\
&\leq& c_1^{-2/2^*} vol_{g_\alpha}(B_{g_\alpha}(x_0,r_{\alpha}))^{2/n} \rightarrow 0
\end{eqnarray*}

\n as $\alpha \rightarrow + \infty$, so that

\begin{equation} \label{E.3}
\lim \limits_{\alpha \rightarrow +\infty} \int_{B_{g_\alpha}(x_0,r_{\alpha})} G_\alpha(x,U_{\alpha})\; dv_{g_\alpha} = 0\; .
\end{equation}

\n We now show that

\begin{equation} \label{E.4}
\lim \limits_{\alpha \rightarrow +\infty} \lambda_\alpha = 1\; .
\end{equation}

\n Given $\varepsilon > 0$, by standard arguments and (\ref{H.1}), one easily finds a real constant $\tilde{B}_\varepsilon$, independent of $\alpha$, such that

\[
\left(\int_{M} |u|^{2^*}\; dv_{g_\alpha} \right)^{2/2^*} \leq (A_0(n) + \varepsilon) \int_{M} |\nabla_{g_\alpha} u|^2\; dv_{g_\alpha} + \tilde{B}_\varepsilon \int_{M} u^2\; dv_{g_\alpha}
\]

\n for all $u \in C^1(M)$ (alternatively, this conclusion also follows with $\varepsilon = 0$ from the work \cite{EM6} about the geometric continuity of $B_0(n,1,g)$). So, mimicking (\ref{des-2}) and using (\ref{H.2}), we obtain a constant $\tilde{\cal B}_\varepsilon$, independent of $\alpha$, such that, for any $U \in C^1_{0,k}(B_{{g_\alpha}}(x_0,r_\alpha))$,

\[
\int_{B_{g_\alpha}(x_0,r_{\alpha})} F_\alpha(U)\; dv_{g_\alpha} \leq (\mathcal{A}^\alpha + \varepsilon) \int_{B_{g_\alpha}(x_0,r_{\alpha})} |\nabla_{g_\alpha} U|^2\; dv_{g_\alpha} + \tilde{\cal B}_\varepsilon \int_{B_{g_\alpha}(x_0,r_{\alpha})} G_\alpha(x,U)\; dv_{g_\alpha}
\]

\n Taking now $U = U_\alpha$ in this inequality and using that $U_\alpha$ is a minimizer for $\lambda_\alpha$, one has

\begin{eqnarray*}
1 &\leq& ( \mathcal{A}^\alpha + \varepsilon ) (\mathcal{A}^\alpha)^{-1} \left( \lambda_\alpha - \mathcal{B}^\alpha \int_{B_{g_\alpha}(x_0,r_{\alpha})}G_\alpha(x,U_{\alpha})\; dv_{g_\alpha} \right) \\
&&+ \tilde{\cal B}_\varepsilon \int_{B_{g_\alpha}(x_0,r_{\alpha})} G_\alpha(x,U_{\alpha})\; dv_{g_\alpha}\,,
\end{eqnarray*}

\n Letting then $\alpha \rightarrow +\infty$ and after $\varepsilon \rightarrow 0$, one concludes from (\ref{lim1}), (\ref{lim2}) and (\ref{E.3}) that

\[
\liminf \limits_{\alpha \rightarrow +\infty} \lambda_{\alpha} \geq 1\; .
\]

\n So, since $\lambda_\alpha < 1$ for all $\alpha$, the assertion (\ref{E.4}) follows readily. Consequently, from the equality

\[
\mathcal{A}^\alpha \int_{B_{g_\alpha}(x_0,r_{\alpha})} |\nabla_{g_\alpha} U_{\alpha}|^2\; dv_{g_\alpha} + \mathcal{B}^\alpha \int_{B_g(x_0,r_{\alpha})} G_\alpha(x,U_{\alpha})\; dv_{g_\alpha} = \lambda_{\alpha},
\]

\n we derive

\begin{equation} \label{E.13}
\lim \limits_{\alpha \rightarrow +\infty} \int_{B_{g_\alpha}(x_0,r_{\alpha})} |\nabla_{g_\alpha} U_{\alpha}|^2\; dv_{g_\alpha} = \mathcal{A}_0(n,F)^{-1}\, .
\end{equation}

In order to establish some fine properties to $(U_\alpha)$, we will re-normalize this sequence as follows. Let $x_{\alpha} \in B_{g_\alpha}(x_0,r_{\alpha})$ be a maximum point of $|U_{\alpha}(x)|$ and set

\[
\mu_\alpha = |U_{\alpha}(x_{\alpha})|^{-2^*/n}\, .
\]

\n By (\ref{H.2}), one has

\[
1 = \int_{B_{g_\alpha}(x_0,r_{\alpha})}F_\alpha(U_{\alpha})\; dv_{g_\alpha} \leq C_1 \int_{B_{g_\alpha}(x_0,r_{\alpha})} |U_{\alpha}|^{2^*}\; dv_g \leq C_1 vol_{g_\alpha}(B_{g_\alpha}(x_0,r_{\alpha})) \mu_{\alpha}^{-n}\; ,
\]

\n so that

\begin{equation} \label{E.17}
\mu_{\alpha} \rightarrow 0
\end{equation}

\n as $\alpha \rightarrow +\infty$. Introduce, for each $\alpha > 0$, the open set in $\mathbb{R}^n$

\[
\Omega_{\alpha} = \mu_{\alpha}^{-1} \exp_{x_{\alpha}}^{-1} \left
(B_{g_\alpha}(x_0,r_{\alpha}) \right )\; ,
\]

\n the metric

\[
h_{\alpha}(x) = \left( \exp_{x_{\alpha}}^*g_\alpha \right)(\mu_{\alpha}x)\; ,
\]

\n and the map

\[
V_{\alpha}(x) = \left\{
\begin{array}{ll}
\mu_{\alpha}^{n/2^*}U_{\alpha}\left
(\exp_{x_{\alpha}}(\mu _{\alpha}x)  \right ) \ \ {\rm if}\ x \in \Omega_{\alpha},\\
0\ \ {\rm if}\ x \in \mathbb{R}^n \setminus \Omega_{\alpha}
\end{array}
\right.\, ,
\]

\n where $\exp_{x_{\alpha}}$ denotes the exponential chart with respect to the metric $g_\alpha$ centered at $x_\alpha$.

\n Clearly, $V_{\alpha} = (v_{\alpha}^1, \ldots, v_{\alpha}^k)$ satisfies

\begin{equation}\label{NE1}
-\mathcal{A}^\alpha \Delta_{h_{\alpha}} v_{\alpha}^i + \frac{\mu_{\alpha}^2}{2} \mathcal{B}^\alpha \frac{\partial
G_\alpha}{\partial t_i}(\exp_{x_{\alpha}}(\mu_{\alpha}x), V_{\alpha}) = \frac{\lambda_{\alpha}}{2^*} \frac{\partial F_\alpha}{\partial t_i}(V_{\alpha}) \ \ {\rm in} \ \Omega_\alpha
\end{equation}

\n and

\begin{equation} \label{E.20}
\int_{\Omega_{\alpha}} F_\alpha(V_{\alpha})\; dv_{h_{\alpha}} = 1\; .
\end{equation}

An important step is to deduce the convergence of the sequence $(V_\alpha)$ in ${\cal D}^{1,2}_k(\R^n)$ by using the convergence of $(g_\alpha)$, $(F_\alpha)$ and $(G_\alpha)$ in, respectively, ${\cal M}^n$, $\mathcal{F}_k$ and $\mathcal{G}_k$.

\begin{lema} \label{L.1}
The sequence $(V_\alpha)$ converges to some map $V_0$ in ${\cal D}^{1,2}_k(\R^n)$. Moreover, we have $V_0 = t_0 v_0$, where $t_0 \in \s^{k-1}$ is a maximum point of $F$ and $v_0$ is an extremal function to $A_0(n)$.
\end{lema}

\n {\bf Proof of Lemma \ref{L.1}.} At first, through a simple change of variable, we have

\[
\int_{\Omega_{\alpha}} |\nabla_{h_{\alpha}} V_\alpha|^2\; dv_{h_{\alpha}} = \int_{B_{g_\alpha}(x_0,r_{\alpha})} |\nabla_{g_\alpha} U_{\alpha}|^2\; dv_{g_\alpha}\, ,
\]

\n so that

\[
\lim \limits_{\alpha \rightarrow +\infty} \int_{\Omega_{\alpha}} |\nabla_{h_{\alpha}} V_\alpha|^2\; dv_{h_{\alpha}} = \mathcal{A}_0(n,F)^{-1}\, .
\]

\n Note also that $|V_{\alpha}| \leq 1$ in $\Omega_\alpha$, $|V_{\alpha}(0)| = 1$ and the metrics $h_{\alpha}$ converges to the Euclidean metric $\xi$ in $C^2_{loc}(\R^n)$. Regarding the Cartan expansion of $h_\alpha$, we find a constant $C > 0$, independent of $\alpha$, such that

\[
(1 - C \mu _{\alpha}^2) dx \leq dv_{h_{\alpha}} \leq (1 + C \mu_{\alpha}^2) dx
\]

\n and

\[
(1 - C \mu _{\alpha}^2) |\nabla V_{\alpha}|^2 \leq |\nabla_{h_{\alpha}} V_{\alpha}|^2 \leq (1 + C \mu_{\alpha}^2 ) |\nabla V_{\alpha}|^2\,.
\]

\n Clearly, these inequalities imply that

\begin{equation} \label{E.5}
\lim \limits_{\alpha \rightarrow +\infty} \int_{\R^n} F(V_{\alpha})\; dx = 1
\end{equation}

\n and

\begin{equation} \label{E.6}
\lim \limits_{\alpha \rightarrow +\infty} \int_{\R^n} |\nabla V_\alpha|^2\; dx = \mathcal{A}_0(n,F)^{-1}\, .
\end{equation}

\n Therefore, the sequence $(V_{\alpha})$ is bounded in ${\cal D}^{1,2}_k(\R^n)$, so that it converges weakly to $V_0$ in ${\cal D}^{1,2}_k(\R^n)$. Evoking now the concentration-compactness principle of Lions to the measure $F_\alpha(V_\alpha) dv_{h_{\alpha}}$, three situations may occur: compactness, vanishing or dichotomy, see, for instance, Lemma 4.3 of \cite{St}. As for vanishing, by applying the De Giorgi-Nash-Moser estimate for potential type elliptic systems given in Proposition \ref{P.4} (this point requires only the convergence of $(F_\alpha)$ and $(G_\alpha)$ in, respectively, $\mathcal{F}_k$ and $\mathcal{G}_k$), we obtain a constant $C_0 > 0$, independent of $\alpha$, such that

\begin{equation} \label{E.9}
1 = \sup \limits_{B(0,1) \cap \Omega_{\alpha}} |V_{\alpha}|^2 \leq C_0 \left( \int_{B(0,2) \cap \Omega_{\alpha}} F_\alpha(V_{\alpha}) \; dv_{h_{\alpha}} \right)^{1/2^*}
\end{equation}

\[
\leq C_0 \left( \int_{B(0,2)} F_\alpha(V_{\alpha}) \; dv_{h_{\alpha}} \right)^{1/2^*}
\]

\n and vanishing can not happen. Dichotomy is classically forbidden by (\ref{E.5}) and (\ref{E.6}) following some ideas of P. L. Lions. We give some details of this last statement.
Setting $\rho_{\alpha}=|\nabla V_{\alpha}|^2+F(V_{\alpha})$, we have $\rho_{\alpha}\geq0$ and

\[
\lim \limits_{{\alpha}\rightarrow +\infty}\int_{\R^n}\rho_{\alpha}\; dx=(\mathcal{A}_0(n,F)^{-1}+1)>0\,.
\]

\n Suppose that dichotomy is underway: there exists $\beta \in ]0,\mathcal{A}_0(n,F)^{-1}+1[$ such that for any $\varepsilon>0$ there are $\alpha_0 \in \mathbb{N}$ and $\rho_{\alpha}^1, \rho_{\alpha}^2\in L^1(\R^n)$  satisfying

\[
||\rho_{\alpha}-\rho_{\alpha}^1-\rho_{\alpha}^2||_{L^1(\R^n)} < \varepsilon\, ,
\]

\[
\left|\int_{\R^n}\rho_{\alpha}^idx-\beta \right| < \varepsilon
\]

\n and

\[
\lim \limits_{{\alpha}\rightarrow +\infty}dist(supp\,\rho_{\alpha}^1,supp\,\rho_{\alpha}^2) = +\infty
\]

\n for all ${\alpha} \geq {\alpha}_0$. Consider the following sets

\[
B_1 = \{x\in \R^n:\,\, dist(x,supp\,\rho_{\alpha}^1) <1\}
\]

\n and

\[
B_2=\{x\in \R^n:\,\, dist(x,supp\,\rho_{\alpha}^2) <1\}\, .
\]

\n Consider also the following smooth cutoff functions

\[
\xi_{\alpha}=\left\{
       \begin{array}{cc}
         1, & \mbox{in $D_1$} \\
         0, & \mbox{in $\R^n\setminus B_1$} \\
       \end{array}
     \right.
\]

\n and

\[
\eta_{\alpha}=\left\{
       \begin{array}{cc}
         1, & \mbox{in $D_2$} \\
         0, & \mbox{in $\R^n\setminus B_2$} \\
       \end{array}
     \right.
\]

\n such that $0\leq\xi_{\alpha},\eta_{\alpha}\leq1$ and $|\nabla \xi_{\alpha}|, |\nabla \eta_{\alpha}|\leq 2$, where $D_1= supp\,\rho_{\alpha}^1$ and $D_2= supp\,\rho_{\alpha}^2$. Define

\[
\tilde{\rho}_{\alpha}^1=\rho_{\alpha}\xi_{\alpha}
\]

\n and

\[
\tilde{\rho}_{\alpha}^2=\rho_{\alpha}\eta_{\alpha}\,.
\]

\n Note that

\[
||\rho_{\alpha}-\tilde{\rho}_{\alpha}^1-\tilde{\rho}_{\alpha}^2||_{L^1(\R^n)}<\varepsilon\,,
\]

\n and

\[
\lim \limits_{{\alpha}\rightarrow +\infty}dist(supp\,\tilde{\rho}_{\alpha}^1,supp\,\tilde{\rho}_{\alpha}^2)=+\infty\,.
\]

\n Set

\[
V_{\alpha}^1=V_{\alpha}\xi_{\alpha}\,,
\]

\[
V_{\alpha}^2=V_{\alpha}\eta_{\alpha}
\]

\n and

\[
E_{\alpha}=\R^n\setminus\{x\in \R^n:\,\,\xi_{\alpha}(x)=1\,\,\mbox{or}\,\,\eta_{\alpha}(x)=1  \}\,.
\]

\n Then, for any $\alpha$ large,

\[
\left|\int_{\R^n}|\nabla V_{\alpha}|^2\; dx-\int_{\R^n}|\nabla V_{\alpha}^1|^2\; dx-\int_{\R^n}|\nabla V_{\alpha}^2|^2\; dx\right|\leq C\int_{E_{\alpha}}|\nabla V_{\alpha}|^2\; dx\leq C\varepsilon
\]

\n and

\[
\left|\int_{\R^n}F(V_{\alpha})\; dx-\int_{\R^n}F(V_{\alpha}^1)\; dx-\int_{\R^n}F(V_{\alpha}^2)\; dx\right|\leq C\int_{E_{\alpha}}F(V_{\alpha})\; dx\leq C\varepsilon\,.
\]

\n Hence,

\[
\mathcal{A}_0(n,F)^{-1}+\varepsilon>\int_{\R^n}|\nabla V_{\alpha}|^2\; dx\geq
\int_{\R^n}|\nabla V_{\alpha}^1|^2\; dx+\int_{\R^n}|\nabla V_{\alpha}^2|^2\; dx- C\varepsilon\,.
\]

\n Set

\[
\lambda^i_{\alpha}=\int_{\R^n}F(V_{\alpha}^i)\; dx\,, \quad \, \mbox{i=1,\,2}\,.
\]

\n Clearly, $0\leq\lambda^i_{\alpha}\leq1$. Thus, there exists a subsequence $(\lambda^i_{{\alpha}_l})$ such that

\[
\lim \limits_{\alpha_l\rightarrow +\infty}\lambda_{\alpha_l}^i=\lambda_i
\]

\n and

\[
\lambda_1+\lambda_2=1\,.
\]

\n If $\lambda_1=0$, then $\lambda_2=1$. Consequently,

\[
\beta -\varepsilon<\int_{\R^n}\rho_{\alpha}^1\; dx=\int_{\R^n}\left(F(V_{\alpha}^1)+|\nabla V_{\alpha}^1|^2\right)\; dx
\]

\[
=\int_{\R^n}|\nabla V_{\alpha}^1|^2\; dx+o(1)\,,
\]

\n since

\[
\lim \limits_{\alpha\rightarrow +\infty}\int_{\R^n}F(V_{\alpha}^1)\; dx=\lambda_1=0\,.
\]

\n Thus,

\[
0<\beta-\varepsilon\leq\int_{\R^n}|\nabla V_{\alpha}^1|^2\; dx\,.
\]

\n and

\[
\mathcal{A}_0(n,F)^{-1}+\varepsilon>\int_{\R^n}|\nabla V_{\alpha}|^2\; dx\geq \int_{\R^n}|\nabla V_{\alpha}^1|^2\; dx+\int_{\R^n}|\nabla V_{\alpha}^2|^2\; dx
-C\varepsilon
\]

\[
\geq (\beta-\varepsilon)+\int_{\R^n}|\nabla V_{\alpha}^2|^2\; dx-C\varepsilon\,.
\]

\n Letting $\alpha\rightarrow+\infty$ and $\varepsilon\rightarrow 0$, we deduce that $\beta \leq0$, contradicting that $\beta>0$. Therefore, $0<\lambda_1<1$ and, similarly, $0<\lambda_2<1$. From what has been done previously, we have

\[
\int_{\R^n}|\nabla V_{\alpha}|^2\; dx\geq \int_{\R^n}|\nabla V_{\alpha}^1|^2\; dx+\int_{\R^n}|\nabla V_{\alpha}^2|^2\; dx\,.
\]

\n Then,

\[
\mathcal{A}_0(n,F)^{-1}\geq I_{\lambda_1}+ I_{\lambda_2}\,,
\]

\n where

\[
I_{\lambda_i}=\inf \left\{ \int_{{\R^n}}|\nabla V|^2\; dx\,:\,\,\,\int_{\R^n}F(V)\; dx=\lambda_i \right\}\,,\quad i=1,\,2\,.
\]

\n But, this contradicts

\[
\mathcal{A}_0(n,F)^{-1}< I_{\lambda_1}+ I_{\lambda_2}\,,
\]

\n since

\[
\mathcal{A}_0(n,F)^{-1}=I_{1}=\inf \left\{ \int_{{\R^n}}|\nabla V|^2\; dx\,:\,\,\,\int_{\R^n}F(V)\; dx=1 \right\}\,,
\]

\[
I_{\lambda_i}=\lambda_i^{\frac{2}{2^*}}I_1\,,\quad\,i=1,\,2\,,
\]

\n and

\begin{eqnarray*}
I_{\lambda_1}+I_{\lambda_2}&=&\left(\lambda_1^{\frac{2}{2^*}}+\lambda_{2}^{\frac{2}{2^*}}\right)I_1\\
&>&\left(\lambda_1+\lambda_2  \right)^{\frac{2}{2^*}}I_1\\
&=&I_1\,.
\end{eqnarray*}

\n Therefore, compactness holds in the following sense: there exists a sequence $(y_\alpha) \subset \R^n$ such that for any $\varepsilon > 0$ there is a radius $R > 0$ with the property that

\begin{equation} \label{E.10}
\int_{B(y_\alpha,R)} F_\alpha(V_{\alpha}) \; dv_{h_{\alpha}} \geq 1 - \varepsilon
\end{equation}

\n for all $\alpha$. We claim that the sequence $(y_\alpha)$ is bounded. Otherwise, we get $B(0,2) \cap B(y_\alpha,R) = \emptyset$ for all $\alpha$ large. In particular,

\[
\int_{B(y_\alpha,R)} F_\alpha(V_{\alpha}) \; dv_{h_{\alpha}} \leq \int_{\R^n \setminus B(0,2)} F_\alpha(V_{\alpha}) \; dv_{h_{\alpha}}
\]

\n for all $\alpha$ large, so that

\[
1 - \varepsilon \leq \liminf  \limits_{\alpha \rightarrow + \infty} \int_{\R^n \setminus B(0,2)} F_\alpha(V_{\alpha}) \; dv_{h_{\alpha}}
\]

\n for any $\varepsilon > 0$. Of course, this implies that

\[
\lim \limits_{\alpha \rightarrow + \infty} \int_{\R^n \setminus B(0,2)} F_\alpha(V_{\alpha}) \; dv_{h_{\alpha}} = 1\, ,
\]

\n or equivalently,

\[
\lim \limits_{\alpha \rightarrow + \infty} \int_{B(0,2)} F_\alpha(V_{\alpha})\; dv_{h_{\alpha}} = 0\, .
\]

\n By (\ref{E.9}), we then derive a contradiction. Assume now that the sequence $(y_\alpha)$ converges to some point $y_0 \in \R^n$. Consider the decomposition

\[
\int_{B(y_\alpha,R)} F_\alpha(V_{\alpha}) \; dv_{h_{\alpha}} = \int_{B(y_\alpha,R) \setminus B(y_0,R)} F_\alpha(V_{\alpha}) \; dv_{h_{\alpha}} - \int_{B(y_0,R) \setminus B(y_\alpha,R)} F_\alpha(V_{\alpha}) \; dv_{h_{\alpha}}
\]

\[
+ \int_{B(y_0,R)} F_\alpha(V_{\alpha}) \; dv_{h_{\alpha}}\, .
\]

\n Since $|V_\alpha|$ is bounded in $L^\infty(\R^n)$ and $h_{\alpha}$ and $F_\alpha$ satisfy, respectively, (\ref{H.1}) and (\ref{H.2}), then the first two integrals on the right-hand side of the equality above converge to $0$. Moreover, the dominated convergence theorem produces

\[
\lim \limits_{\alpha \rightarrow + \infty} \int_{B(y_0,R)} F_\alpha(V_{\alpha}) \; dv_{h_{\alpha}} = \int_{B(y_0,R)} F(V_0) \; dx\, .
\]

\n Therefore, by (\ref{E.10}),

\[
\int_{\R^n} F(V_0) \; dx \geq \int_{B(y_0,R)} F(V_0) \; dx \geq 1 - \varepsilon
\]

\n for any $\varepsilon > 0$, so that

\[
\int_{\R^n} F(V_0) \; dx = 1\, .
\]

\n Clearly, this leads to

\[
\int_{\R^n} |\nabla V_0|^2\; dx \geq \mathcal{A}_0(n,F)^{-1}\, .
\]

\n On the other hand, using (\ref{E.6}) and the convexity of norms, we find

\[
\int_{\R^n} |\nabla V_0|^2\; dx \leq \lim \limits_{\alpha \rightarrow +\infty} \int_{\R^n} |\nabla V_\alpha|^2\; dx = \mathcal{A}_0(n,F)^{-1}\, .
\]

\n Thus,

\[
\int_{\R^n} |\nabla V_0|^2\; dx = \mathcal{A}_0(n,F)^{-1}\, ,
\]

\n so that the sequence $(V_\alpha)$ converges strongly to $V_0$ in ${\cal D}^{1,2}_k(\R^n)$. The rest of the proof then follow from the part (b) of Proposition \ref{P.1}.\bl \\

\begin{rem} \label{R.7} Thanks to Lemma \ref{L.1} and to the De Giorgi-Nash-Moser estimate (\ref{lim3}), we then have

\[
\lim \limits_{R \rightarrow + \infty} \left( \lim \limits_{\alpha \rightarrow + \infty} \sup_{\Omega_\alpha \setminus B(0,R)} |V_\alpha| \right) = 0\, .
\]
\end{rem}

\begin{lema} \label{L.2}
There exists a constant $C_1 > 0$ such that

\[
d_{g_\alpha}(x, x_\alpha)^{n/2^*} |U_{\alpha}(x)| \leq C_1
\]

\n for all $x \in B_{g_\alpha}(x_0,r_{\alpha})$ and $\alpha > 0$ large, where $d_{g_\alpha}$ stands for the distance with respect to $g_\alpha$.

\end{lema}

\n {\bf Proof of Lemma \ref{L.2}.} We proceed by contradiction. Let $y_\alpha \in B_{g_\alpha}(x_0,r_{\alpha})$ be a maximum point of the function

\[
u_\alpha(x) = d_{g_\alpha}(x, x_\alpha)^{n/2^*} |U_{\alpha}(x)|\, ,
\]

\n and suppose that $u_\alpha(y_\alpha)$ blows up as $\alpha \rightarrow + \infty$. This implies that $|U_{\alpha}(y_\alpha)|$ blows up too. Clearly, the sequence $(y_\alpha)$ converges to $x_0$. Extend the map $U_\alpha$ to be zero outside $B_{g_\alpha}(x_0,r_{\alpha})$.

Let $\exp_{y_{\alpha}}$ be the exponential chart centered at $y_{\alpha}$ with respect to the metric $g_\alpha$. Choose $\delta > 0$ small such that the map $\exp_{y_{\alpha}}$ is a diffeomorphism from $B(0,2\delta) \subset \R^n$ onto $B_{g_\alpha}(y_{\alpha},2\delta)$ for all $\alpha > 0$ large. Consider now the open set in $\R^n$

\[
\tilde{\Omega}_\alpha = \nu_\alpha^{-1} \exp_{y_\alpha}^{-1} (B_{g_\alpha}(x_\alpha, \delta))
\]

\n for all $\alpha > 0$, where

\[
\nu_{\alpha} = |U_{\alpha}(y_{\alpha})|^{-2^*/n}\, .
\]

\n Introduce then the following metric and map on $\tilde{\Omega}_\alpha$:

\[
\tilde{h}_{\alpha}(x) = (\exp_{y_{\alpha}}^*g_{\alpha})(\nu_{\alpha}x)
\]

\n and

\[
\tilde{V}_{\alpha}(x) = \nu_{\alpha}^{n/2^*} U_{\alpha} \left( \exp_{y_{\alpha}}(\nu_{\alpha}x) \right)\, .
\]

\n Note that $(\tilde{h}_{\alpha})$ converges to $\xi$ in $C^2_{loc}(\R^n)$. We claim that the sequence $(\tilde{V}_{\alpha})$ is uniformly bounded on $B(0,2)$. In fact, for any $x \in B(0,2)$,

\[
d_{g_\alpha} \left( x_{\alpha}, \exp_{y_{\alpha}} ( \nu_{\alpha}x ) \right) \geq d_{g_\alpha}(x_{\alpha}, y_{\alpha}) - 2 \nu_{\alpha} =  (1 - 2 u_{\alpha}(y_{\alpha})^{-2^*/n} ) d_{g_\alpha}(x_{\alpha}, y_{\alpha})\, ,
\]

\n so that

\begin{equation}\label{E.45}
d_{g_\alpha}(x_{\alpha}, \exp_{y_{\alpha}}(\nu_{\alpha} x) ) \geq
\frac{1}{2} d_{g_\alpha}(x_{\alpha},y_{\alpha})
\end{equation}

\n for $\alpha > 0$ large. Then, from (\ref{E.45}), we get

\[
|\tilde{V}_{\alpha}(x)| = \nu_{\alpha}^{n/2^*} |U_{\alpha} \left( \exp_{y_{\alpha}}(\nu_{\alpha}x) \right)| = \nu_{\alpha}^{n/2^*} d_{g_\alpha}(x_{\alpha}, \exp_{y_{\alpha}}(\nu_{\alpha} x) )^{-n/2^*} u_\alpha(\exp_{y_{\alpha}}(\nu_{\alpha} x))
\]

\[
\leq 2^{n/2^*} \nu_{\alpha}^{n/2^*} d_{g_\alpha}(x_{\alpha},y_{\alpha})^{-n/2^*} u_\alpha(y_{\alpha}) \leq 2^{n/2^*}\, .
\]

\n In conclusion, we have proved that

\begin{equation}\label{E.46}
\sup \limits_{x \in B(0,2)} |\tilde{V}_{\alpha}(x)| \leq 2^{n/2^*}
\end{equation}

\n for all $\alpha > 0$ large. On the other hand, the maps $\tilde{V}_\alpha = (\tilde{v}_{\alpha}^1, \ldots, \tilde{v}_{\alpha}^k)$ satisfy the systems

\[
- \mathcal{A}^\alpha \Delta_{\tilde{h}_{\alpha}} \tilde{v}_{\alpha}^i + \frac{\nu_{\alpha}^2}{2} \mathcal{B}^\alpha \frac{\partial G_{\alpha}}{\partial t_i} (\exp_{y_{\alpha}}(\nu_{\alpha} x), \tilde{V}_{\alpha}) = \frac{\lambda_{\alpha}}{2^*} \frac{\partial F_{\alpha}}{\partial t_i} (\tilde{V}_{\alpha}) \ \ {\rm in} \ \Theta_\alpha
\]

\n and the boundary conditions

\[
\tilde{V}_{\alpha} = 0 \ \ {\rm on} \ \partial \Theta_\alpha\, ,
\]

\n where $0 \in \Theta_\alpha = \nu_\alpha^{-1} \exp_{y_\alpha}^{-1} (B_{g_\alpha}(x_0, r_\alpha)) \subset \Omega_\alpha$. So, thanks to (\ref{E.46}), we can apply Proposition \ref{P.4} to these systems, so that there exists a constant $C_0 > 0$ such that, for $\alpha > 0$ large,

\begin{equation} \label{E.11}
1 = |\tilde{V}_\alpha(0)| \leq \sup_{B(0,1) \cap \Theta_\alpha} |\tilde{V}_\alpha| \leq C_0 \int_{B(0,2) \cap \Theta_\alpha} |\tilde{V}_\alpha|^{2^*} \; dv_{\tilde{h}_{\alpha}} \leq C_0 \int_{B_{\tilde{h}_{\alpha}}(0,3)} |\tilde{V}_\alpha|^{2^*} \; dv_{\tilde{h}_{\alpha}}
\end{equation}

\[
= C_0 \int_{B_{g_\alpha}(y_{\alpha},3 \nu_\alpha)} |U_{\alpha}|^{2^*}\; dv_{g_\alpha}\, .
\]

\n By Lemma \ref{L.1}, we easily derive an absurd if we are able to verify that

\begin{equation} \label{E.12}
B_{g_\alpha}(y_{\alpha},3 \nu_\alpha) \cap B_{g_\alpha}(x_{\alpha}, R \mu_\alpha) = \emptyset
\end{equation}

\n for $\alpha > 0$ large. In fact, plugging (\ref{E.12}) into the inequality (\ref{E.11}), we produce the following contradiction:

\[
1 \leq C_0 \int_{B_{g_\alpha}(y_{\alpha},3 \nu_\alpha)} |U_{\alpha}|^{2^*}\; dv_{g_\alpha} \leq C_0 \int_{M \setminus B_{g_\alpha}(x_{\alpha}, R \mu_\alpha)} |U_{\alpha}|^{2^*}\; dv_{g_\alpha} = C_0 \int_{\R^n \setminus B_{h_\alpha}(0,R)} |V_{\alpha}|^{2^*}\; dv_{h_\alpha}
\]

\[
\leq C_1 \int_{\R^n \setminus B(0,\frac{R}{2})} |V_0|^{2^*}\; dx
\]

\n for any fixed $R > 0$ and $\alpha > 0$ large. Finally, suppose by contradiction that (\ref{E.12}) is false for infinitely many $\alpha$. Then, for these same indexes $\alpha$, we can write

\[
u_\alpha(y_\alpha)^{2^*/n} = d_{g_\alpha}(x_{\alpha},y_{\alpha}) |U_\alpha(y_\alpha)|^{2^*/n} \leq 3 + R |U_\alpha(y_\alpha)|^{2^*/n} \mu_\alpha \leq 3 + R\, ,
\]

\n which clearly contradicts the fact that $(u_\alpha(y_\alpha))$ blows up as $\alpha \rightarrow + \infty$. \bl\\

\begin{rem} \label{R.8} In the same way, Lemma \ref{L.2} combined with (\ref{E.11}) lead to

\begin{equation} \label{E.13}
\lim \limits_{R \rightarrow + \infty} \left( \lim \limits_{\alpha \rightarrow + \infty} \sup_{\Omega_\alpha \setminus B(0,R)} |V_\alpha|\; |x|^{n/2^*} \right) = 0\, .
\end{equation}
\end{rem}

\begin{lema} \label{L.2.1}
Let $s > 0$ be a small number. There exists a constant $C_2 > 0$, independent of $\alpha$, such that

\begin{equation}\label{E.14}
|V_{\alpha}| \leq C_2 |x|^{-n+2+s}
\end{equation}

\n for all $x \in \Omega_{\alpha}$ and $\alpha$ large.

\end{lema}

\n {\bf Proof of Lemma \ref{L.2.1}.} Let $L_{\alpha}$ be the second order elliptic operator

\[
L_{\alpha} u = -\mathcal{A}^\alpha \Delta_{g_{\alpha}} u +2\mathcal{B}^\alpha \mu_{\alpha}^2 m_{G_\alpha} u - 2 \lambda_{\alpha} M_{F_\alpha} |u|^{2^*-2}u\, .
\]

\n Set $\varphi(x) = \left( \frac{R}{|x|} \right)^{n-2-s}$ with $s, R>0$. Thanks to (\ref{E.3}), (\ref{E.4}), (\ref{E.13}) and (\ref{E.12}), easy computations furnish

\[
L_{\alpha} \varphi \geq 0 \ \ {\rm in}\ \ \Omega_{\alpha} \setminus B(0,R)\, .
\]

\n In addition, mimicking the proof of (\ref{R}), we have

\[
L_{\alpha} |V_{\alpha}| \leq 0\, .
\]

\n Now we may apply a maximum principle due to Aubin and Li (see Lemma 3.4 of \cite{AuLi}), to get a constant $C_2 > 0$ such that

\[
|V_{\alpha}| \leq C_2 |x|^{-n+2+s} \ \ {\rm in}\ \ \Omega_{\alpha} \setminus B(0,R)
\]

\n for $\alpha$ large. Since this inequality clearly holds in $B(0,R)$, we conclude its validity in $\Omega_{\alpha}$. \bl \\

\n {\bf Proof of Proposition \ref{P.7}.} In order to become simpler some notations below, we use the shorthand

\[
G_\alpha(V_{\alpha}) := G_\alpha(\exp_{x_{\alpha}}(\mu_\alpha x),V_{\alpha})\, .
\]

\n Our aim now is to find asymptotic expansions for the integrals

\[
\int_{\Omega_{\alpha}} F_\alpha(V_{\alpha})\; dx
\]

\n and

\[
\int_{\Omega _{\alpha}} |\nabla V_{\alpha}|^2\; dx
\]

\n as $\alpha \rightarrow +\infty$. Such expansions will follow from the following limits:

\begin{equation}\label{E.15}
\lim \limits_{\alpha \rightarrow + \infty} \left( \int_{\Omega_{\alpha}} F_\alpha(V_{\alpha}) Ric_{g_\alpha}(x_{\alpha})_{ij}x^ix^j\; dv_{h_{\alpha}} \right) = \int_{\mathbb{R}^n} F(V_0) Ric_{g}(x_0)_{ij}x^ix^j\; dx\, .
\end{equation}

\begin{equation} \label{E.16}
\lim \limits_{\alpha \rightarrow + \infty} \left( \int _{\Omega _{\alpha}}|\nabla _{h_{\alpha}}V_{\alpha}|^2 Ric_{g_\alpha}(x_{\alpha})_{ij}x^ix^j \; dv_{h_{\alpha}} \right) = \int_{\R^n}|\nabla V_0|^2 Ric_{g}(x_0)_{ij}x^ix^j\; dx\, ,
\end{equation}

\n where $Ric_g$ denotes the Ricci curvature tensor of the metric $g$. The limit (\ref{E.15}) directly follows from the dominated convergence theorem with the aid of Lemmas \ref{L.1} and \ref{L.2}. In order to show (\ref{E.16}), we first claim that

\begin{equation} \label{E.21}
\lim \limits_{R \rightarrow + \infty} \left( \lim \limits_{\alpha \rightarrow + \infty} \int_{\Omega_{\alpha} \setminus B(0,R)} |\nabla_{h_{\alpha}}V_{\alpha}|^2 |x|^2\; dv_{h_{\alpha}} \right) = 0\, .
\end{equation}

\n In fact, let $\eta_R$ be a smooth function such that $\eta_R = 1$ in $\R^n \setminus B(0,R)$ and $\eta_R = 0$ in $B(0,R/2)$. Taking $\eta_R^2 v_\alpha^i |x|^2$ as a test function in the $i$th equation of the system (\ref{NE1}), integrating by parts and adding all equations, we get

\[
\mathcal{A}^\alpha \sum_{i} \int_{\Omega_\alpha} \nabla_{h_{\alpha}} v_\alpha^i \cdot \nabla_{h_{\alpha}} (\eta_R^2 v_\alpha^i |x|^2)\; dv_{h_{\alpha}} +  \mathcal{B}^\alpha \mu_{\alpha}^2 \int_{\Omega_\alpha} \eta_R^2 G_\alpha(V_{\alpha}) |x|^2\; dv_{h_{\alpha}}
\]

\[
= \lambda_{\alpha} \int_{\Omega_\alpha} \eta_R^2 F_\alpha(V_{\alpha}) |x|^2\; dv_{h_{\alpha}}\, .
\]

\n So, thanks to (\ref{H.1}), (\ref{H.2}), (\ref{lim1}), (\ref{lim2}), (\ref{E.4}), (\ref{E.17}) and Young's inequality, there exists a constant $C_0 > 0$ such that, for any $\alpha$ large,

\[
\int_{\Omega_{\alpha}} \eta_R^2 |\nabla_{h_{\alpha}} V_{\alpha}|^2 |x|^2\; dv_{h_{\alpha}} \leq C_0 \left( \int_{\Omega_\alpha} \eta_R^2 |V_{\alpha}|^{2^*} |x|^2\; dx + \mu_{\alpha}^2 \int_{\Omega_\alpha} \eta_R^2 |V_{\alpha}|^2 |x|^2\; dx \right.
\]

\[
\left. + \int_{\Omega_\alpha} \eta_R^2 |V_{\alpha}|^2\; dx  \right)\, .
\]

\n But this implies that

\[
\int_{\Omega_{\alpha} \setminus B(0,R)} |\nabla_{h_{\alpha}} V_{\alpha}|^2 |x|^2\; dv_{h_{\alpha}} \leq C_0 \left( \int_{\R^n \setminus B(0,R/2)} |V_{\alpha}|^{2^*} |x|^2\; dx + \mu_{\alpha}^2 \int_{\R^n \setminus B(0,R/2)} |V_{\alpha}|^2 |x|^2\; dx \right.
\]

\[
\left.+ \int_{\R^n \setminus B(0,R/2)} |V_{\alpha}|^2\; dx  \right)\, .
\]

\n Using that $n \geq 5$ and (\ref{E.14}) on the right hand-side above and after letting $R \rightarrow + \infty$, we arrive at (\ref{E.21}). Finally, (\ref{E.16}) follows from (\ref{E.21}) and Lemma \ref{L.1}. Since

\[
\int_{\R^n} f x^ix^j\; dx = \frac{\delta_{ij}}{n} \int_{\R^n} f |x|^2\; dx
\]

\n for radial functions $f$, (\ref{E.15}) and (\ref{E.16}) immediately yield

\begin{equation} \label{E.18}
\lim \limits_{\alpha \rightarrow + \infty} \left( \int_{\Omega_{\alpha}} F_\alpha(V_{\alpha}) Ric_{g_\alpha}(x_{\alpha})_{ij}x^ix^j\; dv_{h_{\alpha}} \right) = \frac{Scal_g(x_0)}{n} \int_{\R^n} F(V_0) |x|^2\; dx
\end{equation}

\begin{equation}\label{E.19}
\lim \limits_{\alpha \rightarrow + \infty} \left( \int _{\Omega _{\alpha}}|\nabla _{h_{\alpha}}V_{\alpha}|^2 Ric_{g_\alpha}(x_{\alpha})_{ij}x^ix^j \; dv_{h_{\alpha}} \right) = \frac{Scal_g(x_0)}{n} \int_{\R^n}|\nabla V_0|^2 |x|^2\; dx
\end{equation}

\n By the Cartan expansion of $h_{\alpha}$ around $0$, we have

\[
dx = \left( 1 + \frac{\mu_{\alpha}^2}{6} Ric_{g_\alpha}(x_{\alpha})_{ij} x^i x^j + o(\mu_{\alpha}^2|x|_2^2) \right) dv_{h_{\alpha}}\, .
\]

\n Thus, by (\ref{E.20}), (\ref{E.14}) and (\ref{E.18}),

\begin{equation} \label{E.24}
\int_{\Omega_{\alpha}} F_\alpha (V_{\alpha})\; dx = 1 + \frac{Scal_g(x_0)}{6n} \mu_{\alpha}^2 \int_{\R^n} F(V_0) |x|^2\; dx + o \left( \mu_{\alpha}^2 \right)\, .
\end{equation}

\n By the Cartan expansion of $h_\alpha$ around $0$, since $r_\alpha \rightarrow 0$, we also have

\[
|\nabla V_{\alpha}|^2 = |\nabla_{h_{\alpha}} V_{\alpha}|^2 \left( 1 - \frac{\mu_{\alpha}^2}{6} |\nabla_{h_{\alpha}} V_{\alpha}|^{-2} \sum_i Rm_{g_\alpha}(x_{\alpha}) (\nabla_{h_{\alpha}} V^i_{\alpha},x,x,\nabla_{h_{\alpha}} V^i_{\alpha}) + o(\mu_{\alpha}^2 |x|^2) \right)\, .
\]

\n Since $v_0^i$ is a radial function and $\nabla_{h_{\alpha}}v_0^i$ e $x$ are pointwise colinear vector fields, by Lemma \ref{L.1}, it follows that (see \cite{D3})

\[
\int_{\Omega_{\alpha}} Rm_{g_\alpha}(x_{\alpha}) \left( \nabla_{h_{\alpha}} V^i_{\alpha},x,x, \nabla_{h_{\alpha}} V^i_{\alpha} \right) dv_{h_{\alpha}}\rightarrow 0\, .
\]

\n Consequently, by (\ref{E.21}) and (\ref{E.19}),

\begin{equation} \label{E.22}
\int_{\Omega_{\alpha}} |\nabla V_{\alpha}|^2\; dx = \int_{\Omega_{\alpha}} |\nabla_{h_{\alpha}} V_{\alpha}|^2\; dv_{h_{\alpha}} + \frac{Scal_g(x_0)}{6n} \mu_{\alpha}^2 \int_{\mathbb{R}^n} |\nabla V_0|^2 |x|^2\; dx + o(\mu_{\alpha}^2)\, .
\end{equation}

\n By (\ref{I.1}) and the system (\ref{NE1}), we have

\[
\mathcal{A}^\alpha \int_{\Omega_{\alpha}} |\nabla_{h_{\alpha}} V_{\alpha}|^2\; dv_{h_{\alpha}} < 1 - \mathcal{B}^\alpha \mu_{\alpha}^2 \int_{\Omega_{\alpha}} G_\alpha(V_{\alpha})\; dv_{h_{\alpha}}\, .
\]

\n By Lemma \ref{L.1}, we can write $V_0 = t_0 v_0$, where $t_0 \in X_F$ and $v_0 \in {\cal D}^{1,2}(\R^n)$ is an extremal function to $A_0(n)$ such that $v_0(0) = 1$ and $||v_0||_{2^*} = 1$. Therefore, again by (\ref{E.14}),

\begin{eqnarray*}
\int_{\Omega_{\alpha}} G_\alpha(V_{\alpha})\; dv_{h_{\alpha}} &=& \int_{\R^n} G(x_0,V_0)\; dx + o(1)\\
&=& \int_{\R^n} G(x_0,t_0)v_0^2\; dx + o(1)\,,
\end{eqnarray*}

\n so that

\begin{equation}\label{E.23}
\mathcal{A}^\alpha \int_{\Omega_{\alpha}} |\nabla_{h_{\alpha}} V_{\alpha}|^2\; dv_{h_{\alpha}} < 1 - m_{F,G}(x_0) \mathcal{B}^\alpha \mu_{\alpha}^2 \int_{\R^n} v_0^2\; dx + o(\mu_{\alpha}^2)\, .
\end{equation}

\n So, by (\ref{E.22}) and (\ref{E.23}),

\begin{equation} \label{E.25}
\mathcal{A}^\alpha \int_{\Omega_{\alpha}} |\nabla V_{\alpha}|^2\; dx \leq
\end{equation}

\[
1 - m_{F,G}(x_0) \mathcal{B}^\alpha \mu_{\alpha}^2 \int_{\R^n} v_0^2\; dx + \mathcal{A}_0(n,F) \frac{Scal_g(x_0)}{6n} \mu_{\alpha}^2 \int_{\R^n} |\nabla v_0|^2 |x|^2\; dx + o(\mu_{\alpha}^2)\, .
\]

\n Note also that (\ref{E.24}) can be rewritten as

\begin{equation} \label{E.26}
\int_{\Omega_{\alpha}} F_\alpha (V_{\alpha})\; dx = 1 + M_F \frac{Scal_g(x_0)}{6n} \mu_{\alpha}^2 \int_{\R^n} v_0^{2^*} |x|^2\; dx + o \left( \mu_{\alpha}^2 \right)\, .
\end{equation}

\n On the other hand, from the equation satisfied by $v_0$,

\[
- A_0(n) \Delta v_0 = v_0^{2^* - 1}\, ,
\]

\n we have

\[
\int_{\R^n} v_0^{2^*} |x|^2\; dx = A_0(n) \int_{\R^n} |\nabla v_0|^2 |x|^2\; dx + 2 A_0(n) \int_{\R^n} v_0 \nabla v_0 \cdot x\; dx\, .
\]

\n Thus, (\ref{E.26}) yields

\begin{equation} \label{E.27}
\int_{\Omega_{\alpha}} F_\alpha (V_{\alpha})\; dx =
\end{equation}

\[
1 + \mathcal{A}_0(n,F) \frac{Scal_g(x_0)}{6n} \mu_{\alpha}^2 \left( \int_{\R^n} |\nabla v_0|^2 |x|^2\; dx + 2 \int_{\R^n} v_0 \nabla v_0 \cdot x\; dx \right) + o \left( \mu_{\alpha}^2 \right)\, .
\]

\n Finally, replacing (\ref{E.25}) and (\ref{E.27}) in the potential type Euclidean Sobolev inequality

\[
\left( \int_{\Omega_{\alpha}} F_\alpha (V_{\alpha})\; dx \right)^{2/2^*} \leq \mathcal{A}^\alpha \int_{\Omega_{\alpha}} |\nabla V_{\alpha}|^2\; dx\, ,
\]

\n we get

\[
\left( \mathcal{B}^\alpha - \frac{n-2}{4(n-1)}\frac{\mathcal{A}_0(n,F)}{m_{F,G}(x_0)}Scal_g(x_0) \right) \mu_\alpha^2 + o(\mu_\alpha^2)
\leq 0\, .
\]

\n But, this contradicts

\[
\lim \limits_{\alpha \rightarrow + \infty} \mathcal{B}^\alpha = \frac{n-2}{4(n-1)}\frac{\mathcal{A}_0(n,F)}{m_{F,G}(x_0)}Scal_g(x_0) + \varepsilon_0
\]
\bl

\section{Continuous dependence of optimal constants}

In Subsection 2.1 we derived the continuity of the first optimal constant $\mathcal{A}_0(n,F,G,g)$ on the parameters $F$, $G$ and $g$, see Proposition \ref{P.2} or Remark 1 on the page 18. There remains to prove the continuity of the second optimal constant $\mathcal{B}_0(n,F,G,g)$ such as stated in Theorem \ref{Teo.1}. This continuity is a strategic one in the proof of Theorems \ref{Teo.2} and \ref{Teo.3} and its verification is partly inspired in some ideas used in the proof of Proposition \ref{P.7}. We should give more emphasis to the main difficult points in each specific context. Recall us also that Propositions \ref{P.4}, \ref{P.6} and \ref{P.7} are essential ingredients in the proof of Theorem \ref{Teo.1}.

We perform the proof by contradiction. Suppose that $\mathcal{B}_0(n,F,G,g)$ is discontinuous at some point $(F,G,g) \in \mathcal{F}_k  \times \mathcal{G}_k \times {\cal M}^n$. It means that there exist $\varepsilon_0 > 0$ and a sequence $(F_\alpha, G_\alpha, g_\alpha) \subset \mathcal{F}_k \times \mathcal{G}_k \times {\cal M}^n$ converging to $(F, G, g)$ such that, for any $\alpha > 0$,

\[
|\mathcal{B}_0(n,F_\alpha,G_\alpha,g_\alpha) - \mathcal{B}_0(n,F,G,g)| > \varepsilon_0\, .
\]

\n Then, at least, one of the following alternatives holds:

\begin{equation} \label{AL}
\mathcal{B}_0(n,F,G,g) - \mathcal{B}_0(n,F_\alpha,G_\alpha,g_\alpha) > \varepsilon_0\ \ {\rm or}\ \ \mathcal{B}_0(n,F_\alpha,G_\alpha,g_\alpha) - \mathcal{B}_0(n,F,G,g) > \varepsilon_0
\end{equation}

\n for infinitely many $\alpha$. If the first one happens, majoring $\mathcal{B}_0(n,F_\alpha,G_\alpha,g_\alpha)$ by $\mathcal{B}_0(n,F,G,g) - \varepsilon_0$ in the second sharp potential type Riemannian $L^2$-Sobolev inequality

\[
\left( \int_{M} F_\alpha(U)\; dv_{g_\alpha} \right)^{2/2^*} \leq \mathcal{A}_0(n,F_\alpha) \int_{M} |\nabla_{g_\alpha} U|^2 \; dv_{g_\alpha} +
\mathcal{B}_0(n,F_\alpha,G_\alpha,g_\alpha) \int_{M} G_\alpha(x,U)\; dv_{g_\alpha}\, ,
\]

\n one gets

\[
\left( \int_{M} F_\alpha(U)\; dv_{g_\alpha} \right)^{2/2^*} \leq \mathcal{A}_0(n,F_\alpha) \int_{M} |\nabla_{g_\alpha} U|^2 \; dv_{g_\alpha} +
\left( \mathcal{B}_0(n,F,G,g) - \varepsilon_0 \right) \int_{M} G_\alpha(x,U)\; dv_{g_\alpha}\, .
\]

\n Letting now $\alpha \rightarrow +\infty$ in the inequality above, we clearly derive a contradiction. If the second alternative holds, from the definition of $\mathcal{B}_0(n,F_\alpha,G_\alpha,g_\alpha)$, we easily deduce that

\begin{equation} \label{C-min}
\lambda_{\alpha} := \inf_{U \in \Lambda_{\alpha}} J_{\alpha }(U) < 1
\end{equation}

\n for any $\alpha > 0$, where $J_\alpha$ denotes the functional

\[
J_{\alpha }(U)= \mathcal{A}_0(n,F_\alpha) \int_M |\nabla_{g_\alpha} U|^2 \; dv_{g_\alpha} +(\mathcal{B}_0(n,F,G,g) + \varepsilon _0)
\int_{M} G_{\alpha}(x,U) \; dv_{g_\alpha}
\]

\n defined on $\Lambda_{\alpha} =\{ U \in H^{1,2}_k(M):\;  \int_{M} F_{\alpha}(U) \; dv_{g_\alpha} = 1\}$. By Remark \ref{R.6} of Proposition \ref{P.6}, we can assume that $F_\alpha$ and $G_\alpha$ are of class $C^1$ and converge, respectively, to $F$ and $G$ and that (\ref{C-min}) remains valid. In particular, $F_\alpha$, $G_\alpha$ and $g_\alpha$ satisfy (\ref{H.1}) and (\ref{H.2}). By Proposition \ref{exist}, (\ref{C-min}) leads to the existence of a minimizer $U_{\alpha} = (u_{\alpha}^1, \ldots, u_{\alpha}^k) \in \Lambda_{\alpha}$ related to $\lambda_{\alpha}$. In particular, each map $U_\alpha$ satisfies

\begin{equation} \label{S.3}
- \mathcal{A}_0(n,F_{\alpha}) \Delta_{g_\alpha} u_{\alpha}^i + \frac{1}{2}\left( \mathcal{B}_0(p,F,G,g) + \varepsilon_0 \right) \frac{\partial G_{\alpha}}{\partial t_i} (x,U_{\alpha}) = \frac{\lambda_{\alpha}}{2^*} \frac{\partial F_{\alpha}}{\partial t_i} (U_{\alpha}) \ \ {\rm on} \ M\
\end{equation}

\n with $i= 1,\ldots, k$. Moreover, by Proposition \ref{Reg}, $U_\alpha$ is of class $C^1$.

Our goal now is to study the behavior of the sequence $(U_{\alpha})$ as $\alpha$ tends to $+ \infty$. Of course, the sequence $(U_{\alpha})$ is bounded in $H^{1,2}_k(M)$, and so we can assume that it converges weakly in $H^{1,2}_k(M)$ and strongly in $L^q_k(M)$ to the map $U_0$ for any $1 \leq q < 2^*$. Since the sequence $(\lambda_\alpha)$ is bounded too, we can assume that it converges to $\lambda_0 \in [ 0, 1]$. We first show that $\lambda_0 = 1$. Otherwise,

\[
\left( \int_{M} F(U_\alpha)\; dv_{g} \right)^{2/2^*}
\]

\[
< \mathcal{A}_0(n,F) \int_{M} |\nabla_{g} U_\alpha|^2 \; dv_{g} +
\left( \mathcal{B}_0(n,F,G,g) + \varepsilon_0 \right) \int_{M} G(x,U_\alpha)\; dv_{g}
\]

\[
\leq (1 + \varepsilon_1) \left( \mathcal{A}_0(n,F_\alpha) \int_{M} |\nabla_{g_\alpha} U_\alpha|^2 \; dv_{g_\alpha} +
\left( \mathcal{B}_0(n,F,G,g) + \varepsilon_0 \right) \int_{M} G_\alpha(x,U_\alpha)\; dv_{g_\alpha} \right)
\]

\[
= (1 + \varepsilon_1) \lambda_\alpha \int_{M} F_\alpha(U_\alpha)\; dv_{g_\alpha} \leq (1 - \varepsilon_0) \int_{M} F(U_\alpha)\; dv_{g}\, ,
\]

\n which implies

\[
\left( \int_{M} F(U_\alpha)\; dv_{g} \right)^{2/n} > \frac{1}{1 - \varepsilon_0}
\]

\n for suitable $\varepsilon_0, \varepsilon_1 > 0$ small and $\alpha > 0$ large. But this last inequality contradicts the fact that

\[
\lim \limits_{\alpha \rightarrow +\infty} \int_{M} F(U_\alpha)\; dv_{g} = \lim \limits_{\alpha \rightarrow +\infty} \int_{M} F_\alpha(U_\alpha)\; dv_{g_\alpha} = 1\, .
\]

\n The next step is to eliminate the case when $U_0 \neq 0$. In fact, if this occurs we prove that the sequence $(U_\alpha)$ is compact in $L^{2^*}_k(M)$ and this leads us easily to a contradiction. Indeed, assuming that $(U_\alpha)$ converges strongly to $U_0$ in $L^{2^*}_k(M)$, we get

\[
\int_{M} F(U_0)\; dv_{g} = 1\, ,
\]

\n and so, letting $\alpha \rightarrow +\infty$ in the equality $J_\alpha(U_\alpha) = \lambda_\alpha$, we arrive at the following contradiction:

\[
\mathcal{A}_0(n,F) \int_{M} |\nabla_{g} U_0|^2 \; dv_{g} +
\left( \mathcal{B}_0(n,F,G,g) + \varepsilon_0 \right) \int_{M} G(x,U_0)\; dv_{g} \leq 1
\]

\[
= \left( \int_{M} F(U_0)\; dv_{g} \right)^{2/2^*}\, .
\]

\n In order to establish the compactness of $(U_\alpha)$, we evoke Propositions \ref{B-L0} and \ref{C-P1}. Let $(z_j)_{j \in {\cal T}} \subset M$ be the concentration points furnished in Proposition \ref{B-L0}. Fix $k \in {\cal T}$ and choose a cutoff function $\varphi_{\varepsilon} \in C_0^{\infty} \left(B_{g}(z_k, 2\varepsilon) \right)$ such that $0 \leq \varphi_{\varepsilon} \leq 1$, $\varphi_{\varepsilon} = 1$ in $B_{g}(z_k, \varepsilon)$ and
$|\nabla_{g} \varphi_{\varepsilon}| \leq \frac{C_1}{\varepsilon}$ for some constant $C_1 > 0$ independent of $\varepsilon$. Taking $\varphi_{\varepsilon} u^i_\alpha$ as a test function in the $i$th equation of the system (\ref{S.3}), integrating by parts and adding all equations, on the one hand,

\begin{equation} \label{E.30}
\lim_{\alpha \rightarrow \infty} \left[ \sum_{i=1}^k \left( \mathcal{A}_0(n,F_{\alpha}) \int_M \nabla_{g_\alpha} u^i_\alpha \cdot \nabla_{g_\alpha}
(\varphi_{\varepsilon} u^i_\alpha)\; dv_{g_\alpha} + \frac{1}{2}\left( \mathcal{B}_0(n,F,G,g) + \varepsilon_0 \right) \right. \right.
\end{equation}

\[
\left. \left. \times \int_M \frac{\partial G_{\alpha}}{\partial t_i} (x,U_{\alpha}) (\varphi_{\varepsilon} u^i_\alpha)\; dv_{g_\alpha} \right) \right]
\]

\[
= \lim_{\alpha \rightarrow \infty} \left( \sum_{i=1}^k \frac{\lambda_{\alpha}}{2^*} \int_M  \frac{\partial F_{\alpha}}{\partial t_i} (U_{\alpha}) (\varphi_{\varepsilon} u^i_\alpha)\; dv_{g_\alpha} \right)
\]

\[
= \lim_{\alpha \rightarrow \infty} \left( \lambda_\alpha \int_M F_\alpha (U_\alpha) \varphi_{\varepsilon}\; dv_{g_\alpha} \right) = \int_M \varphi_{\varepsilon}\; d\nu\, ,
\]

\n and, on the other hand,

\begin{equation} \label{E.31}
\lim_{\alpha \rightarrow \infty} \left[ \sum_{i=1}^k \left( \mathcal{A}_0(n,F_{\alpha}) \int_M \nabla_{g_\alpha} u^i_\alpha \cdot \nabla_{g_\alpha}
(\varphi_{\varepsilon} u^i_\alpha)\; dv_{g_\alpha} + \frac{1}{2}\left( \mathcal{B}_0(n,F,G,g) + \varepsilon_0 \right) \right. \right.
\end{equation}

\[
\left. \left. \times \int_M \frac{\partial G_{\alpha}}{\partial t_i} (x,U_{\alpha}) (\varphi_{\varepsilon} u^i_\alpha)\; dv_{g_\alpha} \right) \right]
\]

\[
= \lim_{\alpha \rightarrow \infty} \left( \frac{1}{2} \mathcal{A}_0(n,F_{\alpha}) \int_M \nabla_{g_\alpha} |U_\alpha|^2 \cdot \nabla_{g_\alpha}
\varphi_{\varepsilon}\; dv_{g_\alpha} + \mathcal{A}_0(n,F_{\alpha}) \int_M |\nabla_{g_\alpha} U_\alpha|^2 \varphi_{\varepsilon}\; dv_{g_\alpha} \right.
\]

\[
\left. + \left( \mathcal{B}_0(n,F,G,g) + \varepsilon_0 \right) \int_M G_{\alpha}(x,U_{\alpha}) \varphi_{\varepsilon}\; dv_{g_\alpha} \right)
\]

\[
= \frac{1}{2} \mathcal{A}_0(n,F) \lim_{\alpha \rightarrow \infty} \left( \int_M \nabla_{g_\alpha} |U_\alpha|^2 \cdot \nabla_{g_\alpha}
\varphi_{\varepsilon}\; dv_{g_\alpha} \right) +  \mathcal{A}_0(n,F) \int_M \varphi_{\varepsilon}\; d\mu
\]

\[
+  \left( \mathcal{B}_0(n,F,G,g) + \varepsilon_0 \right) \int_M G(x,U_0) \varphi_{\varepsilon}\; dv_{g}\, .
\]

\n We now let $\varepsilon \rightarrow 0$ on the right-hand side of (\ref{E.30}) and (\ref{E.31}). By H\"{o}lder's inequality, we find as $\varepsilon \rightarrow 0$,

\[
\lim_{\alpha \rightarrow \infty} \left| \int_M \nabla_{g_\alpha} |U_\alpha|^2 \cdot \nabla_{g_\alpha}
\varphi_{\varepsilon}\; dv_{g_\alpha} \right|
\]

\[
\leq C_0 \limsup_{\alpha \rightarrow \infty} \left[ \left( \int_M |\nabla_{g_\alpha} U_\alpha|^2\; dv_{g_\alpha} \right)^{1/2} \left( \int_{B_{g}(z_k, 2 \varepsilon) \setminus B_{g}(z_k, \varepsilon)} |\nabla_{g} \varphi_{\varepsilon}|^n\; dv_{g} \right)^{1/n} \right.
\]

\[
\left. \times \left( \int_{B_{g}(z_k, 2 \varepsilon) \setminus
B_{g}(z_k, \varepsilon)} F_\alpha(U_\alpha)\; dv_{g_\alpha} \right)^{1/2^*} \right]
\]

\[
\leq C_1 \left[ \frac{1}{\varepsilon^n} vol_{g} \left( B_{g}(z_k, 2 \varepsilon) \setminus B_{g}(z_k, \varepsilon) \right) \right]^{1/n}
\lim_{\alpha \rightarrow \infty} \left( \int_{B_{g}(z_k, 2 \varepsilon) \setminus B_{g}(z_k, \varepsilon)} F_\alpha(U_\alpha)\;
dv_{g_\alpha} \right)^{1/2^*}
\]

\[
\leq C_2 \left( \int_{B(z_k, 2 \varepsilon) \setminus B(z_k, \varepsilon)}  F(U_0)\;dv_{g} + \sum_{j \in {\cal T}} \nu_j \delta_{z_j} ( B_{g}(z_k, 2 \varepsilon) \setminus B_{g}(z_k, \varepsilon)) \right)^{1/2^*} \rightarrow 0 \ .
\]

\n Therefore, from (\ref{E.30}) and (\ref{E.31}), we derive

\[
\nu_k = \mathcal{A}_0(n,F) \mu_k\ .
\]

\n If $\mu_k > 0$, by Proposition \ref{C-P1}, one has

\[
\mu_k \geq \frac{1}{\mathcal{A}_0(n,F)}
\]

\n and so, we arrive at the following contradiction:

\[
1 = \lim_{\alpha \rightarrow \infty} \left( \mathcal{A}_0(n,F_{\alpha}) \int_M |\nabla_{g_\alpha} U_\alpha|^2\; dv_{g_\alpha} + \left( \mathcal{B}_0(n,F,G,g) + \varepsilon_0 \right) \int_M G_{\alpha}(x,U_{\alpha})\; dv_{g_\alpha} \right)
\]

\[
\geq \mathcal{A}_0(n,F) \int_M |\nabla_{g} U_0|^2\; dv_{g} + \left( \mathcal{B}_0(n,F,G,g) + \varepsilon_0 \right) \int_M G(x,U_0)\; dv_{g}
\]

\[
+ \mathcal{A}_0(n,F) \sum_{j \in {\cal T}} \mu _j > \mathcal{A}_0(n,F) \mu_k \geq 1\, .
\]

\n Here, it was used that $U_0 \neq 0$. Therefore, $\mu_j = 0$ and so, $\nu_j = 0$ for all $j \in {\cal T}$. Consequently, $\int_M F_\alpha (U_\alpha) dv_{g_\alpha}$ converges to $\int_M F (U_0) dv_{g}$, so that the strong convergence of $(U_\alpha)$ to $U_0$ in $L^{2^*}_k(M)$ follows from Proposition \ref{B-L0}. Hence, the remaining of proof consists in showing that $U_0 = 0$ can not hold. Clearly, if this happens, then the sequence $(|U_\alpha|)$ blows up in $L^\infty(M)$, since $(U_\alpha)$ converges to $0$ in $L^2_k(M)$ and

\[
1 = \int_M F_\alpha (U_\alpha)\; dv_{g_\alpha} \leq C_0 |U_\alpha(x_\alpha)|^{2^* - 2} \int_M |U_\alpha|^2\; dv_{g_\alpha},
\]

\n where $x_\alpha \in M$ is a maximum point of $|U_\alpha|$. We now perform a concentration study for the sequence $(U_\alpha)$ in order to establish the final contradiction. For this, we will base on part of the strategy used in the proof of Proposition \ref{P.7}. For a concentration point of $(U_\alpha)$ we mean a point $x_0 \in M$ such that

\[
\limsup_{\alpha \rightarrow + \infty} \int_{B_{g_\alpha}(x_0,\delta)} |U_{\alpha}|^{2^*}\; dv_{g_\alpha} > 0
\]

\n for any $\delta > 0$. Note that this notion is a natural extension of the corresponding scalar one.

We now split the main steps of proof into three lemmas.

\begin{lema} \label{L.3}
The concentration points set of the sequence $(U_\alpha)$ is single.
\end{lema}

\n {\bf Proof of Lemma \ref{L.3}.} First it is clear that

\[
\limsup_{\alpha \rightarrow + \infty} \int_{B_g(x_0,\delta)} F_\alpha(U_{\alpha})\; dv_{g_\alpha} > 0
\]

\n for some point $x_0 \in M$ and any $\delta > 0$, since $U_\alpha \in \Lambda_{\alpha}$. This is equivalent to say that $(U_\alpha)$ concentrates at $x_0$, by (\ref{H.1}) and (\ref{H.2}). The interesting part is indeed to prove the uniqueness of concentration points. For this, we fix an arbitrary number $\delta > 0$ and set

\[
\limsup_{\alpha \rightarrow + \infty} \int_{B_g(x_0,\delta)} F_\alpha(U_{\alpha})\; dv_{g_\alpha} = a \in (0,1]\, .
\]

\n Mimicking the proof of (\ref{R}), one easily checks that

\[
- \Delta_{g_\alpha} |U_\alpha| \leq \mathcal{A}_0(n,F_\alpha)^{-1} |U_\alpha|^{-1} F_\alpha(U_\alpha) \ \ {\rm on}\ M\, .
\]

\n We now take a cutoff function $\varphi \in C^{\infty}_0(B_{g}(x_0,\delta))$ such that $0 \leq \varphi \leq 1$, $\varphi = 1$ in $B_{g}(x_0,\delta /2)$ and $|\nabla_{g_\alpha} \varphi| \leq C_0$ with $C_0 > 0$ independent of $\alpha$. Let $m > 0$ be a small number to be chosen later. Regarding $\varphi^2 |U_\alpha|^{m+1}$ as a test function in the preceding inequality, we can write

\begin{equation}\label{E.33}
\int_M \nabla_{g_\alpha} |U_\alpha| \cdot \nabla_{g_\alpha} (\varphi^2 |U_\alpha|^{m + 1}) \; dv_{g_\alpha} \leq \mathcal{A}_0(n,F_\alpha)^{-1} \int_M \varphi^2 |U_\alpha|^{m} F_\alpha(U_\alpha) \; dv_{g_\alpha}\, .
\end{equation}

\n Developing the left-hand side of (\ref{E.33}), we get

\begin{equation}\label{E.34}
(m + 1) \int_M \varphi^2 |U_\alpha|^m | \nabla_{g_\alpha} |U_\alpha| |^2 \; dv_{g_\alpha} \leq \mathcal{A}_0(n,F_\alpha)^{-1} \int_M \varphi^2 |U_\alpha|^{m} F_\alpha(U_\alpha) \; dv_{g_\alpha}
\end{equation}

\[
- \int_M |U_\alpha|^{m + 1} \nabla_{g_\alpha} (\varphi^2) \cdot \nabla_{g_\alpha} |U_\alpha|\; dv_{g_\alpha}\, .
\]

\n On the other hand, given $\varepsilon >0$ one finds a constant $C_{\varepsilon}>0$, independent of $\alpha$, such that

\[
\int_M |\nabla_{g_\alpha} (\varphi |U_\alpha|^{\frac{m+2}{2}})|^2 \; dv_{g_\alpha} \leq (1 + \varepsilon) \frac{(m+2)^2}{4} \int_M \varphi^2 |U_\alpha|^m | \nabla_{g_\alpha} |U_\alpha| |^2 \; dv_{g_\alpha}
\]

\[
+ C_{\varepsilon} ||\nabla_{g_\alpha} \varphi||_{\infty}^2 \int_M |U_\alpha|^{m+2} \; dv_{g_\alpha}\, .
\]

\n Plugging (\ref{E.34}) into this, we get

\begin{equation} \label{E.35}
\int_M |\nabla_{g_\alpha} (\varphi |U_\alpha|^{\frac{m+2}{2}})|^2 \; dv_{g_\alpha} \leq (1 + \varepsilon) \frac{(m+2)^2}{4(m+1)} M_{F_\alpha}^{-2/2^*} A_0(n)^{-1} \int_M \varphi^2 |U_\alpha|^{m} F_\alpha(U_\alpha) \; dv_{g_\alpha}
\end{equation}

\[
- (1 + \varepsilon) \frac{(m+2)^2}{4(m+1)} \int_M |U_\alpha|^{m+1} \nabla_{g_\alpha} (\varphi^2) \cdot \nabla_{g_\alpha} |U_\alpha| \; dv_{g_\alpha}+ C_{\varepsilon} ||\nabla_{g_\alpha} \varphi||_{\infty}^2 \int_M |U_\alpha|^{m+2} \; dv_{g_\alpha}\, .
\]

\n Using now H\"{o}lder's inequality, we derive

\begin{equation} \label{E.36}
\int_M \varphi^2 |U_\alpha|^{m} F_\alpha(U_\alpha) \; dv_{g_\alpha} = \int_M \varphi^2 |U_\alpha|^{m} F_\alpha(U_\alpha)^{2/2^*} F_\alpha(U_\alpha)^{(2^*-2)/2^*} \; dv_{g_\alpha}
\end{equation}

\[
\leq M_{F_\alpha}^{2/2^*} \int_M \varphi^2 |U_\alpha|^{m+2} F_\alpha(U_\alpha)^{(2^*-2)/2^*}\; dv_{g_\alpha}
\]

\[
\leq M_{F_\alpha}^{2/2^*} \left( \int_M (\varphi |U_\alpha|^{\frac{m+2}{2}})^{2^*}\; dv_{g_\alpha} \right)^{2/2^*}  \left( \int_{B_{g_\alpha}(x_0,\delta)} F_\alpha(U_\alpha)\; dv_{g_\alpha} \right)^{(2^* - 2)/2^*}
\]

\n and

\begin{equation} \label{E.37}
\left| \int_M |U_\alpha|^{m+1} \nabla_{g_\alpha} (\varphi^2) \cdot \nabla_{g_\alpha} |U_\alpha|\; dv_{g_\alpha} \right| \leq 2||\nabla_{g_\alpha} \varphi||_{\infty} \int_M |U_\alpha|^{m+1} |\nabla_{g_\alpha} U_\alpha|\; dv_{g_\alpha}
\end{equation}

\[
\leq 2||\nabla_{g_\alpha} \varphi||_{\infty} \left( \int_M |\nabla_{g_\alpha} U_\alpha|^2\; dv_{g_\alpha} \right)^{1/2} \left( \int_M |U_\alpha|^{2m+2}\; dv_{g_\alpha} \right)^{1/2}
\]

\[
\leq 2 M_{F_\alpha}^{-1/2^*} A_0(n)^{-1/2} ||\nabla_{g_\alpha} \varphi||_{\infty} \left( \int_M |U_\alpha|^{2m+2}\; dv_{g_\alpha} \right)^{1/2}\, .
\]

\n Thanks to the asymptotically sharp $L^2$-Sobolev inequality, for any $\varepsilon > 0$ there exists a positive constant $B_\varepsilon$, independent of $\alpha$, such that

\begin{equation} \label{E.38}
\left( \int_M (\varphi |U_\alpha|^{\frac{m+2}{2}})^{2^*}\; dv_{g_\alpha} \right)^{2/2^*} \leq (A_0(n) + \varepsilon) \int_M |\nabla_{g_\alpha} (\varphi |U_\alpha|^{\frac{m+2}{2}})|^2 \; dv_{g_\alpha}
\end{equation}

\[
 + B_\varepsilon \int_{M} |U_\alpha|^{m+2}\; dv_{g_\alpha}\, .
\]

\n Inserting (\ref{E.36}), (\ref{E.37}) and (\ref{E.38}) into (\ref{E.35}), we arrive at

\begin{equation} \label{E.39}
A_{\alpha} \left( \int_M (\varphi |U_\alpha|^{\frac{m+2}{2}})^{2^*}\; dv_{g_\alpha} \right)^{2/2^*} \leq B_{\alpha} \int_{M} |U_\alpha|^{m+2}\; dv_{g_\alpha} + C_{\alpha} \left( \int_M |U_\alpha|^{2m+2}\; dv_{g_\alpha} \right)^{1/2}\, .
\end{equation}

\n where

\[
A_{\alpha } = 1 - (1 + \varepsilon) \frac{(m+2)^2}{4(m+1)} A_0(n)^{-1} (A_0(n) + \varepsilon) \left( \int_{B_g(x_0,\delta)} F_\alpha(U_\alpha)\; dv_{g_\alpha} \right)^{(2^* - 2)/2^*}\ ,
\]

\[
B_{\alpha} = C_{\varepsilon} (A_0(n) + \varepsilon) ||\nabla_{g_\alpha} \varphi||_{\infty}^2  + B_{\varepsilon}
\]

\n and

\[
C_{\alpha} = 2 (1 + \varepsilon) \frac{(m+2)^2}{4(m+1)} M_{F_\alpha}^{-1/2^*} A_0(n)^{-1} (A_0(n) + \varepsilon) ||\nabla_{g_\alpha} \varphi||_{\infty}\ .
\]

\n What we wish prove now is that $a=1$. Suppose by contradiction that $a < 1$. In this case, we choose $\varepsilon > 0$ and $m > 0$ small enough such that

\[
A_\alpha \geq C > 0,\ \ m + 2 \leq 2^*,\ \ 2m + 2 \leq 2^*\ \ {\rm and}\ \ 2 \leq 2^* - \frac{m}{2^* - 2} < 2^*\, .
\]

\n With that choice, (\ref{E.39}) produces

\[
\left( \int_M (\varphi |U_\alpha|^{\frac{m+2}{2}})^{2^*}\; dv_{g_\alpha} \right)^{2/2^*} \leq C_1
\]

\n for $\alpha > 0$ large, where $C_1 > 0$ does not depend on $\alpha$. Using again H\"{o}lder's inequality, we get

\begin{eqnarray*}
\int_{B_{g_\alpha}(x_0,\frac{\delta}{4})} |U_{\alpha}|^{2^*}\; dv_{g_\alpha} &=& \int_{B_g(x_0,\frac{\delta}{2})} |U_{\alpha}|^{m + 2} |U_{\alpha}|^{2^* - 2 - m}\; dv_{g_\alpha}\\
&\leq& \left( \int_M (\varphi |U_\alpha|^{\frac{m+2}{2}})^{2^*}\; dv_{g_\alpha} \right)^{2/2^*} \left( \int_M |U_\alpha|^{2^* - \frac{m}{2^* - 2}}\;
dv_{g_\alpha} \right)^{(2^* - 2)/2^*}\\
&\leq& C_1 \left( \int_M |U_\alpha|^{2^* - \frac{m}{2^* - 2}}\;
dv_{g_\alpha} \right)^{(2^* - 2)/2^*}\, .
\end{eqnarray*}

\n Since $2 \leq 2^* - m/(2^* - 2) < 2^*$ and $(U_\alpha)$ converges to $0$ in $L^2_k(M)$, a simple interpolation scheme reveals that

\[
\limsup_{\alpha \rightarrow + \infty} \int_{B_{g_\alpha}(x_0,\frac{\delta}{4})} |U_{\alpha}|^{2^*}\; dv_{g_\alpha} = 0\, ,
\]

\n contradicting the fact that $x_0$ is a concentration point of $(U_\alpha)$. Therefore, $a = 1$ and this easily implies that $x_0$ is the unique concentration point of $(U_\alpha)$.\bl \\

\begin{lema} \label{L.4}
The sequence $(x_\alpha)$ converges to the unique concentration point $x_0$ of $(U_\alpha)$. Furthermore, for any $\delta > 0$, we have

\[
|U_\alpha| \rightarrow 0 \ \ {\rm in}\ \ C^0(M \setminus B_{g_\alpha}(x_0, \delta))\, .
\]

\end{lema}

\n {\bf Proof of Lemma \ref{L.4}.} Let $\tilde{x}_0 \in M$ be a limit point of $(x_\alpha)$. We first prove that $(U_\alpha)$ concentrates at $\tilde{x}_0$. Otherwise, by Lemma \ref{L.3},

\[
\limsup_{\alpha \rightarrow + \infty} \int_{B_{g_\alpha}(\tilde{x}_0, 2\delta)} |U_{\alpha}|^{2^*}\; dv_{g_\alpha} = 0
\]

\n for any $\delta > 0$ small enough. Thanks to the estimate (\ref{E.39}), there exist constants $m_1 > 0$ and $C_1 > 0$, independent of $\alpha$, such that

\[
\int_{B_{g_\alpha}(\tilde{x}_0, 2\delta)} |U_\alpha|^{\frac{(m_1+2)2^*}{2}} \; dv_{g_\alpha} \leq C_1
\]

\n for $\alpha > 0$ large. Applying Proposition \ref{P.4}, we deduce that

\[
\sup_{B_{g_\alpha}(\tilde{x}_0,\delta)} |U_\alpha| \leq C_0 \delta^{-n/p}  \left( \int_{B_{g_\alpha}(\tilde{x}_0, 2\delta)} |U_\alpha|^{2^*}\; dv_{g_\alpha} \right)^{1/2^*} \rightarrow 0\, .
\]

\n But this clearly contradicts the fact that $(|U_\alpha(x_\alpha)|)$ blows up as $\alpha$ tends to $\infty$, so that $\tilde{x}_0 = x_0$. For the remaining conclusion, since for any $\delta > 0$,

\[
\limsup_{\alpha \rightarrow + \infty} \int_{M \setminus B_{g_\alpha}(x_0, \frac{\delta}{2})} |U_{\alpha}|^{2^*}\; dv_{g_\alpha} = 0\, ,
\]

\n going back to (\ref{E.39}), we obtain other positive constants $m_2$ and $C_2$, independent of $\alpha$, such that

\[
\int_{M \setminus B_{g_\alpha}(x_0, \frac{\delta}{2})} |U_\alpha|^{\frac{(m_2+2)2^*}{2}} \; dv_{g_\alpha} \leq C_2\, .
\]

\n Evoking again Proposition \ref{P.4}, we end up with

\[
\sup_{M \setminus B_{g_\alpha}(x_0, \delta)} |U_\alpha| \leq C_3  \left( \int_{M \setminus B_{g_\alpha}(x_0, \frac{\delta}{2})} |U_\alpha|^{2^*}\; dv_{g_\alpha} \right)^{1/2^*} \rightarrow 0\, .
\]

\n \bl \\

With this lemma at hand, we establish the following concentration estimate of kind $L^2$:

\begin{lema} \label{L.5}
For any $\delta > 0$,

\[
\lim \limits_{\alpha \rightarrow +\infty} \frac{\int_{M \setminus B_{g_\alpha}(x_\alpha,\delta)} |U_{\alpha}|^2 \; dv_{g_\alpha}}{\int_M |U_{\alpha}|^2 \; dv_{g_\alpha}} = 0\, .
\]
\end{lema}

\n {\bf Proof of Lemma \ref{L.5}.} Thanks to Lemma \ref{L.4}, Proposition \ref{P.4} applied to the system (\ref{S.3}) produces

\begin{equation} \label{E.40}
\int_{M \setminus B_{g_\alpha}(x_\alpha,\delta)} |U_{\alpha}|^2 \; dv_{g_\alpha} \leq \left( \sup_{M \setminus B_{g_\alpha}(x_\alpha,\delta)} |U_{\alpha}| \right) \left( \int_{M \setminus B_{g_\alpha}(x_\alpha,\delta)} |U_{\alpha}|\; dv_{g_\alpha} \right)
\end{equation}

\[
\leq C_0 ||U_\alpha||_{L^2_k(M)} \int_{M} |U_{\alpha}|\; dv_{g_\alpha}\, .
\]

\n Mimicking the computation of (\ref{R}) to (\ref{S.3}), one has

\[
- \mathcal{A}_0(n,F_{\alpha}) \Delta_{g_\alpha} |U_\alpha| + \left( \mathcal{B}_0(n,F,G,g) + \varepsilon_0 \right) |U_\alpha|^{-1} G_{\alpha}(x,U_{\alpha}) \leq |U_\alpha|^{-1} F_{\alpha}(U_{\alpha})\, .
\]

\n Integrating now this inequality over $M$ and using the fact that $F_\alpha$ and $G_\alpha$ satisfy (\ref{H.2}), we get

\[
\int_M |U_{\alpha}|\; dv_{g_\alpha} \leq C_1 \int_M |U_{\alpha}|^{2^* - 1}\; dv_{g_\alpha}\, ,
\]

\n where $C_1 > 0$ does not depend on $\alpha$. Replacing this inequality in (\ref{E.40}), one arrives at

\[
\int_{M \setminus B_{g_\alpha}(x_\alpha,\delta)} |U_{\alpha}|^2 \; dv_{g_\alpha} \leq C_0 C_1 ||U_\alpha||_{L^2_k(M)} ||U_\alpha||_{L^{2^*-1}_k(M)}^{2^*-1}
\]

\n We now analyze two possibilities. If $2^* - 1 \leq 2$, H\"{o}lder's inequality provides

\[
\int_{M \setminus B_{g_\alpha}(x_\alpha,\delta)} |U_{\alpha}|^2 \; dv_{g_\alpha} \leq C_2 ||U_\alpha||_{L^2_k(M)} ||U_\alpha||_{L^2_k(M)}^{2^*-1}\, ,
\]

\n so that

\[
\frac{\int_{M \setminus B_{g_\alpha}(x_\alpha,\delta)} |U_{\alpha}|^2 \; dv_{g_\alpha}}{\int_M |U_{\alpha}|^2 \; dv_{g_\alpha}} \leq C_2 ||U_\alpha||_{L^2_k(M)}^{2^*-2} \rightarrow 0\, .
\]

\n Otherwise, if $2^* - 1 > 2$, an interpolation inequality combined with the fact that $U_\alpha \in \Lambda_{\alpha}$ readily yield

\[
\frac{\int_{M \setminus B_{g_\alpha}(x_\alpha,\delta)} |U_{\alpha}|^2 \; dv_{g_\alpha}}{\int_M |U_{\alpha}|^2 \; dv_{g_\alpha}} \leq C_3 ||U_\alpha||_{L^2_k(M)}^{\frac{2^*-3}{2^*-2}} \rightarrow 0\, .
\]

\n This ends the proof. \bl\\

\n {\bf Conclusion of the proof of Theorem \ref{Teo.1}.} Our aim here is to complete the proof of Theorem \ref{Teo.1} by seeking a final contradiction. By Proposition \ref{P.7}, for each $\varepsilon > 0$ there exists a constant $r_0 > 0$ such that, for any map $U \in C^1_{0,k}(B_{g_\alpha}(x_0,r_0))$ and $\alpha > 0$ large,

\begin{equation} \label{E.47}
\left ( \int_M F_\alpha(U)\; dv_{g_\alpha}  \right )^{2/2^*} \leq  \mathcal{A}_0(n,F_\alpha) \int_M |\nabla_{g_\alpha} U|^2\; dv_{g_\alpha} + \mathcal{B}_\varepsilon(F_\alpha,G_\alpha,g_\alpha) \int_M G_\alpha(x,U)\; dv_{g_\alpha}\, ,
\end{equation}

\n with

\[
\mathcal{B}_\varepsilon(F_\alpha,G_\alpha,g_\alpha) := \frac{n-2}{4(n-1)} \frac{\mathcal{A}_0(n,F_\alpha)}{m_{F_\alpha,G_\alpha}(x_0)} Scal_{g_\alpha}(x_0) + \varepsilon\, ,
\]

\n where

\[
m_{F_\alpha,G_\alpha}(x) := \min \limits_{t \in X_{F_\alpha}} G_\alpha(x,t)\, ,
\]

\n and $X_{F_\alpha} = \{ t \in \mathbb{S}^{k-1}:\; F_\alpha(t) = M_{F_\alpha}\}$. Choose now $0 < \varepsilon < \varepsilon_0$ and consider a cutoff function $\varphi_\alpha \in C^\infty(B_{g_\alpha}(x_\alpha, r_0))$ with $0 \leq \varphi_\alpha \leq 1$, $\varphi_\alpha = 1$ in $B_{g_\alpha}(x_\alpha, r_0/2)$ and $|\nabla_{g_\alpha} \varphi_\alpha| \leq C_0$, where $C_0 > 0$ does not depend on $\alpha$. Taking $\varphi_\alpha^2 u_\alpha^i$ as a test function in the $i$th equation of (\ref{S.3}), integrating by parts and adding all equations, one easily checks that

\begin{equation} \label{E.48}
\mathcal{A}_0(n,F_{\alpha}) \int_M |\nabla_{g_\alpha} (\varphi_\alpha U_\alpha)|^2\; dv_{g_\alpha} + \left( \mathcal{B}_0(n,F,G,g) + \varepsilon_0 \right) \int_M G_{\alpha}(x,U_{\alpha}) \varphi_\alpha^2\; dv_{g_\alpha}
\end{equation}

\[
= \lambda_\alpha \int_M F_{\alpha}(U_{\alpha}) \varphi_\alpha^2\; dv_{g_\alpha} + \mathcal{A}_0(n,F_{\alpha}) \int_M |\nabla_{g_\alpha} \varphi_\alpha|^2 |U_\alpha|^2\; dv_{g_\alpha}\, .
\]

\n Plugging (\ref{E.48}) into (\ref{E.47}), we derive

\begin{equation} \label{E.49}
\left( \int_M F_\alpha(U_\alpha) \varphi_\alpha^{2^*}\; dv_{g_\alpha}  \right)^{2/2^*} + \left( \mathcal{B}_0(n,F,G,g) - \mathcal{B}_\varepsilon(F_\alpha,G_\alpha,g_\alpha) + \varepsilon_0 \right) \int_M G_{\alpha}(x,U_{\alpha}) \varphi_\alpha^2\; dv_{g_\alpha}
\end{equation}

\[
\leq \int_M F_{\alpha}(U_{\alpha}) \varphi_\alpha^2\; dv_{g_\alpha} + \mathcal{A}_0(n,F_{\alpha}) \int_M |\nabla_{g_\alpha} \varphi_\alpha|^2 |U_\alpha|^2\; dv_{g_\alpha}\, .
\]

\n On the other hand, by H\"{o}lder's inequality,

\[
\int_M F_\alpha(U_\alpha) \varphi_\alpha^2\; dv_{g_\alpha} = \int_M F_\alpha(U_\alpha)^{(2^*-2)/2^*} F_\alpha(U_\alpha)^{2/2^*} \varphi_\alpha^2\; dv_{g_\alpha}
\]

\[
\leq \left( \int_M F_\alpha(U_\alpha)\; dv_{g_\alpha}  \right)^{(2^*-2)/2^*} \left( \int_M F_\alpha(U_\alpha) \varphi_\alpha^{2^*}\; dv_{g_\alpha} \right)^{2/2^*}
\]

\[
= \left( \int_M F_\alpha(U_\alpha) \varphi_\alpha^{2^*}\; dv_{g_\alpha} \right)^{2/2^*}
\]

\n and, by Proposition \ref{P.3},

\[
\mathcal{B}_0(n,F,G,g) m_{F,G}(x_0) \geq \frac{n-2}{4(n-1)} \mathcal{A}_0(n,F) S_g(x_0)\, .
\]

\n So, (\ref{E.49}) yields

\[
\left[ \frac{n-2}{4(n-1)} \left( \frac{\mathcal{A}_0(n,F)}{m_{F,G}(x_0)} Scal_{g}(x_0) - \frac{\mathcal{A}_0(n,F_\alpha)}{m_{F_\alpha,G_\alpha}(x_0)} Scal_{g_\alpha}(x_0) \right) + \varepsilon_0 - \varepsilon \right] \int_M G_{\alpha}(x,U_{\alpha}) \varphi_\alpha^2\; dv_{g_\alpha}
\]

\[
\leq \mathcal{A}_0(n,F_{\alpha}) \int_M |\nabla_{g_\alpha} \varphi_\alpha|^2 |U_\alpha|^2\; dv_{g_\alpha}\, ,
\]

\n which implies that

\[
\frac{\int_{M \setminus B_{g_\alpha}(x_\alpha, r_0/2)} |U_\alpha|^2\; dv_{g_\alpha}}{\int_{M} |U_{\alpha}|^2\; dv_{g_\alpha}} \geq C_1
\]

\n for $\alpha > 0$ large, where $C_1 > 0$ does not depend on $\alpha$. Of course, this fact contradicts Lemma \ref{L.5} and so we finish the proof of Theorem \ref{Teo.1}. \bl \\

\section{The vector duality}

In this section we prove Theorem \ref{Teo.2} which concerns the duality between the existence of extremal maps to (\ref{B-opt-v}) and the exact value of $\mathcal{B}_0(n,F,G,g)$ for general potential functions. The proof of this result goes in a similar spirit of the preceding proof, except that Theorem \ref{Teo.1} plays a crucial role.

Let $(F,G,g) \in \mathcal{F}_k  \times \mathcal{G}_k \times {\cal M}^n$. The conclusion readily follows if

\[
\mathcal{B}_0(n,F,G,g) G(x_0,t_0) = \frac{n-2}{4(n-1)} \mathcal{A}_0(n,F) S_g(x_0)
\]

\n for some $x_0 \in M$ and a maximum point $t_0$ of $F$ on $\s^{k-1}$. Otherwise, we can assume that, for any $x \in M$,

\begin{equation} \label{E.50}
\mathcal{B}_0(n,F,G,g) \geq \frac{n-2}{4(n-1)} \frac{\mathcal{A}_0(n,F)}{m_{F,G}(x)} S_g(x) + \varepsilon_0
\end{equation}

\n for some $\varepsilon_0 > 0$ small. In this case, we choose a sequence $(\alpha)$ converging to $\mathcal{B}_0(n,F,G,g)$ with $\alpha < \mathcal{B}_0(n,F,G,g)$. By Proposition \ref{P.6}, there exist sequences $(F_\alpha) \subset \mathcal{F}_k$ and $(G_\alpha) \subset \mathcal{G}_k$ of potential functions of $C^1$ class converging, respectively, to $F$ and $G$ as $\alpha \rightarrow \mathcal{B}_0(n,F,G,g)$. In particular, $(F_\alpha)$ and $(G_\alpha)$ satisfy (\ref{H.2}) for $\alpha$ near enough $\mathcal{B}_0(n,F,G,g)$. By Theorem \ref{Teo.1}, we can also assume

\begin{equation} \label{E.51}
\alpha < \mathcal{B}_0(n,F_{\alpha},G_{\alpha},g)\, .
\end{equation}

\n Given $\alpha > 0$, we consider the functional

\[
J_{\alpha}(U) = \mathcal{A}_0(n,F_{\alpha}) \int_M |\nabla_g U|^2\; dv_g + \alpha \int_{M} G_{\alpha}(x,U)\; dv_g
\]

\n defined on $\Lambda_{\alpha} =\{ U \in H^{1,2}_k(M):\;  \int_M F_{\alpha}(U)\; dv_g = 1\}$, and its respective infimum

\begin{equation} \label{C-min1}
\lambda_{\alpha} := \inf_{U \in \Lambda_{\alpha}} J_{\alpha}(U)\, .
\end{equation}

\n Since (\ref{E.51}) leads readily to $\lambda_{\alpha} < 1$, it follows that $\lambda_{\alpha}$ is achieved by a map $U_{\alpha} = (u_{\alpha}^1, \ldots, u_{\alpha}^k) \in \Lambda_{\alpha}$, by Proposition \ref{exist}. In particular, the map $U_\alpha$ satisfies

\begin{equation} \label{S.4}
- \mathcal{A}_0(n,F_{\alpha}) \Delta_g u_{\alpha}^i + \frac{\alpha}{2} \frac{\partial G_{\alpha}}{\partial t_i} (x,U_{\alpha}) = \frac{\lambda_{\alpha}}{2^*} \frac{\partial F_{\alpha}}{\partial t_i} (U_{\alpha}) \ \ {\rm on} \ M\, ,
\end{equation}

\n which in turn implies $C^1$ regularity of $U_\alpha$, by Proposition \ref{Reg}. Proceeding exactly as in the proof of Theorem \ref{Teo.1}, up to a subsequence, $(\lambda_\alpha)$ converges to $1$ and $(U_\alpha)$ converges weakly in $H^{1,2}_k(M)$ and strongly in $L^q_k(M)$ to a map $U_0$ for any $1 \leq q < 2^*$. Moreover, we claim that $U_0 \neq 0$. In fact, as in the previous proof, the nullity of $U_0$ lead to conclusions of Lemmas \ref{L.3}, \ref{L.4} and \ref{L.5} for $(U_\alpha)$. In particular, if $x_\alpha \in M$ is a maximum point of $|U_\alpha|$, we then can assume that $(x_\alpha)$ converges to a point $x_0$ and, for any $\delta > 0$,

\begin{equation} \label{E.52}
\lim \limits_{\alpha \rightarrow \mathcal{B}_0(n,F,G,g)} \frac{\int_{M \setminus B_g(x_\alpha,\delta)} |U_{\alpha}|^2 \; dv_g}{\int_M |U_{\alpha}|^2 \; dv_g} = 0\, .
\end{equation}

\n We now fix $0 < \varepsilon < \varepsilon_0$ and apply Proposition \ref{P.7} to $F_\alpha$, $G_\alpha$ and $g$ in order to obtain the local sharp inequality provided there for some radius $r_0 > 0$. We then choose a cutoff function $\varphi_\alpha \in C^\infty(B_g(x_\alpha, r_0))$ with $0 \leq \varphi_\alpha \leq 1$, $\varphi_\alpha = 1$ in $B_g(x_\alpha, r_0/2)$ and $|\nabla_g \varphi_\alpha| \leq C_0$, where $C_0 > 0$ does not depend on $\alpha$, and work with the test function $\varphi_\alpha^2 u_\alpha^i$ in (\ref{S.4}). After some manipulations as those ones done in the final part of the proof of Theorem \ref{Teo.1}, we deduce that

\[
\left( \alpha - \mathcal{B}_\varepsilon(F_\alpha,G_\alpha,g)\right) \int_M G_{\alpha}(x,U_{\alpha}) \varphi_\alpha^2\; dv_g
\leq \mathcal{A}_0(n,F_{\alpha}) \int_M |\nabla_g \varphi_\alpha|^2 |U_\alpha|^2\; dv_g\, .
\]

\n From (\ref{E.50}) and the definition of $\mathcal{B}_\varepsilon(F_\alpha,G_\alpha,g)$ in Proposition \ref{P.7}, we find a constant $C_1 > 0$ such that

\[
\frac{\int_{M \setminus B_g(x_\alpha, r_0/2)} |U_\alpha|^2\; dv_g}{\int_{M} |U_{\alpha}|^2\; dv_g} \geq C_1
\]

\n for $\alpha$ near enough $\mathcal{B}_0(n,F,G,g)$. But this last inequality obviously contradicts (\ref{E.52}) and, consequently, one has $U_0 \neq 0$. We wish to prove that $U_0$ is an extremal map to (\ref{B-opt-v}). Evoking Propositions \ref{B-L0} and \ref{C-P1} in exactly as done in the proof of Theorem \ref{Teo.1}, we derive the compactness of $(U_\alpha)$ in $L^{2^*}_k(M)$. In particular,

\[
\int_M F(U_0)\; dv_g = 1\, ,
\]

\n since $U_\alpha \in \Lambda_\alpha$. Letting now $\alpha \rightarrow \mathcal{B}_0(n,F,G,g)$ in the equality

\[
\mathcal{A}_0(n,F_{\alpha}) \int_M |\nabla_g U_\alpha|^2\; dv_g + \alpha \int_{M} G_{\alpha}(x,U_\alpha)\; dv_g = \lambda_\alpha\, ,
\]

\n we discover that

\[
\mathcal{A}_0(n,F) \int_M |\nabla_g U_0|^2\; dv_g + \mathcal{B}_0(n,F,G,g) \int_{M} G_{\alpha}(x,U_0)\; dv_g \leq 1\, .
\]

\n Therefore, $U_0$ is an extremal map to (\ref{B-opt-v}) and the proof is concluded. \bl \\

\section{Compactness of extremal maps}

This section is devoted to the compactness problem of extremal maps to (\ref{B-opt-v}). We shall establish a general result of $L^{2^*}_k$ and $C^0_k$ compactness such as stated in Theorem \ref{Teo.3}. The proof of the first part of this result is quite delicate due to the complete absence of smoothness assumptions. In particular, our extremal maps do not solve systems of kind (\ref{S.1}). Our great challenging here is to prove that versions of Lemmas \ref{L.3}, \ref{L.4} and \ref{L.5} hold in the present context.

We recall that, by assumption,  $((F_\alpha,G_\alpha,g_\alpha))$ converges to $(F,G,g)$ in $\mathcal{F}_k  \times \mathcal{G}_k \times {\cal M}^n$ as $\alpha \rightarrow + \infty$ and $(F,G,g)$ satisfies, for any $x \in M$,

\begin{equation} \label{E.50}
\mathcal{B}_0(n,F,G,g) > \frac{n-2}{4(n-1)} \frac{\mathcal{A}_0(n,F)}{m_{F,G}(x)} S_g(x)\, .
\end{equation}

\n So, by Theorem \ref{Teo.1}, we can take a constant $\varepsilon_0 > 0$ such that

\begin{equation} \label{E.50}
\mathcal{B}_0(n,F_\alpha,G_\alpha,g_\alpha) \geq \frac{n-2}{4(n-1)} \frac{\mathcal{A}_0(n,F_\alpha)}{m_{F_\alpha,G_\alpha}(x)} S_{g_\alpha}(x) + \varepsilon_0
\end{equation}

\n for any $x \in M$ and $\alpha > 0$ large. Consider a sequence $(U_\alpha)$ of extremal maps with $U_\alpha \in {\cal E}_k(F_\alpha,G_\alpha,g_\alpha)$, namely $U_\alpha$ satisfies $\int_M F_\alpha (U_\alpha) dv_{g_\alpha} = 1$ and

\[
\mathcal{A}_0(n,F_{\alpha}) \int_M |\nabla_{g_\alpha} U_\alpha|^2\; dv_{g_\alpha} + \mathcal{B}_0(n,F_\alpha,G_\alpha,g_\alpha) \int_{M} G_{\alpha}(x,U_\alpha)\; dv_{g_\alpha} = 1\, .
\]

Our initial target is to show the $L^{2^*}_k$-compactness of $(U_\alpha)$. From the equality above and from the continuity of $\mathcal{A}_0(n,F)$ on $F$ and $\mathcal{B}_0(n,F,G,g)$ on $(F,G,g)$, we easily deduce that $(U_\alpha)$ is bounded in $H_k^{1,2}(M)$. In particular, $(U_\alpha)$ converges weakly in $H_k^{1,2}(M)$ and strongly in $L^q_k(M)$ to $U_0$ for any $1 \leq q < 2^*$. We first discard the case when $U_0 = 0$. Note that the maps $U_\alpha$ minimize the non-smooth functionals

\begin{equation} \label{E.53}
J_\alpha(U) = \mathcal{A}_0(n,F_\alpha) \int_M |\nabla_{g_\alpha} U|^2\; dv_{g_\alpha} + \mathcal{B}_0(n,F_\alpha,G_\alpha,g_\alpha) \int_M G_\alpha(x,U)\; dv_{g_\alpha}
\end{equation}

\n on $\Lambda_\alpha = \{U \in H^{1,2}_k(M) :\, \int_M F_\alpha(U) dv_{g_\alpha} = 1\}$. In particular, $t = 0$ is a minimum point of

\[
f(t) := J_\alpha (\frac{(1 + t \varphi) U_\alpha}{\left( \int_M F_\alpha((1 + t \varphi)U_\alpha) dv_{g_\alpha} \right)^{1/2^*}})\, .
\]

\n Using now the $2^*$-homogeneity of $F_\alpha$ and the $2$-homogeneity of $G_\alpha$, we can differentiate at $t = 0$. Straightforward computation furnishes

\begin{eqnarray}
0 = f'(0) &=& \int_M \nabla_{g_\alpha} U_\alpha \cdot \nabla_{g_\alpha} (\varphi U_\alpha)\; dv_{g_\alpha} + \mathcal{B}_0(n,F_\alpha,G_\alpha,g_\alpha) \int_M G_\alpha(x,U_\alpha) \varphi\; dv_{g_\alpha} \label{E.55}\\
 && - \int_M F_\alpha(U_\alpha) \varphi\; dv_{g_\alpha} \nonumber
\end{eqnarray}

\n for all $\varphi \in C^1(M)$. Mimicking the proof of Proposition \ref{P.5}, we also have

\begin{equation} \label{E.54}
\int_{M} \nabla_{g_\alpha} |U_\alpha| \cdot \nabla_{g_\alpha} \varphi\; dv_{g_\alpha} \leq \int_M |U_\alpha|^{-1} F_\alpha(U) \varphi\; dv_{g_\alpha}
\end{equation}

\n for all nonnegative function $\varphi \in C^1(M)$.

If $U_0 = 0$, thanks to (\ref{E.54}), we can develop a concentration study on $(U_\alpha)$ in a similar spirit as in the proof of Theorem \ref{Teo.1} without using maximum points of $|U_\alpha|$. Precisely, Lemmas \ref{L.3}, \ref{L.4} and \ref{L.5} hold for $(U_\alpha)$ in the sense following: the sequence $(U_\alpha)$ concentrates at a unique point $x_0$, $(|U_\alpha|)$ converges to 0 in $L^\infty_{loc}(M \setminus \{x_0\})$ as $\alpha \rightarrow + \infty$ and, for any $\delta > 0$,

\begin{equation} \label{E.56}
\lim \limits_{\alpha \rightarrow + \infty} \frac{\int_{M \setminus B_{g_\alpha}(x_0,\delta)} |U_{\alpha}|^2 \; dv_{g_\alpha}}{\int_M |U_{\alpha}|^2 \; dv_{g_\alpha}} = 0\, .
\end{equation}

\n We now argue as in the conclusion of the proof of Theorem \ref{Teo.2}. Let $0 < \varepsilon < \varepsilon_0$ fixed and $r_0 > 0$ be the radius provided in Proposition \ref{P.7} when applied to $\varepsilon$, $F_\alpha$, $G_\alpha$ and $g_\alpha$. Consider a cutoff function $\varphi_\alpha \in C^\infty(B_{g_\alpha}(x_0, r_0))$ with $0 \leq \varphi_\alpha \leq 1$, $\varphi_\alpha = 1$ in $B_{g_\alpha}(x_0, r_0/2)$ and $|\nabla_{g_\alpha} \varphi_\alpha| \leq C_0$, where $C_0 > 0$ does not depend on $\alpha$. Choosing $\varphi = \varphi_\alpha^2$ in (\ref{E.55}) and repeating some simple manipulations, we easily check that

\[
\left( \mathcal{B}_0(n,F_\alpha,G_\alpha,g_\alpha) - \mathcal{B}_\varepsilon(F_\alpha,G_\alpha,g_\alpha)\right) \int_M G_{\alpha}(x,U_{\alpha}) \varphi_\alpha^2\; dv_{g_\alpha} \leq \mathcal{A}_0(n,F_{\alpha}) \int_M |\nabla_{g_\alpha} \varphi_\alpha|^2 |U_\alpha|^2\; dv_{g_\alpha}\, .
\]

\n By (\ref{E.50}) and the definition of $\mathcal{B}_\varepsilon(F_\alpha,G_\alpha,g_\alpha)$, we arrive at

\[
\frac{\int_{M \setminus B_{g_\alpha}(x_\alpha, r_0/2)} |U_\alpha|^2\; dv_{g_\alpha}}{\int_{M} |U_{\alpha}|^2\; dv_{g_\alpha}} \geq C_1
\]

\n for some constant $C_1 > 0$ and $\alpha > 0$ large, which contradicts (\ref{E.56}). Thus, $U_0 \neq 0$. Again, thanks to (\ref{E.55}) and Theorem \ref{Teo.1}, we can proceed in exactly the same way as in the proof of Theorem \ref{Teo.1} with the aid of Propositions \ref{B-L0} and \ref{C-P1}, in order to establish the compactness of $(U_\alpha)$ in $L^{2^*}_k(M)$.

Suppose further that $((F_\alpha,G_\alpha,g_\alpha))$ converges to $(F,G,g)$ in $C^1_{loc}(\R^k) \times C^0(M, C^1_{loc}(\R^k)) \times {\cal M}^n$. In this case, the maps $U_\alpha$ are critical points of the functionals $J_\alpha$. In other words, these maps satisfy the systems

\[
- \mathcal{A}_0(n,F_{\alpha}) \Delta_{g_{\alpha}} u_{\alpha}^i + \frac{1}{2} \mathcal{B}_0(n,F_\alpha,G_\alpha,g_\alpha) \frac{\partial G_\alpha}{\partial t_i}(x, U_{\alpha}) = \frac{1}{2^*} \frac{\partial F_\alpha}{\partial t_i}(U_{\alpha})
\]

\n and so are of class $C^1$, by Proposition \ref{Reg}. Moreover, by standard elliptic estimates, the $C_k^0$-compactness of $(U_\alpha)$ clearly follows from a uniform estimate for the sequence $(|U_\alpha|)$. In what follows, we focus our attention on this last point.

Let $x_\alpha \in M$ be a maximum point of $|U_\alpha|$ and assume by contradiction that $(|U_\alpha(x_\alpha)|)$ blows up as $\alpha \rightarrow + \infty$. Consider exponential charts $\exp_{x_{\alpha}}$ centered at $x_{\alpha}$ with respect to the metrics $g_\alpha$ and let $\delta > 0$ be a small number such that $\exp_{x_{\alpha}}$ are diffeomorphisms from $B(0,\delta) \subset \R^n$ onto $B_{g_\alpha}(x_{\alpha},\delta)$ for any $\alpha > 0$ large. Define now the metrics $h_{\alpha}$ and the maps $V_\alpha$ on the open ball $\Omega_\alpha = B(0, \mu_{\alpha}^{-1} \delta)$ by

\[
h_{\alpha}(x) = \left( \exp_{x_{\alpha}}^*g_\alpha \right)(\mu_{\alpha}x)
\]

\n and

\[
V_{\alpha}(x) = \mu_{\alpha}^{n/2^*} U_{\alpha} \left( \exp_{x_{\alpha}}(\mu _{\alpha}x) \right)\, ,
\]

\n where

\[
\mu_\alpha = |U_\alpha(x_\alpha)|^{-2^*/n}\, .
\]

\n Clearly, the maps $V_{\alpha} = (v_{\alpha}^1, \ldots, v_{\alpha}^k)$ satisfy the systems

\[
- \mathcal{A}_0(n,F_{\alpha}) \Delta_{h_{\alpha}} v_{\alpha}^i + \frac{\mu_{\alpha}^2}{2} \mathcal{B}_0(n,F_\alpha,G_\alpha,g_\alpha) \frac{\partial G_\alpha}{\partial t_i}(\exp_{x_{\alpha}}(\mu_{\alpha}x), V_{\alpha}) = \frac{1}{2^*} \frac{\partial F_\alpha}{\partial t_i}(V_{\alpha})\ \ {\rm on}\ B(0, \mu_{\alpha}^{-1} \delta)\, .
\]

\n Some facts that deserve attention are: $(h_{\alpha})_{ij}$ converges to $\xi_{ij}$ in $C^1_{loc}(\R^n)$ and $\mathcal{A}_0(n,F_{\alpha})$, $\mathcal{B}_0(n,F_\alpha,G_\alpha,g_\alpha)$, $|V_\alpha|$, $\frac{\partial G_\alpha}{\partial t_i}(\exp_{x_{\alpha}}(\mu_{\alpha}x), V_{\alpha})$ and $\frac{\partial F_\alpha}{\partial t_i}(V_{\alpha})$ are bounded for $\alpha > 0$ large. So, usual elliptic estimates applied to the above system imply that $V_\alpha$ converges to $V_0$ in $C^1_{k,loc}(\R^n)$. Letting then $\alpha \rightarrow +\infty$ in the system above, we derive

\begin{equation} \label{E.58}
- \mathcal{A}_0(n,F) \Delta v_0^i = \frac{1}{2^*} \frac{\partial F}{\partial t_i}(V_0)\ \ {\rm on}\ \R^n\, .
\end{equation}

\n In addition, we have $|V_0(0)| = 1$ since $|V_\alpha(0)| = 1$ for all $\alpha$. On the other hand, since $(U_\alpha)$ converges to $U_0$ in $L^{2^*}_k(M)$, there exists a nonnegative function $f_0 \in L^{2^*}(M)$ such that $|U_\alpha| \leq f_0$ on $M$ for $\alpha > 0$ large. Consequently, for any $R > 0$,

\[
\int_{B(0,R)} |V_0|^{2^*}\; dx = \lim_{\alpha \rightarrow + \infty} \int_{B(0,R)} |V_\alpha|^{2^*}\; dv_{h_{\alpha}} = \lim_{\alpha \rightarrow + \infty} \int_{B_{g_\alpha}(x_\alpha, R\mu_{\alpha})} |U_\alpha|^{2^*}\; dv_{g_{\alpha}}
\]

\[
\leq \lim_{\alpha \rightarrow + \infty} \int_{B_{g_\alpha}(x_\alpha, R\mu_{\alpha})} f_0^{2^*}\; dv_{g_{\alpha}} = 0\, ,
\]

\n so that $V_0 = 0$, which contradicts the fact that $|V_0(0)| = 1$. This finally completes the proof.\bl \\

\section{Examples and counter-examples}

This section is devoted to some examples of existence and
non-existence of extremal maps, lack of compactness of the set
$\mathcal{E}_k(F,G,g)$ and the lost of continuity of the second best
constant. We start this section given bounds for the second best
constant $\mathcal{B}_0(n,F,G,g)$ in relation to the scalar one
$B_0(n,1,g)$. These bounds are used in the building of examples and
counter-examples.

\begin{propo}\label{Ex}
Let $(M,g)$ be a compact Riemannian manifold of dimension $n \geq 3$
and let $t_0 \in \mathbb{S}^{k-1}$ be such that $F(t_0)=M_F$. Then,

\begin{equation}\label{Des1}
\frac{M_F^{2/2^*}B_0(n,1,g)}{\max_{x\in M}G(x,t_0)} \leq
\mathcal{B}_0(n,F,G,g) \leq \frac{M_F^{2/2^*}B_0(n,1,g)}{m_G}
\end{equation}

\n where

\[
m_G = \min \limits_{M \times \mathbb{S}^{k-1}} G\; .
\]

\n In particular, if $t_0 \in \mathbb{S}^{k-1}$ is such that
$F(t_0)=M_F$ and $m_{G} = \max \limits_{x\in M}G(x,t_0)$, then

\[
\mathcal{B}_0(n,F,G,g) = \frac{M_F^{2/2^*}B_0(n,1,g)}{m_G}\,.
\]

\end{propo}

\n {\bf Proof of Proposition \ref{Ex}.}  Proposition \ref{P.2}
implies that

\[
\left(\int_{M} F(U)\; dv_g\right)^{\frac{2}{2^*}} \leq \mathcal{A}_0(n,F) \int_{M} |\nabla_g U|^2 \; dv_g  + \mathcal{B} _0(2,F,G,g)
\int_{M} G(x,U)\; dv_g\,,
\]

\n for all $U \in H^{1,2}_k(M)$. Hence, taking $U = u t_0$ with $u \in
H^{1,2}(M)$, we obtain

\[
\left(\int_{M} |u|^{2^*}\; dv_g\right)^{\frac{2}{2^*}} \leq A_0(n)
\int_{M} |\nabla_g u|^2 \; dv_g + \mathcal{B}
_0(n,F,G,g)M_F^{-2/2^*}\max _{x\in M}G(x,t_0) \int_{M} |u|^2\;
dv_g\; .
\]

\n From definition of $B_0(n,1,g)$, we find

\begin{equation}\label{Des2}
\mathcal{B}_0(n,F,G,g) \geq \frac{M_F^{2/2^*}B_0(n,1,g)}{\max_{x \in
M} G(x,t_0)}\;.
\end{equation}

\n Another hand, from the proof of Proposition \ref{P.2}, we have

\[
\left(\int_{M} F(U)\; dv_g\right)^{\frac{2}{2^*}} \leq
\mathcal{A}_0(n,F) \int_{M} |\nabla_g U|^2 \; dv_g +
\frac{M_F^{2/2^*}B_0(n,1,g)}{m_G}\int_{M} G(x,U)\; dv_g
\]

\n for all $U \in H^{1,2}_k(M)$. Then, from definition of $
\mathcal{B}_0(n,F,G,g)$,

\begin{equation}\label{Des3}
\mathcal{B} _0(n,F,G,g)\leq \frac{M_F^{2/2^*}B_0(n,1,g)}{m_G}\;.
\end{equation}

\n Combining (\ref{Des2}) and (\ref{Des3}), we find (\ref{Des1}). \bl

The first examples are on the existence and non-existence of
extremal maps.

\n {\bf Example 1.} Let $(M,g)$ be a compact Riemannian manifold of
dimension $n\geq 3$ such that

\[
B_0(n,1,g) = \frac{n-2}{4(n-1)} A_0(n) \max \limits_{M} S_g
\]

\n and (\ref{B-opt}) admit an extremal function. Consider the function
$G:M \times \mathbb{R}^k \rightarrow \mathbb{R}$, $G(x,t) =
\sum_{i,j = 1}^k A_{ij}(x) |t_i| |t_j|$, where the functions
$A_{ij}$'s are nonnegative and continuous. Assume also that, for
some $i_0$, $A_{i_0i_0}
> 0$ does not depends of $x$ and $A_{ii} \geq A_{i_0i_0}$ for all $i$. Clearly,

\[
A_{i_0i_0}|t|^2 \leq \sum _{i=1}^k A_{ii}(x)|t_i|^2 \leq
\sum_{i,j}A_{ij}(x) |t_i| |t_j|\;.
\]

\n Thus, $m_G = A_{i_0i_0}$. Let $F:\mathbb{R}^k \rightarrow
\mathbb{R}$ be a $2^*$-homogeneous, continuous and positive function
such that $F(e_{i_0}) = M_F$, where $e_{i_0}$ is the $i_0$-element
of the standard bases of $\R^k$ . Hence, from the Proposition
\ref{Ex},

\[
\mathcal{B}_0(n,F,G,g) = \frac{M_F^{2/2^*}B_0(n,1,g)}{A_{i_0i_0}}\;.
\]

\n If $u_0\in H^{1,2}(M)$ is an extremal function of (\ref{B-opt}),
then $U=u_0e_{i_0}$ is an extremal map of (\ref{B-opt-v}). Then,
this example shows that $(I)^*$ and $(II)^*$, in the Theorem
\ref{Teo.2}, can happen in the same time. \bl

\n {\bf Example 2.} Let $(M,g)$ be a compact Riemannian manifold of
dimension $n\geq 3$. Consider a $2$-homogeneous, continuous and
positive $G:M \times \mathbb{R}^k \rightarrow \mathbb{R}$ that does
not depends of $x$. Choose $t_0 \in \mathbb{S}^{k-1}_2$ such that
$m_G = G(t_0)$, and a linear injective transformation $A:
\mathbb{R}^k \rightarrow \mathbb{R}^k$ such that $A(t_1) = t_0$,
where $t_1 \in \mathbb{S}^{k-1}$ satisfies $F(t_1) = M_F$. Hence,
from Proposition \ref{Ex}, the inequality (\ref{B-opt-v}) with
$F\circ A$, possesses an extremal map, if (\ref{B-opt}) possesses an
extremal function. \bl

The following example is on the non-existence of extremal maps.

\n {\bf Example 3.} Let $(M,g)$ be a compact Riemannian manifold of
dimension $n\geq 3$ such that

\[
B_0(n,1,g) = \frac{n-2}{4(n-1)} A_0(n) \max \limits_{M} S_g
\]

\n and (\ref{B-opt}) possesses an extremal function. Consider the
function $G:M \times \mathbb{R}^k \rightarrow \mathbb{R}$, $G(x,t) =
\sum_{i,j = 1}^k A_{ij}(x) |t_i| |t_j|$, where the functions
$A_{ij}$'s are nonnegative and continuous. Assume also that, for
some $i_0$, $A_{i_0i_0}
> 0$ does not depends of $x$ and $A_{ii} \geq A_{i_0i_0}$ for all $i$. Clearly,

\[
A_{i_0i_0}|t|^2 \leq \sum _{i=1}^k A_{ii}(x)|t_i|^2 \leq
\sum_{i,j}A_{ij}(x) |t_i| |t_j|\;.
\]

\n Thus, $m_G = A_{i_0i_0}$. Let $F:\mathbb{R}^k \rightarrow
\mathbb{R}$ be a $2^*$-homogeneous, continuous and positive function
such that $F(e_{i_0}) = M_F$, where $e_{i_0}$ is the $i_0$-element
of the standard bases of $\R^k$ . Hence, from the Proposition
\ref{Ex},

\[
\mathcal{B}_0(n,F,G,g) = \frac{M_F^{2/2^*}B_0(n,1,g)}{A_{i_0i_0}}\;.
\]

\n Suppose, by contradiction, that exists an extremal map $U_0$ for the
inequality (\ref{B-opt-v}). Hence,

\[
\mathcal{B} _0(n,F,G,g) \int_{M} \sum \limits
_{i,j=1}^{k}A_{ij}|u_0^i| |u_0^j|\; dv_g=\left(\int_{M} F(U_0)\;
dv_g\right)^{\frac{2}{2^*}}  - \mathcal{A}_0(n,F) \int_{M} |\nabla_g
U_0|^2 \; dv_g
\]

\begin{eqnarray*}
&\leq&  M_F^{2/2^*}\sum \limits _{i=1}^{k}\left(\int _M|u_0^i|^{2^*}\;
dv_g\right)^{\frac{2}{2^*}} - \mathcal{A}_0(n,F) \int_{M} |\nabla_g
U_0|^2 \; dv_g\\
&\leq& \mathcal{A}_0(n,F) \sum \limits _{i=1}^{k}\int _M|\nabla
u_0^i|^{2}\; dv_g + M_F^{2/2^*} B_0(n,1,g)\sum \limits
_{i=1}^{k}\int _M|u_0^i|^{2}\; dv_g\\
&&-\mathcal{A}_0(n,F) \int_{M} |\nabla_g U_0|^2 \; dv_g \leq
\frac{M_F^{2/2^*} B_0(n,1,g)}{A_{i_0i_0}}\int _M\sum \limits
_{i,j=1}^{k}A_{ij}|u_0^i| |u_0^j|\; dv_g\,,
\end{eqnarray*}

\n since

\[
A_{i_0i_0} \sum \limits _{i=1}^{k}|u_0^i|^2 \leq  \sum \limits
_{i,j=1}^{k}A_{ij}|u_0^i| |u_0^j| \,.
\]

\n This implies that

\begin{equation}\label{a1}
\sum \limits _{i=1}^{k}\left(\int _M|u_0^i|^{2^*}\;
dv_g\right)^{\frac{2}{2^*}}= A_0(n) \sum \limits _{i=1}^{k}\int_{M}
|\nabla_g u_0^i|^2 \; dv_g + B _0(n,1,g) \int_{M} \sum \limits
_{i=1}^{k}|u_0^i|^2\; dv_g\,.
\end{equation}

\n Independently, it follows, from the inequality (\ref{B-opt}), that

\begin{equation}\label{b1}
\left(\int _M|u_0^i|^{2^*}\; dv_g\right)^{\frac{2}{2^*}}\leq
A_0(n)\int_{M} |\nabla_g u_0^i|^2 \; dv_g + B _0(n,1,g)
\int_{M}|u_0^i|^2\; dv_g\,,
\end{equation}

\n for all $i$. Then, from (\ref{a1}) and (\ref{b1}), there exists
$j\in \{ 1,...,k\}$ such that $u_0^j\neq 0$ and

\[
\left(\int _M|u_0^j|^{2^*}\; dv_g\right)^{\frac{2}{2^*}}=
A_0(n)\int_{M} |\nabla_g u_0^j|^2 \; dv_g + B_0(n,1,g)
\int_{M}|u_0^j|^2\; dv_g\,.
\]

\n This is not possible, since (\ref{B-opt}) does not possesses
extremal function. \bl

Next we give an example on lack of compactness of the set
$\mathcal{E}_k(F,G,g)$.

\n {\bf Example 4.} Consider the standard sphere $(\mathbb{S}^n,h)$.
Note that

\[
B_0(n,1,h)= \frac{n-2}{4(n-1)}A_0(n)\max
\limits_{\mathbb{S}^n}S_h\,.
 \]

\n Given $p_0\in \mathbb{S}^n$ and $\beta >1$ real, define the function

\[
 u_{p_0,\beta }(p)=(\beta^2-1)^{\frac{n-2}{4}}\omega
 _n^{-\frac{1}{2^*}}\left( \beta - \cos r(p)   \right)^{1-\frac{n}{2}},
\]

\n where $r(p)=d_{h}(p,p_0)$. It is easily to see that $u_{p_0,\beta }$
is an extremal function for the inequality (\ref{B-opt}) and
$||u_{p_0,\beta }||_{2^*}=1$. Observe that

\[
u_{p_0,\beta }(p_0)\rightarrow \infty
\]

\n as $\beta\rightarrow 1$. Thus, if we choose $F$ and $G$ as in the
Example 1 such that $M_F=1$, the compactness does not holds. Just
use the sequence of maps $U_{\beta}=u_{p_0,\beta }t_0$. This lack of
compactness is possible because here we have

\[
A_{i_0i_0}\mathcal{B}_0(n,F,G,h) = \frac{n-2}{4(n-1)}
\mathcal{A}_0(n,F) \max \limits_{M} S_h\,.
\]

\bl

Next we give an example on the lost of continuity of the second best
constant.

\n {\bf Example 5.}  Let $(M,g)$ be a smooth compact Riemannian
manifold of dimension $n\geq 4$. Consider a sequence
$(f_\alpha)_\alpha \subset C^{\infty }(M)$ of positive functions
converging to the constant function $f_0 = 1$ in $L^p(M)$, $p>n$,
such that $\max_{M}f_\alpha \rightarrow + \infty$. Let $u_\alpha\in
C^{\infty }(M)$, $u_\alpha>0$, be the unique solution of

\[
- \frac{4(n-1)}{n-2} \Delta _{g}u + u = f_\alpha \ .
\]

\n From the classical elliptic $L^p$ theory, it follows that
$(u_\alpha)_\alpha$ is bounded in $H^{p}_2(M)$, where $H^{p}_2(M)$
stands for the second order $L^p$-Sobolev space on $M$, so that
$u_\alpha$ converges to $u_0$ in $C^{1,\beta}(M)$ for some $0 <
\beta < 1$. Moreover, $u_0=1$, since $f_\alpha$ converges to $1$ in
$L^p(M)$ and the constant function $1$ is the unique solution of the
limit problem. Therefore, $g_\alpha = u_\alpha^{2^*-2}g$ is a smooth
Riemannian metric converging to $g$ in the $C^{1,\beta}$-topology.
Note also that there exists a constant $c > 0$, independent of
$\alpha$, such that

\[
S_{g_\alpha} = (- \frac{4(n-1)}{n-2} \Delta _{g}u_\alpha + S_{g}
u_\alpha ) u_\alpha^{1-2^*} \geq f_\alpha u_\alpha^{1-2^*} -
cu_\alpha^{2-2^*},
\]

\n so that $\max_{M}S_{g_\alpha} \rightarrow + \infty$. On the other
hand, for $n \geq 4$, we have the lower bound

\[
B_0(n,1,g_\alpha) \geq \frac{n-2}{4(n-1)} A_0(n)
\max_{M}S_{g_\alpha},
\]

\n so that $B_0(n,1,g_\alpha)\rightarrow + \infty$. In particular,
$B_0(n,1,g_\alpha) \not \rightarrow B_0(n,1,g)$. Thus, if we choose
$F$ and $G$ as in the Example 1, we have

\[
A_{i_0i_0}\mathcal{B}_0(n,F,G,g_{\alpha}) \geq \frac{n-2}{4(n-1)}
\mathcal{A}_0(n,F) \max \limits_{M} S_{g_{\alpha}}
\]

\n and, consequently, $\mathcal{B}_0(n,F,G,g_{\alpha}) \not \rightarrow
\mathcal{B}_0(n,F,G,g)$.

\bl

\n {\bf Acknowledgments:}  Both authors were partially supported by CNPq and Fapemig. The authors are indebted to Prof. Hebey for his valuable comments
concerning this work.

\appendix

\section{Appendix}

\subsection {A generalized Br\'{e}zis-Lieb lemma}

We will prove a vector version of the Br\'{e}zis-Lieb Lemma \cite{BrLi} with low regularity. Indeed, below we assume no regularity conditions on the
function $F$ besides continuity. The result is precisely

\begin{propo}\label{B-L0}
Let $(M,g)$ be a smooth Riemannian manifold of dimension $n \geq 2$ and $F:\R^k\rightarrow \R$ be a
$p$-homogeneous, positive and continuous function with $0 < p < \infty$. Suppose that $(\int_M |U_{\alpha}|^pdv_g)$ is a bounded sequence and $U_{\alpha} \rightarrow U$ a.e. on $M$. Then

\[
\lim _{n\rightarrow \infty}\int _M\left(F(U_n) - F(U_n -U)
\right)dv_g=\int _MF(U)dv_g\;.
\]

\end{propo}

\n {\bf Proof of Proposition \ref{B-L0}.} First, we claim that given $\varepsilon>0$ there exists a constant $C(\varepsilon)$ such that for all $t$, $s\in\R^k$, we have

\begin{equation}\label{r2}
\left|F(s+t)-F(s)\right|\leq\varepsilon|s|^p+C(\varepsilon)|t|^p.
\end{equation}

\n Indeed, the continuity of $F$ implies that for each $\varepsilon>0$ there exists $\delta(\varepsilon)>0$ such that

\[
\left|F(s+t)-F(s)\right|\leq\varepsilon, \quad \forall|t|<\delta(\varepsilon)\,,
\]

\n where we restrict $F$ to the ball $B_2(0)$. Without loss of generality we can assume
$\delta(\varepsilon)<1$. Define $M=\max\limits_{\overline{B_2}(0)}F$ and $C(\varepsilon)=\frac{M}{\delta(\varepsilon)^p}$.
If $\frac{|s|}{|t|}\leq1$,

\begin{eqnarray*}
\left|F\left(\frac{s}{|t|}+\frac{t}{|t|}\right)-F\left(\frac{s}{|t|}\right)\right|&\leq& M\\
&\leq&\frac{M}{\delta(\varepsilon)^p}+\varepsilon\frac{|s|^p}{|t|^p}\,.
\end{eqnarray*}

\n If $\delta(\varepsilon)\leq\frac{|t|}{|s|}\leq1$,

\begin{eqnarray*}
\left|F\left(\frac{s}{|s|}+\frac{t}{|s|}\right)-F\left(\frac{s}{|s|}\right)\right|&\leq&M\\
&=&\frac{M}{\delta(\varepsilon)^p}\delta(\varepsilon)^p\\
&\leq&\frac{M}{\delta(\varepsilon)^p}\frac{|t|^p}{|s|^p}+\varepsilon\,.
\end{eqnarray*}

\n If $\frac{|t|}{|s|}\leq\delta(\varepsilon)$,

\begin{eqnarray*}
\left|F\left(\frac{s}{|s|}+\frac{t}{|s|}\right)-F\left(\frac{s}{|s|}\right)\right|&\leq&\varepsilon\\
&\leq&\varepsilon+\frac{M}{\delta(\varepsilon)^p}\frac{|t|^p}{|s|^p}\,.
\end{eqnarray*}

\n Thus, from the $p$-homogeneity of $F$ we obtain (\ref{r2}). Now, we define $V_{\alpha}=U_{\alpha}-U$ and the functional

\[
H_{\varepsilon}^{\alpha}(x)=(|F(U_{\alpha}(x))-F(V_{\alpha}(x))-F(U(x))|-\varepsilon|V_{\alpha}(x)|^p)_+\,,
\]

\n where $(u)_+=\max(u,0)$. Observe that $V_{\alpha}\rightarrow0$ a.e. on $M$, and $F(0)=0$: if $\lambda >0$ and $x\in\R^k$, we have:

\[
F(0)=\lim_{\lambda\rightarrow 0^+}F(\lambda x)=\lim_{n\rightarrow 0}\lambda^pF(x)=0\,.
\]

\n Hence, if $\alpha\rightarrow\infty$, we obtain

\[
F(V_{\alpha}(x))\rightarrow F(0)=0\,,
\]

\n since $F$ is a continuous function, and

\[
H_{\varepsilon}^{\alpha}\rightarrow0
\]

\n a.e. on $M$. Another hand,

\[
|F(U_{\alpha})-F(V_{\alpha})-F(U)|\leq|F(U_{\alpha})-F(V_{\alpha})|+|F(U)|\,.
\]

\n Putting $s=V_{\alpha}$ and $t=U$ in (\ref{r2}), we find

\[
\left|F(U_{\alpha})-F(V_{\alpha})\right|=\left|F(s+t)-F(s)\right|\leq \varepsilon|V_{\alpha}|^p+C(\varepsilon)|U|^p\,.
\]

\n Thus,

\[
\left|F(U_{\alpha})-F(V_{\alpha})-F(U)\right|\leq\varepsilon|V_{\alpha}|^p+C(\varepsilon)|U|^p+|F(U)|
\]

\n and consequently

\[
H_{\varepsilon}^{\alpha}\leq C(\varepsilon)|U|^p+|F(U)|\,.
\]

\n Since $(C(\varepsilon)|U|^p+|F(U)|)\in L^1(M)$, it follows from the Lebesgue's convergence Theorem that

\[
\int_M H_{\varepsilon}^{\alpha} dv_g\rightarrow0
\]

\n as $\alpha\rightarrow\infty$. From definition of $H_{\varepsilon}^{\alpha}$, we have

\[
H_{\varepsilon}^{\alpha}(x)\geq |F(U_{\alpha}(x))-F(V_{\alpha}(x))-F(U(x))|-\varepsilon|V_{\alpha}(x)|^p\; ,
\]

\n so that

\[
\left|F(U_{\alpha})-F(V_{\alpha})-F(U)\right|\leq\varepsilon|V_{\alpha}|^p+H_{\varepsilon}^{\alpha}\,.
\]

\n Hence,

\[
I_{\alpha}=\int_M \left|F(U_{\alpha})-F(V_{\alpha})-F(U)\right|dx\leq\int_M(\varepsilon|V_{\alpha}|^p+
H_{\varepsilon}^{\alpha})dv_g
\]

\n and

\[
\lim_{\alpha\rightarrow\infty}I_{\alpha}\leq\lim_{\alpha\rightarrow\infty}\int_M\varepsilon|V_{\alpha}|^pdv_g+
\lim_{\alpha\rightarrow\infty}\int_M H_{\varepsilon}^{\alpha} dv_g=\varepsilon\lim_{\alpha\rightarrow\infty}\int_M|V_{\alpha}|^pdv_g= \varepsilon \widetilde{C}\,.
\]

\n Note that the constant $\widetilde{C}$ comes from the limitation
of the sequence $(U_{\alpha})$ in $L^p_k(M)$. Letting
$\varepsilon\rightarrow0$ we obtain the desired conclusion. \bl

\subsection{ A concentration-compactness principle}

In this section, we prove a version of the Lions's
concentration-compactness principle \cite{Li} for vector valued maps. Such a
result is an application of the generalized Br\'{e}zis-Lieb lemma proved
in the preceeding section. Let $((F_\alpha,G_\alpha,g_\alpha))$ be a
sequence converging to $(F,G,g)$ in $\mathcal{F}_k  \times
\mathcal{G}_k \times {\cal M}^n$.

\begin{propo}\label{C-P1}
Let $(M,g)$ be a smooth compact Riemannian manifold of dimension $n \geq
3$ and $((F_\alpha, G_\alpha))$ be a sequence converging to $(F,G)$ in ${\cal F}_k \times {\cal G}_k$. If $(U_{\alpha})_{\alpha} \subset H^{1,2}(M)$ is a sequence such
that

\[
U_{\alpha} \rightharpoonup U \ \mbox{in} \ H^{1,2}(M),
\]

\[
|\nabla_g U_{\alpha}|^2 \; dv_g \rightharpoonup \mu,
\]

\[
F_{\alpha}(U_{\alpha}) \; dv_g \rightharpoonup \nu
\]

\n as $\alpha \rightarrow \infty$, where $\mu$ and $\nu$ are bounded
nonnegative measures, then there exist at most a countable set
$\{x_j\}_{j \in {\cal T}}$ and positive numbers $\{\mu_j\}_{j \in
{\cal T}}$ e $\{\nu_j\}_{j \in {\cal T}}$ such that

\[
\mu \geq |\nabla _gU|^2 \; dv_g + \sum_{j \in {\cal T}} \mu_j
\delta_{x_j}, \ \ \nu = F(U) \; dv_g + \sum_{j \in {\cal T}} \nu_j
\delta_{x_j}
\]

\n with $\mathcal{A}_0(n,F) \mu_j \geq \nu_j^{2/2^*}$ for all $j \in
{\cal T}$, where $\delta_{x_j}$ denotes the Dirac mass centered at
$x_j$.
\end{propo}

\n {\bf Proof of Proposition \ref{C-P1}.}  Set $W_{\alpha} =
U_{\alpha} - U$, so that $W_{\alpha} \rightharpoonup 0$ in
$H^{1,2}(M)$. Define $\theta = \nu - F(U)\; dv_g$. By the
generalized Brezis-Lieb lemma,

\[
F(W_{\alpha}) dv_g \rightharpoonup \theta\ .
\]

\n In addition, up to a subsequence, we have

\[
|\nabla W_{\alpha}|_g^2 \; dv_g \rightharpoonup \lambda
\]

\n for some bounded nonnegative measure $\lambda$. We have only to
show that there hold reverse H\"{o}lder inequalities for the measure
$\theta$ with respect to $\lambda$. The rest of the proof is
standard. Taking $V = \varphi W_{\alpha}$, where $\varphi \in
C^1(M)$, one finds

\[
\left( \int_M |\varphi|^{2^*} F(W_{\alpha})\; dv_g
\right)^{2/2^*} \leq \mathcal{A}_0(n,F)\int_M |\nabla
_g(\varphi W_{\alpha})|^2 \; dv_g +C \int_M |\varphi
W_{\alpha}|^{2}\; dv_g
\]

\[
\leq \mathcal{A}_0(n,F) \int_M |\varphi|^2 |\nabla _g W_{\alpha}|^2
\; dv_g + C \sum _{i=1}^k\int_M |\varphi \nabla _gW^i_{\alpha}|
|W^i_{\alpha} \nabla _g\varphi| \; dv_g
\]

\[
+ C \int_M |\nabla _g\varphi|^2|W_{\alpha}|^2  \; dv_g + D \int_M
|\varphi|^2| W_{\alpha}|^{2}\; dv_g
\]

\[
\leq \mathcal{A}_0(n,F) \int_M |\varphi|^2 |\nabla _gW_{\alpha}|^2
\; dv_g + \tilde{C} \max_M \left\{ |\nabla _g \varphi|^2 +
|\varphi|^2 \right\} \int_M |W_{\alpha}|^{2}\; dv_g\ .
\]

\n Since $||W_{\alpha}||_2 \rightarrow 0$ as $\alpha \rightarrow
\infty$, it follows from the previous inequality that

\[
\left( \int_M |\varphi|^{2^*}\; d\theta \right)^{2/2^*} \leq
\mathcal{A}_0(n,F) \int_M |\varphi|^{2}\; d\lambda\,.
\]

\n Thus,

\[
\left( \int_M |\varphi|^{2^*}\; d\theta \right)^{2/2^*} \leq
\mathcal{A}_0(n,F) \int_M |\varphi|^{2}\; d\lambda\,.
\]

\n Then, for all $\Phi \in C_k^1(M)$, the reverse
H\"{o}lder inequality

\[
\left( \int_{M} F(\Phi)\; d\theta \right)^{2/2^*} \leq
\mathcal{A}_0(n,F) \int_{M} F(\Phi)^{2/2^*}\; d\lambda
\]

\n holds. \bl

\subsection{Minimizing solutions}

Here we assume that the potential functions $F$ and $G$ are of class
$C^1$. Consider the functional $I:H^{1,2}_k(M)\rightarrow \R$ given
by

\[
I(U)=\int _M|\nabla U|^2dv_g + \int _MG(x,U)dv_g
\]

\n for all $U\in H^{1,2}_k(M)$. Define

\[
\mu_{g,k}=\inf _{U\in S^{1,2}_k(M)}I(U)\,,
\]

\n where

\[
E^k_F=\lbrace U\in H^{1,2}_k(M): \; \int _MF(U)dv_g=1 \rbrace
\;.
\]

\n This section is devoted to the systems

\begin{gather} \label{S}
-\Delta_g u_i + \frac{\lambda}{2} \frac{\partial G(x,U)}{\partial t_i} = \frac{1}{2^*} \frac{\partial F(U)}{\partial t_i} \ \ {\rm on}\ \ M,\ \ i=1,\ldots,k \tag{$S$}
\end{gather}

\n We give a sufficient condition for the existence of a solution to
this system.

\begin{propo}\label{exist}
Assume that $\mu_{g,k}<\mathcal{A}_0(n,F)^{-1}$. Then there exists a
minimizer $U\in E^k_F$ for  $\mu_{g,k}$. Namely, there
exists $U\in E^k_F$ such that $I(U)=\mu_{g,k}$. In
particular, $U=(u_1,...,u_k)$ is a solution of the system (\ref{S}) with
$\lambda=\mu_{g,k}$.
\end{propo}

\n {\bf Proof of Proposition \ref{exist}.}: Let $(U_{\alpha
})_{\alpha }$ be a minimizing sequence for $\mu_{g,k}$ and write
$U_{\alpha }=(u^1_{\alpha},...,u^k_{\alpha})$ for all $\alpha $. It
is easily to see that $(U_{\alpha })_{\alpha }$ is bounded in
$H^{1,2}_k(M)$. Thus, up to a subsequence, $U_{\alpha
}\rightharpoonup U$ in $H^{1,2}_k(M)$, $U_{\alpha }\rightarrow U$ in
$L^2_k(M)$, $U_{\alpha }\rightharpoonup U$ in $L_k^{2^*}(M)$, and
$U_{\alpha }\rightarrow U$ almost everywhere in $M$. By the weakly
convergence in $H^{1,2}_k(M)$, one obtain that

\begin{equation}\label{eq1}
\int _M|\nabla U_{\alpha }|^2dv_g=\int _M|\nabla (U_{\alpha } -
U)|^2dv_g +\int _M|\nabla U|^2dv_g +o(1)
\end{equation}

\n where $o(1)\rightarrow 0$ as $\alpha \rightarrow \infty $. By the
Proposition \ref{B-L0}, it follows that

\begin{equation}\label{eq2}
\int _MF(U_{\alpha })dv_g=\int _M F(U_{\alpha } - U)dv_g + \int _M
F(U)dv_g +o(1)\;
\end{equation}

\n where $o (1)\rightarrow 0$ as $\alpha \rightarrow \infty $. By the
Theorem \ref{Teo.2}, we obtain

\begin{equation}\label{eq3}
\left( \int _MF(U_{\alpha } - U)dv_g\right) ^{2/2^*}\leq
\mathcal{A}_0(n,F) \int _M |\nabla (U_{\alpha }-U)|^2dv_g
\end{equation}

\[
+M_G\mathcal{B}_0(n,F,G,g) \int _M |U_{\alpha } - U|^2dv_g \,.
\]

\n By (\ref{eq1}), (\ref{eq2}) and (\ref{eq3}), we get

\begin{eqnarray*}
\left( 1-\int _MF(U)dv_g \right)^{2/2^*}&=&
\left( \int _MF(U_{\alpha } - U)dv_g \right)^{2/2^*} + o (1)\\
&\leq &  \mathcal{A}_0(n,F) \int _M |\nabla (U_{\alpha }-U)|^2dv_g\\
&&+ M_G\mathcal{B} _0(n,F,G,g) \int _M |U_{\alpha } - U|^2dv_g +o (1)\\
&=& \mathcal{A}_0(n,F) \left( \int _M |\nabla U_{\alpha }|^2dv_g -
\int _M |\nabla U|^2dv_g \right)\\
&&+ M_G\mathcal{B} _0(n,F,G,g) \int _M |U_{\alpha } - U|^2dv_g +o(1)\\
&=& \mathcal{A}_0(n,F) \left( \int _M |\nabla U_{\alpha }|^2dv_g -
\int_M |\nabla U|^2dv_g \right)\\
&& + o(1)
\end{eqnarray*}

\n where  $o(1)\rightarrow 0$ as $\alpha \rightarrow \infty $. We also
have

\[
I(U_{\alpha })=\mu_{g,k} +o(1)\;.
\]

\n See that

\[
\int_M G(x,U_{\alpha })dv_g=\int _MG(x,U)dv_g +o (1)\;.
\]

\n Thus,

\[
\int _M |\nabla U_{\alpha }|^2dv_g  -  \int_M |\nabla
U|^2dv_g = I(U_{\alpha }) - I(U) + o(1)\;.
\]

\n Hence,

\begin{eqnarray*}
\mathcal{A}_0(n,F) \left( \int _M |\nabla U_{\alpha }|^2dv_g  - \int
_M |\nabla U|^2dv_g \right)=
\end{eqnarray*}

\[
\mathcal{A}_0(n,F) \left(  \mu_{g,k} -I(U) \right) +o(1)\;,
\]

\n from which we derive

\begin{equation*}
\left( 1-\int _MF(U)dv_g \right)^{2/2^*} \leq
\mathcal{A}_0(n,F) \left(  \mu_{g,k} -I(U) \right) +o (1)\;.
\end{equation*}

\n Note that

\[
I(U)\geq \mu_{g,k} \left( \int _MF(U)dv_g
\right)^{2/2^*}\;.
\]

\n Thus,

\begin{eqnarray*}
\left( 1-\int _MF(U)dv_g \right)^{2/2^*} &\leq &
\mathcal{A}_0(n,F)
\left(  \mu_{g,k} -\mu_{g,k} \left( \int _MF(U)dv_g \right)^{2/2^*} \right) + o(1)\\
&=& \mathcal{A}_0(n,F)  \mu_{g,k}\left( 1 - \left( \int _MF(U)dv_g
\right)^{2/2^*} \right) +o(1).
\end{eqnarray*}

\n By the Fatou's Lemma,

\[
\int _MF(U)dv_g\leq \liminf_{\alpha \rightarrow + \infty} \int _MF(U_{\alpha })dv_g=1 \;.
\]

\n This implies that

\[
1 - \left( \int _MF(U)dv_g\right) ^{2/2^*}\geq 0 \;.
\]

\n It follows from $1 - \int _MF(U)dv_g\geq 0$ and $2/2^*<1$
that

\begin{eqnarray*}
1&=&\left(1 -  \int _MF(U)dv_g + \int _MF(U)dv_g \right)^{2/2^*}\\
 &\leq &\left(1 -  \int _MF(U)dv_g +  \right)^{2/2^*} + \left( \int _MF(U)dv_g \right)^{2/2^*}.
\end{eqnarray*}

\n Hence,

\[
\int _MF(U)dv_g=1\;,
\]

\n since $\mu_{g,k}< \mathcal{A}_0(n,F)^{-1}$. Another hand, one has

\begin{eqnarray*}
\int _M |\nabla U_{\alpha }|^2dv_g  -  \int _M |\nabla U|^2dv_g  &=&I(U_{\alpha }) - I(U) +o (1)\\
&=& \mu_{g,k} -I(U)
 + o(1)\\
&\leq & \mu_{g,k} -\mu_{g,k} \left( \int _MF(U)dv_g
\right)^{2/2^*} +o(1)\\
&=&o(1)\; .
\end{eqnarray*}

\n From this last inequality and from (\ref{eq1}), one obtains

\[
\int _M |\nabla (U_{\alpha } - U)|^2dv_g =o (1)\;.
\]

\n Thus, $U_{\alpha }\rightarrow U$ in $H^{1,2}_k(M)$ and  $I(U)=\mu_{g,k}$. \bl

\subsection{Regularity of solutions}

In what follows we discuss the regularity of the weak solutions $U =
(u_1, \ldots, u_k) \in H^{1,k}_k(M)$ of the system

\begin{gather} \label{syss}
-\Delta_g u_i + \frac{1}{2} \frac{\partial G(x,U)}{\partial t_i} =
\frac{1}{2^*} \frac{\partial F(U)}{\partial t_i}\ \ \mbox{on}\ M\,,\ i=1,\cdots,k,
\tag{$S_1$}
\end{gather}

\n where the potential functions $F$ and $G$ are of class $C^1$ and belong to ${\cal F}_k$ and ${\cal G}_k$.

\begin{propo}\label{Reg}
Let $(M,g)$ be a compact Riemannian manifold of dimension $n \geq
3$. If $U = (u_1,...,u_k) \in H^{1,2}_k(M)$ is a weak solution of
the system (\ref{syss}), then $U \in C_k^1(M)$.
\end{propo}

\n {\bf Proof of Proposition \ref{Reg}.}: Given real numbers $l>0$ and $\beta>1$, we consider
the functions $\tau$ and $\sigma$ given by

\[
\tau(s) = \left\{
\begin{array}[c]{ll}%
s^\beta & {\rm if}\ 0 \leq s \leq l\\
\beta l^{\beta - 1} (s - l) + l^{\beta} & {\rm if} \ s > l
\end{array}
\right.,
\]

\[
\sigma(s) = \left\{
\begin{array}[c]{ll}%
s^{2(\beta - 1) + 1} & {\rm if}\ 0 \leq s \leq l,\\
\left( 2(\beta - 1) + 1 \right) l^{2(\beta - 1)} (s - l) + l^{2(\beta
- 1) + 1} & {\rm if} \ s > l\ .
\end{array}
\right.
\]

\n Since $u_i\in H^{1,2}_k(M)$, we have that the positive and negative parts of $u_i$, $u_i^+$ and $u_i^-$, are also in the
Sobolev space $u_i\in H^{1,2}_k(M)$. In what follows, we assume that $u_i\geq0$, in other case we can
use the same argument for the positive and negative parts of $u_i$.
Also, it is convenient to choose $\beta>1$ such that $2\beta  \leq
2^*$. Let $U = (u_1, \ldots, u_k)$. Since $\tau$ and $\sigma$ are
Lipschitz functions, taking $\psi = (\sigma(u_1), \ldots,
\sigma(u_k)) \in H^{1,2}_k(M)$ as a test function for the system (\ref{syss}), we get

\[
\sum_{i=1}^k \int_{M} \nabla u_i \cdot \nabla
\sigma(u_i) \; dx = \sum_{i=1}^k \int_{M} h_i(U) \sigma(u_i)\; dv_g\; ,
\]

\n where $h(U) = \frac{1}{2^*} \frac{\partial F(U)}{\partial t_i}-\frac{1}{2} \frac{\partial G(x,U)}{\partial t_i}$. Using now that there exists a constant $c_0 > 0$, independent of $l$, such that

\[
c_0 \tau'(s)^2 \leq \sigma'(s)
\]

\n for all $s\in \R$ and

\[
|h(t)| \leq c_0 (|t|^{2^* - 1} + 1)
\]

\n for all $t \in \mathbb{R}^k$, we obtain a constant $c_1>0$,
independent of $l$, such that

\[
\int_{M} |\nabla \left( \tau(u_i) \right) |^2 \; dv_g =
\int_{M} |\nabla u_i|^2 \tau'(u_i)^2\; dv_g \leq c_1
\int_{M} |u|^{2^* - 1} \sigma(u_i)\; dv_g + c_1
\]

\[
\leq c_1 \int_{\Omega} |U|^{2^* - 1} \sigma(|U|)\; dv_g + c_1\ .
\]

\n Note that there is a constant $c>0$, independent of $l$, such that

\[
\sigma(s) \leq c \tau(s)^2
\]

\n for all $s > 0$. Thus, we obtain

\[
\int_{M} |\nabla \left( \tau(u_i) \right) |^2 \; dv_g \leq
cc_1 \int_{M} |U|^{2^* - 2} \tau(|U|)^2\; dv_g + c_1\ .
\]

\n Hence, applying the classical Sobolev inequality, we find a constant
$c_2> 0$, independent of $l$, such that

\[
\left( \int_{M} \tau(u_i)^{2^*} \; dv_g \right)^{2/2^*}
\leq c_2\int _{M}|\nabla (\tau(u_i))|^2dv_g
\]

\[
\leq c_2 \int_{M} |U|^{2^* - 2} \tau(|U|)^2\; dv_g + c_2\,.
\]

\n Moreover, it follows easily from the definition of $\tau$ that there
is a constant $c_3>0$, independent of $l$, such that

\[
\tau(|U|) \leq c_3 \sum_{i=1}^k \tau(u_i)\ .
\]

\n Hence, using H\"{o}lder and Young inequalities, we find a constant
$c_4>0$, independent of $l$, such that

\[
\left( \int_{\Omega} \tau(|U|)^{2^*} \; dv_g
\right)^{2/2^*} \leq c_4 \int_{\Omega} \tau(|U|)^{2} \;
dv_g + c_4\ .
\]

\n Taking then the limit $l\rightarrow + \infty$ in the inequality
above, we find

\[
\left( \int_{M} |U|^{\beta 2^*} \; dv_g \right)^{2/2^*}
\leq c_4 \int_{M} |U|^{\beta 2} \; dv_g + c_4\ .
\]

\n Therefore, we get $|U|\in L^p(M)$ with $p=\beta 2^*>2^*$, and the
conclusion follows from the classical regularity theory of elliptic PDEs applied to each equation of (\ref{syss}).  \bl \\

 \end{document}